\documentclass{tac}

\usepackage[all,2cell]{xy}
\UseAllTwocells
\usepackage{amsfonts}
\usepackage{amsmath}
\usepackage{amssymb}
\usepackage{lscape} 
\usepackage{multicol}
\usepackage{enumerate}
\usepackage{rotating}
\usepackage{flafter}

\usepackage[colorlinks=true]{hyperref}
\hypersetup{allcolors=[rgb]{0.1,0.1,0.4}}

\author{Rose Kudzman-Blais}

\thanks{The author would like to thank Richard Blute for suggesting this topic as part of their doctoral thesis and for his support during the development of the theory, as well as Shayesteh Naimabadi for many instructive discussions on the topic. The author also wishes to thank Jean-Simon Pacaud Lemay, whose help was instrumental in developing the Examples and Properties Section of the paper, as well as the entire category theory group at Macquarie University for the opportunity to visit and present on the topic. The author also acknowledges the financial support of the Natural Sciences and Engineering Research Council of Canada (NSERC), through a grant awarded to Richard Blute, as well as the Japan Society for the Promotion of Science (JSPS), through a JSPS Postdoctoral Fellowship for Research in Japan (ID P25730).}

\address{Department of Mathematics and Statistics, University of Ottawa\\
150 Louis-Pasteur Pvt, Ottawa, ON, Canada K1N 6N5 \\[5pt]
Research Institute for Mathematical Sciences, Kyoto University\\
Kitashirakawa Oiwake-cho, Sakyo-ku, Kyoto 606-8502 Japan\\
}

\title{Linearly Distributive Fox Theorem}

\copyrightyear{2026}


\eaddress{rosekblais@gmail.com}

\newenvironment{fullwidth}{%
  \everydisplay=\expandafter{%
    \the\everydisplay
    \displayindent=0pt
    \displaywidth=\columnwidth
  }%
  \ignorespaces
}{\ignorespacesafterend}

\newcommand{\rarr}{\rightarrow}
\newcommand{\Rarr}{\Rightarrow}
\newcommand{\op}{\oplus}
\newcommand{\ot}{\otimes}
\newcommand{\os}{\oslash}

\renewcommand{\c}{\,\colon}

\newcommand{\nat}{\mathrm{\bf nat}}

\newcommand{\exch}{\mathrm{\bf exch}}

\newcommand{\mix}{\mathrm{mix}}

\mathbfdef[bzero]{0}
\mathbfdef[bone]{1}

\mathcaldef[cX]{X}

\mathcaldef[cY]{Y}
\mathcaldef[cZ]{Z}

\newcommand{\bP}{\mbox{${\mathbb P}$}}
\newcommand{\bX}{\mbox{${\mathbb X}$}}
\newcommand{\bbX}{\mbox{$\scriptscriptstyle{\mathbb X}$}}
\newcommand{\bbY}{\mbox{$\scriptscriptstyle{\mathbb Y}$}}
\newcommand{\bY}{\mbox{${\mathbb Y}$}}
\newcommand{\bZ}{\mbox{${\mathbb Z}$}}

\newcommand{\PCOH}{\mbox{$\bP$-${\sf Coh}$}}

\mathrmdef[B]{B}
\mathrmdef[inc]{inc}

\newcommand{\LDC}{\mbox{\small{$\mathbf{LDC}$}}}
\newcommand{\SMLDC}{\mbox{\small{$\mathbf{SMLDC}$}}}
\newcommand{\CLDC}{\mbox{\small{$\mathbf{CLDC}$}}}
\newcommand{\SLDC}{\mbox{\small{$\mathbf{SLDC}$}}}
\newcommand{\FLDC}{\mbox{\small{$\mathbf{FLDC}$}}}
\newcommand{\FCLDC}{\mbox{\small{$\mathbf{FCLDC}$}}}
\newcommand{\DUO}{\mbox{\small{$\mathbf{DUO}$}}}
\newcommand{\CART}{\mbox{\small{$\mathbf{CART}$}}}
\newcommand{\SMON}{\mbox{\small{$\mathbf{SMON}$}}}

\newcommand{\MLL}{$\mathrm{MLL}$}
\newcommand{\MLLp}{$\mathrm{MLL}^+$}

\begin{document}

\maketitle
\begin{abstract}
Linearly distributive categories (LDC), introduced by Cockett and Seely to model multiplicative linear logic, are categories equipped with two monoidal structures that interact via linear distributivity. A seminal result in monoidal category theory is the Fox theorem, which characterizes cartesian categories as symmetric monoidal categories whose objects are equipped with canonical comonoid structures. The aim of this work is to extend the Fox theorem to LDCs and characterize the subclass of cartesian linearly distributive categories (CLDC). To do so, we introduce medial linearly distributive categories (MLDC), medial linear functors, and medial linear transformations. The former are LDCs which respect the logical medial rule, appearing frequently in deep inference, and are the appropriate categorical structure at the intersection of LDCs and duoidal categories. 
\end{abstract}

\bigskip
\contentsline {section}{\numberline {1}Introduction}{1}{section.1}%
\contentsline {section}{\numberline {2}Fox Theorem}{5}{section.2}%
\contentsline {subsection}{\numberline {}Monoidal Categories}{5}{subsection.2.1}%
\contentsline {subsection}{\numberline {}Comonoids}{8}{subsection.2.3}%
\contentsline {section}{\numberline {3}Linearly Distributive Categories}{10}{section.3}%
\contentsline {subsection}{\numberline {}Definition, Mix and Examples}{11}{subsection.3.1}%
\contentsline {subsection}{\numberline {}Linear Functors and Transformations}{16}{subsection.3.18}%
\contentsline {subsection}{\numberline {}Cartesian Linearly Distributive Categories}{20}{subsection.3.28}%
\contentsline {section}{\numberline {4}Duoidal Categories}{21}{section.4}%
\contentsline {subsection}{\numberline {}Definition, Symmetry and Examples}{21}{subsection.4.1}%
\contentsline {subsection}{\numberline {}Duoidal Functors and Transformations}{23}{subsection.4.6}%
\contentsline {subsection}{\numberline {}Symmetry}{26}{subsection.4.13}%
\contentsline {subsection}{\numberline {}Cartesian and Cocartesian Duoidal Categories}{28}{subsection.4.19}%
\contentsline {section}{\numberline {5}Medial Linearly Distributive Categories}{30}{section.5}%
\contentsline {subsection}{\numberline {}Medial Rule and Main Definition}{32}{subsection.5.3}%
\contentsline {subsection}{\numberline {}Symmetry}{35}{subsection.5.5}%
\contentsline {subsection}{\numberline {}Adding Negation}{39}{subsection.5.10}%
\contentsline {subsection}{\numberline {}Examples}{40}{subsection.5.14}%
\contentsline {subsection}{\numberline {}Properties}{42}{subsection.5.20}%
\contentsline {section}{\numberline {6}Medial Linear Functors and Transformations}{45}{section.6}%
\contentsline {subsection}{\numberline {}2-Category \mbox {\small {$\mathbf {SMLDC}$}}}{49}{subsection.6.5}%
\contentsline {subsection}{\numberline {}Frobenius Medial Linear Functors and Transformations}{53}{subsection.6.8}%
\contentsline {section}{\numberline {7}Linearly Distributive Fox Theorem}{56}{section.7}%
\contentsline {subsection}{\numberline {}Medial Bimonoids}{56}{subsection.7.1}%
\contentsline {subsection}{\numberline {}The Right Adjoint 2-Functor {\ensuremath {\rm B}}[-]}{61}{subsection.7.13}%
\contentsline {section}{Appendices}{67}{section*.40}%
\contentsline {subsection}{\numberline {A}Large Commuting Diagrams}{67}{subsection.Alph0.1}%

\section{Introduction}\label{sec:Introduction}

Within monoidal category theory, \emph{cartesian categories} refer to categories whose monoidal product is the categorical product and whose monoidal unit is the terminal object. Perhaps the most well-known result concerning cartesian categories is the {\it Fox theorem}. In his seminal 1976 paper ``Coalgebras and cartesian categories'', Fox demonstrated that the inclusion functor from cartesian categories into symmetric monoidal categories has a right adjoint \cite{Fox_1976}. This right adjoint functor maps a symmetric monoidal category $(\cX, \os, I)$ to its category of cocommutative comonoids and their morphisms, meaning it takes objects $A$ in \cX\ equipped with coherent diagonal and counit maps:
\[\Delta_A\c A\rarr A\os A \qquad\qquad  e_A\c A\rarr I\]
As a corollary of this theorem, we see that cartesian categories are precisely the symmetric monoidal categories whose objects have a canonical comonoid structure and whose morphisms preserve these structures. 

Cartesian categories have played a key role in the field of categorical logic. Indeed, one of the greatest successes of categorical logic was Lambek’s equivalence between intuitionistic propositional calculus and cartesian closed categories \cite{Lambek_1972}. If we extend this paradigm to capture the semantics of classical logic and consider cartesian closed categories with a notion of involutive negation, we get a collapse to posetal categories \cite{Lambek_Scott_1988}. This result is now often known as {\it Joyal’s paradox}. As such, there has not yet been a definitive framework for the categorical semantics of classical logic.

The field of categorical logic nonetheless continued to thrive, in particular with the development of linear logic by Girard \cite{Girard_1987}. The categorical semantics for linear logic were explored by Seely, who identified that the appropriate framework was given by {\it $*$-autonomous categories} \cite{Seely_1989}, previously introduced by Barr \cite{Barr_1979}. {\it Linearly distributive categories} (LDC) were later introduced by Cockett and Seely as alternative categorical semantics for linear logic, which take the notions of multiplicative conjunction and disjunction, and the linear distributivity between them, as primitive \cite{Cockett_Seely_1997_LDC}. Specifically, LDCs are categories $\bX$ with two monoidal structures $(\bX, \ot, \top)$ and $(\bX, \op, \bot)$, the former known as the tensor structure and the latter as the par structure, equipped with linear distributivity natural transformations:
\[\delta^L_{A,B,C}\c A\ot (B\op C)\rarr (A\ot B)\op C \quad\qquad \delta^R_{A,B,C}\c (A\op B)\ot C\rarr A\op (B\ot C) \]

The theory of categorical linear logic has proved useful in studying categorical semantics for classical logic. LDCs provide the minimal framework necessary to discuss a logic with two connectives and constants, with left and right (multiplicative-style) introduction rules, and cut as the only structural rule. Since most reasonable models of classical logic include such structure, categorical models for classical logic can be built upon the foundations of LDCs (or $*$-autonomous categories, if negation is instead taken as primitive). 

As the structural rules of contraction and weakening are central to classical logic and are traditionally modeled by cartesian monoidal structures, a recent paradigm of categorical classical logic has been to add ``cartesian'' structure to LDCs (or $*$-autonomous categories). In this vein, F\"{u}hrmann  and Pym developed the notions of classical categories and Dummett categories: poset-enriched $*$-autonomous categories and LDCs respectively whose objects have coherent notions of $\ot$-diagonals and counits, and $\op$-multiplications and units \cite{Fuhrmann_Pym_2007}. Further, Stra{\ss}burger and Lamarche both separately considered $*$-autonomous categories whose objects also have coherent $\ot$-diagonals and counits, and $\op$-multiplications and units, but whose maps do not necessarily strictly preserve these structures \cite{Strassburger_2007_1, Lamarche_2007}. 

More recently, Naeimabadi developed the theory of cartesian linear bicategories \cite{Naeimabadi_2024} and recent work has investigated De Morgan and Peirce's calculus of relations, providing an axiomatisation based on the concepts of cartesian and linear bicategories \cite{Bonchi_DiGiorgio_Haydon_Sobocinski_2024, Bonchi_DiGiorgio_Trotta_2024}. Indeed, the central examples of linear bicategories, the bicategorical analogues of LDCs introduced by Cockett, Koslowski, and Seely \cite{Cockett_Koslowski_Seely_2000}, are also cartesian bicategories. Recall that the Fox theorem is, in fact, central to the theory of cartesian bicategories, as introduced by Carboni and Walters \cite{Carboni_Walters_1987}. 

It was within this context that the idea arose to investigate \emph{cartesian linearly distributive categories} (CLDC) and to develop a linearly distributive version of Fox’s theorem. Precisely, CLDCs are LDCs whose tensor structure is cartesian and par structure is cocartesian. They were first introduced alongside LDCs by Cockett and Seely, though the notion received little attention beyond the first paper \cite{Cockett_Seely_1997_LDC}. 

In developing a linearly distributive analogue of the Fox theorem, it quickly became apparent that no such construction could be adapted to all LDCs. For the same reason that the Fox theorem is an adjunction between cartesian categories and {\it symmetric} monoidal categories, we need to restrict to LDCs equipped with additional structure maps:
\[\Delta_{\bot}\c \bot\rarr \bot\ot\bot \qquad \nabla_{\top}\c \top\op\top\rarr \top \qquad m\c \bot\rarr \top\]
\[\mu_{A,B,C,D}\c (A\ot B)\op (C\ot D)\rarr (A\op C)\ot (B\op D) \]
Viewed through the lens of logical entailment, the latter map is known as the \emph{medial rule}, first introduced by Br\"unnler and Tiu for a system of classical logic within the deep inference framework \cite{Brunnler_Tiu_2001}. Medial maps were also key in both Stra{\ss}burger’s and Lamarche’s categorical models for classical logic. 

Within monoidal category theory, these medial maps are better known as instances of the interchange maps of \emph{duoidal categories}. Duoidal categories are categories $\cX$ with two monoidal structures $(\cX, \diamond, I)$ and $(\cX, \star, J)$, equipped with precisely the structure maps introduced above \cite{Aguiar_Mahajan_2010}. The theory of duoidal categories first appeared in the work of Joyal and Street while studying braided monoidal categories \cite{Joyal_Street_1993}. 

With this perspective in mind, we introduce the notion of \emph{medial linearly distributive categories} (MLDC), inspired both by the work of Stra{\ss}burger and Lamarche in categorical classical logic and by the theory of duoidal categories. MLDCs can be regarded both as potential categorical semantics for (negation-free) multiplicative linear logic with the medial rule, and as an appropriate intersection point between LDCs and duoidal categories. 

Leveraging the well-developed theory of LDCs and duoidal categories, we define the appropriate functors and transformations between MLDCs. Within this framework, we consider an analogue to comonoids, \emph{medial bimonoids}: objects $A$ equipped with coherent $\ot$-diagonal, $\ot$-counit, $\op$-multiplication and $\op$-unit maps,
\[\Delta_A\c A\rarr A\ot A \qquad e_A\c A\rarr \top \qquad \nabla_A\c A\op A\rarr A \qquad u_A\c \bot\rarr A \]
Considering the bicommutative medial bimonoids and their morphisms in a symmetric MLDC yields the construction needed to define a right adjoint to the inclusion 2-functor from CLDCs to symmetric MLDCs. This yields the desired {\it linearly distributive Fox theorem} and, as a corollary, a characterization of CLDCs as symmetric MLDCs whose objects have coherent medial bimonoid structure. 

In the process of developing this paper, it became evident how poorly understood CLDCs are. This prompted the author to collaborate with Pacaud Lemay to explore the nature of CLDCs beyond their characterization via the linearly distributive Fox theorem. This resulted in a separate paper, ``Cartesian linearly distributive categories: Revisited'' \cite{Kudzman-Blais_Lemay_2025}. While both papers are written to be independently comprehensible, we recommend reading them in tandem to gain a complete perspective on cartesian and cocartesian structure in the context of LDCs. \\

{\bf Outline.} Given that the results developed in this paper are based on the theories of the Fox theorem, LDCs, and duoidal categories, the paper begins with three substantial preliminary sections. Section~\ref{sec:Fox_thm} presents the basics of monoidal categories, cartesian categories, and the Fox theorem. Section~\ref{sec:LDC} reviews the literature on LDCs, including linear functors and linear transformations. Section~\ref{sec:duoidal_categories} introduces the necessary background on duoidal categories. Then, Section~\ref{sec:MLDC} begins with a motivation for the introduction of the medial rule and subsequently defines MLDCs. The theory of MLDCs is developed: symmetric MLDCs are defined, adding negation is discussed, examples are introduced, and properties are presented. Section~\ref{sec:medial_linear_functor} defines the appropriate functors and transformations between symmetric MLDCs. Section~\ref{sec:LD_Fox} introduces medial bimonoids and shows that considering such structures and their maps does, in fact, yield a CLDC. Finally,  the main result of this paper is proved: the linearly distributive Fox theorem. Appendix \ref{app:commuting_diagrams} contains commuting diagrams which were too large to include in the body of the text. \\

{\bf Conventions.} There are a few notational choices made in this paper and we take a moment to detail the most important.
\begin{itemize}
    \item Composition of maps in a category is denoted by $;$ and is in diagrammatic order. Objects are denoted by capital letters, while the maps are lowercase letters. Identity morphisms are denoted by $1_A$.
    \item Monoidal products and monoidal units will be denoted with various different symbols depending on the context:
    \begin{itemize}
    \item $\os$ and $I$ are used for standard monoidal categories \cX, 
    \item $\ot$ and $\top$ are used for the tensor monoidal structure, and $\op$ and $\bot$ for the par in a LDC \bX, in agreement with the notation introduced by Cockett and Seely in \cite{Cockett_Seely_1997_LDC}, as opposed to Girard's notation \cite{Girard_1987},
    \item $\diamond$ and $I$ are used for the first monoidal structure and $\star$ and $J$ for the second in a duoidal category \cX, in agreement with the notation introduced by Aguiar and Mahajan \cite{Aguiar_Mahajan_2010}, and 
    \item $\times$ and $\bone$ denote the binary categorical product and the terminal object respectively, while $+$ and $\bzero$ denote the binary categorical coproduct and the initial object in a category.
    \end{itemize}
\end{itemize}

\section{Fox Theorem}\label{sec:Fox_thm}

We begin by a description of the Fox theorem for cartesian categories, as the traditional result will be a guiding thread to generalizing it to the context of linearly distributive categories. We invite any reader who is well versed in the Fox theorem to skip this section. \\

\subsection{Monoidal Categories}  \hfill\

\vspace{0.5\baselineskip}
We briefly outline the foundational definitions of monoidal category theory, primarily to establish notation, a standard reference for which is MacLane \cite[Sec III, VII, XI]{MacLane_1998}.  

Recall that a {\bf monoidal category} $(\cX, \os, I)$ is a category \cX\ equipped with a {\bf product} functor $\os\c\cX\times\cX\rarr\cX$, a distinguished {\bf unit} object $I$, and {\bf associator}, {\bf right unitor} and {\bf left unitor} natural isomorphisms
\[ \alpha_{X,Y,Z}\c (X\os Y)\os Z\rarr X\os (Y\os Z)  \qquad \rho_{X}\c X\rarr X\os I\qquad \lambda_{X}\c X\rarr I\os X\]
satisfying the {\it associativity pentagon diagram} \mbox{\bf\footnotesize{(MC.1)}}\label{cc:monoidal_cat_asso_pentagon} and the {\it unit triangle identities} \mbox{\bf\footnotesize{(MC.2)}}\label{cc:monoidal_cat_unit_triangle}. A monoidal category is {\bf braided} if it has a {\bf braiding} natural isomorphism 
\[ \sigma_{X, Y}\c X\os Y\rarr Y\os X \]
which satisfies the {\it hexagon identities} \mbox{\bf\footnotesize{(BMC.1)}}\label{cc:braided_monoidal_cat_asso_hexagon}. Note that it then follows that the braiding is compatible with the right and left unitors \mbox{\bf\footnotesize{(BMC.2)}}. Further, a monoidal category is {\bf symmetric}, otherwise known as a SMC, if the braiding is {\it self-inverse} \mbox{\bf\footnotesize{(SMC)}}\label{cc:braiding_self_inverse}.

Given a braided monoidal category, there is a {\bf canonical flip} natural isomorphism given by the following equivalent composites
\begin{equation}\begin{gathered}\label{eqn:canonical_flip}\tag{\bf\footnotesize{CF}}
\xymatrixrowsep{1.75pc}\xymatrixcolsep{3.25pc}\xymatrix@L=0.5pc{ 
(W\os X)\os (Y\os Z)\ar[d]_{\alpha_{W,X,Y\os Z}} \\
W\os (X\os (Y\os Z))\ar[r]^{1_{W}\os\alpha_{X,Y,Z}^{-1}} \ar[d]_{1_{W}\os\sigma_{X, Y\os Z}} & W\os ((X\os Y)\os Z) \ar[r]^{1_{W}\os(\sigma_{X, Y}\os 1_{Z})} & W \os ((Y\os X)\os Z)\ar[d]^{1_{W}\os\alpha_{ Y,X,Z}} \\
W\os ((Y\os Z)\os X) \ar[r]_{1_W \os \alpha_{Y,Z,X}} & W\os (Y\os (Z\os X))\ar[r]_{1_W\os (1_Y \os \sigma_{X,Z})}& W\os (Y\os (X\os Z))\ar[d]^{\alpha^{-1}_{W, Y, X\os Z}}\\
&& (W\os Y)\os (X\os Z)
}
\end{gathered}\end{equation}
This natural isomorphism is key to the Fox theorem and, therefore, to improve readability, we will denote it by \[ \tau_{W,X,Y,Z}\c (W\os X)\os (Y\os Z) \Rarr (W\os Y)\os (X\os Z) \]

Given monoidal categories $(\cX, \os, I)$ and $(\cY, \os, I)$, recall that a {\bf lax monoidal functor} $(F, m_{I},m_\os)\c(\cX,\os,I)\rarr(\cY,\os, I)$ is a functor $F\c\cX\rarr\cY$ equipped with the following map and natural transformation 
\[ m_{I}\c I\rarr F(I)\qquad {m_\os}_{A, B}\c F(A)\os F(B)\rarr F(A\os B)\]
satisfying {\it associativity} \mbox{\bf\footnotesize{(MF.1)}}\label{cc:mon_functor_associativity} and {\it unitality} \mbox{\bf\footnotesize{(MF.2)}}\label{cc:mon_functor_unit} coherence conditions. A {\bf colax monoidal functor} $(F, n_{I}, n_\os)\c(\cX,\os,I)\rarr(\cY,\os, I)$ is a lax monoidal functor $(\cX^{op},\os,I)\rarr(\cY^{op},\os, I)$. Moreover, a monoidal functor is a (co)lax monoidal functor whose the structure maps are isomorphisms. Further, if the monoidal categories are symmetric, then the (co)lax monoidal functor is {\bf symmetric} if the structure maps {\it interact coherently with the braidings} \mbox{\bf\footnotesize{(SMF)}}\label{cc:mon_functor_braiding}.

\begin{remark} Given a SMC $(\cX, \os, I)$ and a symmetric lax monoidal functor \\\mbox{$(F, m_{I},m_\os)\c(\cX,\os,I)\rarr(\cY,\os, I)$}, we note that the following commuting diagrams involving the canonical flip \mbox{\bf\footnotesize{(CF)}}.
\begin{equation}\begin{gathered}\label{cc:canonical_flip_unitor}
\xymatrixrowsep{1.75pc}\xymatrixcolsep{3.25pc}\xymatrix{ 
A\os B\ar[r]^-{\rho_A\os \rho_B}\ar[d]_-{\rho_{A\os B}} & (A\os I)\os (B\os I)\ar[d]^-{\tau_{A,I,B,I}} \\
(A\os B)\os I \ar[r]_-{{1_{A\os B}\os \rho_{I}}} & (A\os B)\os (I\os I)
}
\quad
\xymatrixrowsep{1.75pc}\xymatrixcolsep{3.25pc}\xymatrix{ 
A\os B\ar[r]^-{\lambda_A\os \lambda_B}\ar[d]_-{\lambda_{A\os B}} & (I\os A)\os (I\os B)\ar[d]^-{\tau_{I,A,I,B}} \\
I\os (A\os B) \ar[r]_-{{\lambda_{I}\os 1_{A\os B}}} & (I\os I)\os (A\os B)
}
\end{gathered}\end{equation}
\begin{equation}\begin{gathered}\label{cc:canonical_flip_braiding}
\xymatrixrowsep{1.75pc}\xymatrixcolsep{3.75pc}\xymatrix{ 
(A\os B)\os (C\os D)\ar[r]^-{\sigma_{A,B}\os \sigma_{C,D}}\ar[d]_-{\tau_{A,B,C,D}} & (B\os A)\os (D\os C)\ar[d]^-{\tau_{B, A, D, C}} \\
(A\os C)\os (B\os D)\ar[r]_-{\sigma_{A\os C, B\os D}} & (B\os D)\os (A\os C)
}\end{gathered}\end{equation}
\begin{equation}\begin{gathered}\label{cc:canonical_flip_associator}
\xymatrixrowsep{1.75pc}\xymatrixcolsep{5pc}\xymatrix{ 
((A\os B)\os C)\os ((D\os E)\os F)\ar[r]^-{\alpha_{A,B,C}\os \alpha_{D, E,F}}\ar[d]_-{\tau_{A\os B, C, D\os E, F}} & (A\os (B\os C))\os (D\os (E\os F))\ar[d]^-{\tau_{A,B\os C, D, E\os F}} \\
((A\os B)\os (D\os E)) \os (C\os F)\ar[d]_-{\tau_{A,B,D,E}\os 1_{C\os F}} & (A\os D) \os ((B\os C)\os (E\os F))\ar[d]^-{1_{A\os D}\os\tau_{B,C,E,F}}\\
((A\os D)\os (B\os E))\os (C\os F)\ar[r]_-{\alpha_{A\os D, B\os E, C\os F}} & (A\os D)\os ((B\os E)\os (C\os F))
}
\end{gathered}\end{equation}
\begin{equation}\begin{gathered}\label{cc:strong_symmetric_monoidal_flip}
\xymatrixrowsep{1.75pc}\xymatrixcolsep{5pc}\xymatrix{ 
(F(W)\os F(X)) \os (F(Y)\os F(Z)) \ar[r]^-{{m_\os}_{W,X}\os {m_\os}_{Y, Z}} \ar[d]_-{\tau_{F(W), F(X), F(Y), F(Z)}} 
& F(W\os X) \os F(Y\os Z) \ar[d]^-{{m_\os}_{W\os X, Y\os Z}} \\
(F(W)\os F(Y)) \os (F(X)\os F(Z)) \ar[d]_-{{m_\os}_{W,Y}\os {m_\os}_{X,Z}} 
& F((W\os X)\os (Y\os Z)) \ar[d]^-{F(\tau_{W, X, Y, Z})} \\
F(W\os Y) \os F(X\os Z) \ar[r]_-{{m_\os}_{W\os Y, X\os Z}} 
& F((W\os Y)\os (X\os Z))
} 
\end{gathered}\end{equation}
\end{remark}

Given lax monoidal functors $(F, m^{F}_{I},m^{F}_\os), (G, m^{G}_{I},m^{G}_\os)\c(\cX,\os, I)\rarr(\cY, \os, I)$, recall that a {\bf monoidal transformation} $\alpha\c(F, m^{F}_{I},m^{F}_\os)\Rarr (G, m^{G}_{I},m^{G}_\os)$ is a natural transformation $\alpha\c F\Rarr G$ whose component maps {\it behave coherently with the monoidal structures} \mbox{\bf\footnotesize{(MT)}}\label{cc:mon_trans}. Then, a monoidal transformation between colax monoidal functors is simply a monoidal transformation when viewing the functors as lax monoidal between the op-categories. 

Altogether, we can define the 2-category \SMON\ of SMCs, symmetric monoidal functors and monoidal transformations. 

Now, the Fox theorem is concerned with characterizing the SMCs which are cartesian. In this context, a {\bf cartesian category} is a category with finite categorical products and a {\bf cocartesian category} is a category with finite categorical coproducts, where binary products, terminal objects, binary coproducts and initial objects, with the corresponding unique maps, are denoted as follows
\begin{equation*}
\xymatrix{
& C\ar[ld]_{f}\ar[dr]^{g}\ar@{-->}[d]_(.6){\langle f,g\rangle} & && A\ar@{-->}[d]_{t_A} && A\ar[r]^-{\iota^0_{A,B}}\ar[dr]_{h} & A+ B\ar@{-->}[d]_(0.4){[h,k]}& B\ar[l]_-{\iota^1_{A,B}}\ar[dl]^{k} && \bzero\ar@{-->}[d]_{b_A}\\
A& A\times B\ar[l]^{\pi^0_{A,B}}\ar[r]_{\pi^1_{A,B}} & B && \bone && & C &&& A
}
\end{equation*}

A cartesian category \cX\ has a canonical symmetric monoidal structure $(\cX, \times, \bone)$ whose monoidal product is the categorical product
\[ (f\c X\rarr X', g\c Y\rarr Y')\mapsto f\times g = \langle \pi^0_{X, Y};f, \pi^1_{X,Y}; g\rangle\c X\times Y\rarr X'\times Y'\]
and whose monoidal unit is the terminal object, with structure isomorphisms:
\begin{align*}
& \alpha_{X,Y,Z} = \langle \pi^0_{X\times Y, Z}; \pi^0_{X,Y}, \pi^1_{X,Y}\times 1_{Z}\rangle & \sigma_{X, Y} = \langle \pi^1_{X,Y}, \pi^0_{X,Y}\rangle \\
& \rho_{X} = \langle 1_{X}, t_X\rangle & \lambda_{X} = \langle t_X, 1_X\rangle
\end{align*}
Similarly, cocartesian categories \cX\ are canonically symmetric monoidal categories $(\cX, +, \bzero)$. 

Then, there is a 2-category \CART\ of cartesian categories, symmetric monoidal functors and monoidal transformations, known as {\bf cartesian functors} and {\bf cartesian transformations} in this context, which is full sub-2-category of \SMON.

\subsection{Comonoids}  \hfill\

\vspace{0.5\baselineskip}
With the relevant 2-categories now appropriately defined, we can introduce the construction which provides the adjoint to the inclusion. This entire section is based on Fox's original article \cite{Fox_1976}, unless otherwise indicated.

\begin{definition} Let $(\cX, \os, I)$ denote a monoidal category. A {\bf comonoid} in \cX\ is a triple $\langle A, \Delta_A, e_A\rangle$ consisting of an object $A$ and two maps, the {\bf diagonal} and the {\bf counit}
\[ \Delta_A\c A\rarr A\os A \qquad\qquad e_A\c A\rarr I \]
satisfying
\begin{multicols}{2}\raggedcolumns
\noindent
\begin{equation}\begin{gathered}\label{cc:comonoid_assoc}\tag{\bf\footnotesize{CM.1}}
\xymatrixrowsep{1.25pc}\xymatrixcolsep{2.5pc}\xymatrix{ 
A \ar[r]^-{\Delta_A} \ar[d]_-{\Delta_A} & A\os A \ar[dd]^-{1_A \os \Delta_A}  \\
A\os A \ar[d]_-{\Delta_A\os 1_A} \\
(A\os A)\os A\ar[r]_-{\alpha_{A,A,A}} & A\os(A\os A)
}
\end{gathered}\end{equation}
\break
\noindent
\begin{equation}\begin{gathered}\label{cc:comonoid_units}\tag{\bf\footnotesize{CM.2}}
\xymatrixrowsep{1.25pc}\xymatrixcolsep{1.75pc}\xymatrix{ 
A\ar@{=}[rdd]\ar[r]^-{\Delta_A}\ar[d]_-{\Delta_A} & A\os A \ar[d]^-{1_A\os e_A} \\
A\os A\ar[d]_-{e_A\os 1_A} & A\os I\ar[d]^-{\rho_A^{-1}}\\
I\os A\ar[r]_-{\lambda_A^{-1}} & A 
}
\end{gathered}\end{equation}
\end{multicols}
\noindent Further, if the monoidal category is symmetric (braided), a comonoid is {\bf cocommutative} ({\bf cobraided}) if 
\begin{equation}\begin{gathered}\label{cc:cocommutative_comonoid}\tag{\bf\footnotesize{CCM}}
\xymatrixrowsep{1.25pc}\xymatrixcolsep{2.75pc}\xymatrix{ 
A\ar[r]^-{\Delta_A}\ar[dr]_-{\Delta_A} & A\os A \ar[d]^-{\sigma_{A,A}}\\
& A\os A
}
\end{gathered}\end{equation}

A {\bf comonoid morphism} $f\c \langle A, \Delta_{A}, e_A\rangle\rarr\langle B, \Delta_{B}, e_B\rangle$ is a map $f\c A\rarr B$ such that 
\begin{equation}\begin{gathered}\label{cc:comonoid_morphism}\tag{\bf\footnotesize{CMM}}
\xymatrixrowsep{1.25pc}\xymatrixcolsep{2.75pc}\xymatrix{ 
A\ar[r]^-{f}\ar[d]_-{\Delta_A} & B\ar[d]^-{\Delta_B} & A\ar[r]^-{f}\ar[dr]_-{e_A} & B\ar[d]^-{e_B}\\
A\os A\ar[r]_-{f\os f} & B\os B & & I
}
\end{gathered}\end{equation}

We denote the category of cocommutative comonoids and comonoid morphisms in \cX\ by $\mathrm{C}[\cX]$.
\end{definition}

\begin{lemma}\label{lem:unique_comonoid_cartesian_cat}
Given a cartesian category $(\cX, \times, \bone)$, every object $A$ in \cX\ has a canonical unique cocommutative \mbox{$\times$-comonoid} structure 
\[ \langle A, \langle 1_A, 1_A\rangle \c A\rarr A\times A, t_A\c A\rarr\bone\rangle \]
and every arrow $f\c A\rarr B$ in \cX\ is a comonoid morphism with respect to these structures.
\end{lemma}

The key to the Fox theorem is that $\mathrm{C}[\cX]$ inherits a monoidal structure from \cX, due to the existence of the canonical flip \mbox{\bf\footnotesize{(CF)}}, which is in fact cartesian in $\mathrm{C}[\cX]$.

\begin{lemma}\label{lem:constructing_comonoids}
Given a pair of cocommutative comonoids $\langle A, \Delta_A, e_A\rangle$ and $\langle A', \Delta_{A'}, e_{ A'}\rangle$, the triple $\langle A\os A', \Delta_{A\os A'}, e_{A\os A'}\rangle$ defined by
\[ \Delta_{A\os A'} = A\os A' \xrightarrow{\Delta_{A} \os \Delta_{A'}} (A\os A)\os (A'\os A') \xrightarrow{\tau_{A,A,A',A'}} (A\os A')\os (A\os A')\]
\[ e_{A\os A'} = A\os A' \xrightarrow{e_A\os e_{A'}} I \os I \xrightarrow{\rho^{-1}_{I}} I \]
is a cocommutative comonoid. Consequently, $\mathrm{C}[\cX]$ is a SMC and, in particular, $\mathrm{C}[\cX]$ is a cartesian category.
\end{lemma}

This extends canonically to symmetric monoidal functors and transformations.

\begin{lemma}
Given a symmetric monoidal functor $F\c\cX\rarr\cY$ and a cocommutative comonoid $\langle A, \Delta_{A}, e_{A}\rangle$, the triple $\langle F(A), \Delta_{F(A)}, e_{F(A)}\rangle$ defined by
\[ \Delta_{F(A)} = F(A) \xrightarrow{F(\Delta_{A})} F(A\os A) \xrightarrow{{m_\os}^{-1}_{A, A}} F(A) \os F(A)\] \[e_{F(A)} = F(A) \xrightarrow{F(e_{A})} F(I) \xrightarrow{m^{-1}_{I}} I \]
is a cocommutative comonoid. As such, $F$ canonically extends to a cartesian functor $\mathrm{C}[F]\c\mathrm{C}[\cX]\rarr\mathrm{C}[\cY]$.
\end{lemma}

With this, we are finally ready state that the the functor $\mathrm{C}[-]\c\SMON\rarr\CART$, mapping a SMC to its category of cocommutative comonoids and comonoid morphisms is right adjoint to the inclusion functor:

\begin{theorem}\label{thm:Foxthm}\emph{({\bf Fox Theorem})}
$\inc\vdash \mathrm{C}[-]\c\CART\rarr\SMON$ is an adjunction.
\end{theorem}

While the original Fox theorem is only an adjunction between categories, it is well-known that it easily extends to include transformations.

As $\inc$ is fully faithful, the unit of the adjunction is an isomorphism and, as a corollary, we get the following characterization of cartesian categories:
\begin{corollary}\label{cor:SMC_cartesian_iff}
Given a SMC $\cX$, there exists a cartesian category $\cY$ such that $\cX\cong \inc(\cY)$ if and only if there exists an isomorphism $\cX\cong \inc(\mathrm{C}[\cX])$.
\end{corollary}
In other words, a SMC is cartesian if and only if it is isomorphic to its category of cocommutative comonoids. This previous corollary can be unwrapped as a statement about the existence of certain natural transformations.

\begin{corollary}\cite[Thm 4.28]{Heunen_Vicary_2019}\label{cor:char_cartesian_cat}
A SMC $(\cX,\os,I)$ is cartesian if and only if there are natural transformations
\[ \Delta_{A}\c A\rarr A\os A\qquad e_{A}\c A\rarr I\]
such that, $\forall A\in\cX, (A, \Delta_{A}, e_{A})$ is a cocommutative comonoid and, $\forall A,B\in\cX$,
\begin{align*}
&e_{A\os B} = (e_{A}\os e_{B}); \rho_{I}^{-1} &e_{I} = 1_{I} \\
&\Delta_{A\os B} = (\Delta_{A}\os\Delta_{B}); \tau_{A, A, B, B} & \Delta_{I} = \rho_{I}
\end{align*}
\end{corollary}

Of course, there is a dual statement for cocartesian categories, by taking the category of commutative monoids and monoid morphisms in a SMC. We quickly outline these definitions as they will be used throughout.

\begin{definition}
Let $(\cX, \os, I)$ denote a monoidal category. A {\bf monoid} in \cX\ is a comonoid in $\cX^{op}$, i.e. a triple $\langle A, \nabla_A, u_A\rangle$ consisting of an object $A$ and two maps, the {\bf multiplication} and the {\bf unit}
\[ \nabla_A\c A\os A\rarr A \qquad\qquad u_A\c I\rarr A \]
such that the structure maps {\it interact coherently with the associator} \mbox{\bf\footnotesize{(M.1)}}\label{cc:monoid_assoc} and {\it with the unitors} \mbox{\bf\footnotesize{(M.2)}}\label{cc:monoid_units}. Further, if the monoidal category is symmetric (braided), a monoid is {\bf commutative} ({\bf braided}) if the {\it multiplication commutes with the braiding} \mbox{\bf\footnotesize{(CM)}}\label{cc:commutative_monoid}.

A {\bf monoid morphism} $f\c \langle A, \nabla_{A}, e_A\rangle\rarr\langle B, \nabla_{B}, u_B\rangle$ is a comonoid morphism in $\cX^{op}$, i.e. a map $f\c A\rarr B$ which interacts coherently with multiplication and units \mbox{\bf\footnotesize{(MM)}}\label{cc:monoid_morphism}.
\end{definition}

\section{Linearly Distributive Categories}\label{sec:LDC}

Girard introduced linear logic in 1987 to capture the constructive nature of intuitionistic logic and the symmetry of classical logic \cite{Girard_1987}. To do so, Girard omitted contraction and weakening from the structural rules, making linear logic sub-structural, and introduced controlled contraction and weakening via sequents explicitly marked by exponential connectives. The absence of general contraction and weakening led to the existence of two versions of conjunction and disjunction: multiplicative and additive. {\bf Multiplicative linear logic} (\MLL) refers to the fragment concerned only with multiplicative conjunction, disjunction, truth and falsity, as well as linear negation and linear implication. It is sometimes the case that MLL is used to refer to a sub-fragment, which does not include linear negation and linear implication. We shall refer to this more minimal fragment as {\bf negation-free multiplicative linear logic} (\MLLp).

The categorical semantics for linear logic were explored by Seely, who identified that the categorical semantics for MLL are given by $*$-autonomous categories \cite{Seely_1989}, previously introduced by Barr \cite{Barr_1979}. In this context, the monoidal product models multiplicative conjunction, while the involutive functor $*$ models linear negation. The additive and exponential connectives can then be modeled by additional structure on the $*$-autonomous categories, namely finite products and linear exponential comonads respectively. 

Linearly distributive categories (LDC) were introduced by Cockett and Seely as an alternative framework for the semantics of linear logic to $*$-autonomous categories. \cite{Cockett_Seely_1997_LDC}. Unlike $*$-autonomous categories, LDCs take multiplicative conjunction and disjunction as the primitive categorical notions. They consist of two monoidal structures, which interact via linear distributivity. They provide the semantics for the more minimal \MLLp. Of course, monoidal products connected through negation and de Morgan dualities, modeling \MLL, inherently possess such a linear distributivity. Thus, $*$-autonomous categories are LDCs and the latter are strictly more general, as there are significant examples of LDCs without negation.

The current work generalizes the Fox theorem to the context of LDCs in order to characterize the notable subclass of cartesian LDCs (CLDC). In this section, we introduce all the necessary background on LDCs needed, as developed by co-authors Cockett and Seely in the series of papers \cite{Cockett_Seely_1997,Cockett_Seely_1997_LDC, Cockett_Seely_1999}. 

\subsection{Definition, Mix and Examples}  \hfill\

\begin{definition}\cite[Sec 2.1]{Cockett_Seely_1997_LDC}\label{def:LDC}
A {\bf linearly distributive category}, or a LDC, $(\bX,\ot,\top,\op,\bot)$ is a category \bX\ equipped with two monoidal structures
\begin{itemize}
\item a {\bf tensor} monoidal structure $(\bX, \ot, \top)$
\[ {\alpha_\ot}_{A,B,C}\c (A\ot B)\ot C\rarr A\ot (B\ot C) \quad  {u_\ot^R}_{A}\c A\rarr A\ot \top \quad {u_\ot^L}_{A}\c A\rarr \top\ot A \]
\item a {\bf par} monoidal structure $(\bX, \op, \bot)$
\[{\alpha_\op}_{A,B,C}\c A\op(B\op C)\rarr (A\op B)\op C \quad {u_\op^R}_{A}\c A\op\bot\rarr A \quad {u_\op^L}_{A}\c\bot\op A\rarr A\] 
\end{itemize}
and natural transformations, known as the left and right {\bf linear distributors},
\begin{align*}
\delta^L_{A,B,C}\c A\ot(B\op C)\rarr (A\ot B)\op C && \delta^R_{A,B,C}\c (A\op B)\ot C\rarr A\op (B\ot C),
\end{align*}
satisfying coherence conditions between

\begin{enumerate}[{\bf \footnotesize (LDC.1)}]
    \item unitors and linear distributivitors,\\
\begin{minipage}{.315\textwidth}
\xymatrixrowsep{1.75pc}\xymatrixcolsep{2.75pc}\xymatrix{
A\op B\ar[r]^-{{u_\ot^L}_{A\op B}}\ar[rd]_-{{u_\ot^L}_{A}\op 1_{B}} & \top\ot(A\op B)\ar[d]^-{\delta^L_{\top, A, B}}\\ 
& (\top\ot A)\op B 
}
\end{minipage}
\begin{minipage}{.65\textwidth}
\begin{align*}\label{cc:unit_lineardist}
&{u_\ot^L}_{A\op B}; \delta^L_{\top, A, B} = {u_\ot^L}_{A}\op 1_{B}\\
&{u_\ot^R}_{A\op B}; \delta^R_{A, B,\top} = 1_{A}\op {u_\ot^R}_{B}\\
&\delta^R_{\bot, A, B};{u_\op^L}_{A\ot B} = {u_\op^L}_{A}\ot 1_{B}\\
&\delta^L_{A, B, \bot};{u_\op^R}_{A\ot B} = 1_{A}\ot {u_\op^R}_{B}
\end{align*}
\end{minipage}

\item associators and linear distributors, and 
\begin{fullwidth}
\begin{equation*}\begin{gathered}\label{cc:assoc_lineardist}
\xymatrixrowsep{1.75pc}\xymatrixcolsep{5.25pc}\xymatrix{
(A\ot B) \ot (C\op D)\ar[dd]_-{\delta^L_{A\ot B, C, D}}\ar[r]^-{{\alpha_\ot}_{A,B,C\op D}}& A\ot (B\ot (C\op D))\ar[d]^-{1_{A}\ot \delta^L_{B,C,D}} \\
& A\ot ((B\ot C)\op D)\ar[d]^-{\delta^L_{A,B\ot C, D}} \\
((A\ot B)\ot C) \op D\ar[r]_-{{\alpha_\ot}_{A,B,C}\op 1_{D}} & (A\ot (B\ot C))\op D
}
\end{gathered}\end{equation*}
\vspace{-0.5\baselineskip}
\begin{align*}
&\delta^L_{A\ot B, C, D};({\alpha_\ot}_{A,B,C}\op 1_{D}) = {\alpha_\ot}_{A,B,C\op D}; (1_{A}\ot \delta^L_{B,C,D}); \delta^L_{A,B\ot C, D} \\
&{\alpha_\ot}_{A\op B, C, D};\delta^R_{A, B, C\ot D} = (\delta^R_{A,B,C}\ot 1_{D}); \delta^R_{A,B\ot C, D}; (1_{A}\op {\alpha_\ot}_{B,C,D})\\
&\delta^L_{A, B, C\op D};{\alpha_\op}_{A\ot B, C, D} =  (1_{A}\ot {\alpha_\op}_{B,C,D}); \delta^L_{A,B\op C, D};(\delta^L_{A,B,C}\op 1_{D})\\
&({\alpha_\op}_{A,B,C}\ot 1_{D});\delta^R_{A\op B, C, D} =\delta^R_{A,B\op C, D}; (1_{A}\op \delta^R_{B,C,D});{\alpha_\op}_{A,B,C\ot D}
\end{align*}
\end{fullwidth}
\item left and right linear distributors.
\begin{fullwidth}
\begin{equation*}\begin{gathered}\label{cc:leftright_lineardist}
\xymatrixrowsep{1.75pc}\xymatrixcolsep{1pc}\xymatrix@L=0.35pc{
& (A\op B) \ot (C\op D)\ar[dl]_(.55){\delta^R_{A,B,C\op D}}\ar[dr]^(.55){\delta^L_{A\op B, C, D}} \\
A\op (B\ot (C\op D)\ar[d]_-{1_{A}\op \delta^L_{B,C,D}} & & ((A\op B)\ot C) \op D\ar[d]^-{\delta^R_{A,B,C}\op 1_{D}} \\
A \op ((B\ot C)\op D)\ar[rr]_-{{\alpha_\op}_{A,B\ot C, D}} & &(A\op (B\ot C))\op D
}
\end{gathered}\end{equation*}
\vspace{-0.5\baselineskip}
\begin{align*}
& \delta^L_{A\op B, C, D}; (\delta^R_{A,B,C}\op 1_{D}) = \delta^R_{A,B,C\op D}; (1_{A}\op \delta^L_{B,C,D}); {\alpha_\op}_{A,B\ot C, D} \\
&(\delta^L_{A,B,C}\ot 1_{D});\delta^R_{A\ot B, C, D} = {\alpha_\ot}_{A, B\op C, D}; (1_{A}\ot \delta^R_{B,C,D}); \delta^L_{A,B,C\ot D}
\end{align*}
\end{fullwidth}
\end{enumerate}
\end{definition}

LDCs provide the categorical semantics for non-commutative \MLLp: the tensor product models multiplicative conjunction and the par product models multiplicative disjunction. As in the case of $*$-autonomous categories, the additive and exponential connectives are modeled by additional structure on LDCs, as considered in \cite{Cockett_Seely_1999}. 

We shall mostly be concerned with semantics for commutative \MLLp, given by symmetric monoidal products. 

\begin{definition}\cite[Sec 3]{Cockett_Seely_1997_LDC}\label{def:SLDC}
A LDC is {\bf symmetric}, or a SLDC, if $(\bX, \ot, \top)$ and $(\bX,\op,\bot)$ are SMCs, and the linear distributors interact coherently with the braidings: 
\vspace{-0.5\baselineskip}
\begin{equation}\begin{gathered}\label{cc:lin_dist_braiding}\tag{\bf\footnotesize{SLDC}}
\xymatrixrowsep{1.25pc}\xymatrixcolsep{3.25pc}\xymatrix{
(A\op B) \ot C\ar[r]^-{\delta^R_{A,B,C}} \ar[d]_-{{\sigma_\ot}_{A\op B, C}} & A\op (B\ot C) \\
C\ot (A\op B)\ar[d]_-{1_{C} \ot {\sigma_\op}_{A,B}} & (B\ot C) \op A\ar[u]_-{{\sigma_\op}_{B\ot C, A}}\\
C\ot (B\op A)\ar[r]_-{\delta^L_{C,B,A}} & (C\ot B)\op A \ar[u]_-{{\sigma_\ot}_{C,B}\op 1_{A}} 
}
\end{gathered}\end{equation}
\end{definition}

\begin{remark}
In the symmetric case, instead of requiring \mbox{\bf\footnotesize{(SLDC)}}, we can alternatively drop the requirement that the right (or left) linear distributor exists, drop the conditions from \mbox{\bf\footnotesize{(LDC.1), (LDC.2)}} involving the right (or left) linear distributor, and replace \mbox{\bf\footnotesize{(LDC.3)}} by \mbox{\bf\footnotesize{(LDC.3')}} involving only the left linear distributor given below (or a dual version with right linear distributor). 
\begin{equation*}\begin{gathered}\tag{\bf\footnotesize{LDC.3'}}
\xymatrixrowsep{1pc}\xymatrixcolsep{2.1pc}\xymatrix@L=0.5pc{
(A\ot B)\ot (C\op D)\ar[r]^-{{\alpha_\ot}_{A,B,C\op D}}\ar[d]_-{{\sigma_\ot}_{A,B}\ot {\sigma_\op}_{C,D}} & A\ot (B\ot (C\op D))\ar[r]^-{1_A\ot \delta^L_{B,C,D}} & A\ot ((B\ot C)\op D)\ar[d]^-{1_A\ot {\sigma_\op}_{B\ot C, D}} \\
(B\ot A)\ot (D\op C)\ar[d]_-{{\alpha_\ot}_{B,A,D\op C}} && A\ot (D\op (B\ot C))\ar[d]^-{\delta^L_{A, D, B\ot C}} \\
B\ot (A\ot (D\op C))\ar[d]_-{1_B\ot \delta^L_{A,D,C}} && (A\ot D)\op (B\ot C)\ar[d]^-{{\sigma_\op}_{A\ot D, B\ot C}} \\
B\ot ((A\ot D)\op C)\ar[r]_-{1_B\ot {\sigma_\op}_{A\ot D, C}} & B\ot (C\op (A\ot D))\ar[r]_-{\delta^L_{B,C,A\ot D}} & (B\ot C)\op (A\ot D) 
}
\end{gathered}\end{equation*}
\vspace{-0.5\baselineskip}
\begin{align*}
{\alpha_\ot}_{A,B,C\op D};& (1_A\ot \delta^L_{B,C,D}); (1_A\ot {\sigma_\op}_{B\ot C, D}); \delta^L_{A, D, B\ot C}; {\sigma_\op}_{A\ot D, B\ot C} \\
&=({\sigma_\ot}_{A,B}\ot {\sigma_\op}_{C,D}); {\alpha_\ot}_{B,A,D\op C}; (1_B\ot \delta^L_{A,D,C}); (1_B\ot {\sigma_\op}_{A\ot D, C}); \delta^L_{B,C,A\ot D} \\[0.5pc]
{\sigma_\ot}_{A\op B, C\op D}; & \delta^L_{C\op D, A, B}; ({\sigma_\ot}_{C\op D, A}\op 1_B); (\delta^L_{A, C, D}\op 1_B); {\alpha_\op^{-1}}_{A\ot C, D, B} \\
&=\delta^L_{A\op B, C, D}; ({\sigma_\op}_{A\op B, C}\op 1_D); (\delta^L_{C,A,B}\op 1_D); {\alpha_\op^{-1}}_{C\ot A, B, D}; ({\sigma_\ot}_{C,A}\op {\sigma_\op}_{B, D})
\end{align*}
This remains a SLDC as we can define the missing linear distributor via \mbox{\bf\footnotesize{(SLDC)}} and subsequently prove all the missing coherence conditions. In particular, we use \mbox{\bf\footnotesize{(LDC.3'})} to prove \mbox{\bf\footnotesize{(LDC.3)}}.
\end{remark}

Many models of linear logic also satisfy an additional structural rule: MIX. With the cut rule present, satisfying MIX is equivalent to the additional axiom $\bot \vdash \top$ \cite[Lem 6.1]{Cockett_Seely_1997}. In terms of LDCs, this amounts to the existence of a {\bf nullary mix map} $m\c \bot\rarr\top$ such that the two possible induced $A\ot B\rarr A\op B$ maps are equal.

\begin{definition}\cite[Def 6.2]{Cockett_Seely_1997}
A LDC is {\bf mix} if there is a map $m\c\bot\rarr\top$ such that 
\begin{equation}\begin{gathered}\label{cc:mix_LDC}\tag{\bf\footnotesize{MIX}}
\xymatrixrowsep{1.25pc}\xymatrixcolsep{3.5pc}\xymatrix{
A\ot B\ar[r]^-{1_A\ot {u^L_\op}_B^{-1}}\ar[d]_-{{u^R_\op}_A^{-1}\ot 1_B}& A\ot (\bot \op B)\ar[r]^-{1_A\ot(m\op 1_B)} & A\ot (\top \op B)\ar[d]^-{\delta^L_{A, \top, B}} \\
(A\op \bot)\ot B\ar[d]_-{(1_A\op m)\ot 1_B} & & (A\ot \top)\op B \ar[d]^-{{u^R_\ot}_A^{-1}\op 1_B} \\
(A\op \top)\ot B \ar[r]_-{\delta^R_{A,\top,B}} & A\op (\top \ot B)\ar[r]_-{1_A\op {u^L_\ot}_B^{-1}} & A\op B
}
\end{gathered}\end{equation}
\end{definition}


In the process of developing categorical semantics for classical logic based on LDCs, F\"{u}hrmann and Pym gave a useful way to determine whether or not a LDC is mix. 

\begin{lemma}\cite[Lem 3.10]{Fuhrmann_Pym_2007}\label{lem:prove_mix}
A LDC with a map $m\c\bot\rarr\top$ is mix if and only if
\begin{equation}\label{equation_mix_LDC}
    (1_{\bot}\ot m);{u_\ot^R}_{\bot}^{-1} =(m\ot 1_{\bot});{u_\ot^L}_{\bot}^{-1}\c\bot\ot\bot\rarr\bot
\end{equation}
\end{lemma}

If the LDC is symmetric, the above lemma was then utilized to determine that, if $\bot$ has a $\ot$-comonoid structure (or if $\top$ has a $\op$-monoid structure), the LDC is mix.

\begin{theorem}\cite[Thm 3.11]{Fuhrmann_Pym_2007}\label{thm:prove_mix}
Every SLDC with a $\ot$-comonoid $\langle \bot, \Delta_\bot\c \bot\rarr\bot\ot\bot , e_\bot\c \bot\rarr \top\rangle$ is mix with $m=e_\bot$. Dually, every SLDC with a $\op$-monoid $\langle \top, \nabla_\top\c \top\op\top\rarr\top, u_\top\c \bot\rarr\top\rangle$ is mix with $m=u_\top$. 
\end{theorem}

This MIX rule can be strengthened further, by adding the axiom $\top\vdash \bot$. Categorically, this means $\bot\cong\top$ and we call such LDCs isomix.

\begin{definition}{\cite[Def 6.5]{Cockett_Seely_1997}}
A LDC is {\bf isomix} if it is mix and $m$ is an isomorphism.
\end{definition}

In fact, we do not need to check \mbox{\bf\footnotesize{(MIX)}} if the two monoidal units are isomorphic:

\begin{lemma}{\cite[Lem 6.6]{Cockett_Seely_1997}}
A LDC where $\top\cong\bot$ is isomix.
\end{lemma}

A further level of ``degeneracy'' is obtained by considering LDCs where tensor and par are isomorphic to one another.

\begin{definition}\cite[Sec 2.3]{Cockett_Comfort_Srinivasan_2021}\label{def:compact_LDC}
A LDC is {\bf compact} if it is isomix and the mix maps are isomorphisms. 
\end{definition}

While LDCs take multiplicative conjunction and disjunction as primitive, it is also important to consider linear negation in this categorical framework, in order to be able to model \MLL\ and recover $*$-autonomous categories. 

\begin{definition}\cite[Def A.5]{Cockett_Seely_1999}
A {\bf complementation pair} $(A,A^c, \gamma, \tau)$ in a LDC consists of objects $A$ and $A^c$, equipped with maps
\[ \gamma\c A\ot A^c \rarr \bot \qquad\qquad \tau\c \top\rarr A^c\op A\]
such that 
\begin{equation*}\label{cc:complement_LDC}
\begin{gathered}\xymatrixrowsep{1.25pc}\xymatrixcolsep{2.75pc}\xymatrix{
A\ar[d]_-{{u^R_\ot}_A}\ar[r]^-{1_A} & A & A^c\ar[d]_-{{u^L_\ot}_A^c}\ar[r]^-{1_{A^c}} & A^c\\
A\ot \top\ar[d]_-{1_A\ot \tau} & \bot\ot A\ar[u]_-{{u^L_\op}_A} & \top\ot A^c\ar[d]_-{\tau\ot 1_{A^c}}  & A^c\op \bot\ar[u]_-{{u^R_\op}_{A^c}} \\
A\ot (A^c\op A)\ar[r]_-{\delta^L_{A,A^c,A}} & (A\ot A^c)\op A\ar[u]_-{\gamma\op 1_A} & (A^c\op A)\ot A^c\ar[r]_-{\delta^R_{A^c,A,A^c}} & A^c\op (A\ot A^c)\ar[u]_-{1_{A^c}\op \gamma}
}
\end{gathered}\end{equation*}
\end{definition}

Complementation pairs were renamed {\bf linear duals} by Cockett, Comfort and Srinivasan \cite{Cockett_Comfort_Srinivasan_2021}, to echo duals in autonomous and compact closed categories, and linear adjoints in linear bicategories, as defined by Cockett, Koslowski and Seely \cite{Cockett_Koslowski_Seely_2000}. Indeed, a complementation pairs behaves intuitively like duals in compact closed categories whose unit and counit go between two monoidal structures. The coherence conditions can then be thought of as linearly distributive snake equations. 

Recall that $*$-autonomous categories can be understood as a variant of compact closed categories, whose $*$ does not commute with the monoidal product, instead inducing a par product by de Morgan duality. As such, they are recovered from SLDCs by coherently adding complementation pairs.

\begin{definition}\cite[Def 4.1, 4.3]{Cockett_Seely_1997_LDC}\label{def:SLDC_negation}
A LDC has {\bf negation} if there are object functions $(-)^\perp$ and ${^\perp(-)}$, together with the following parametrized family of maps
\begin{align*}
    &\gamma_A^R\c A\ot A^\perp \rarr \bot & \tau^R_A\c \top\rarr A\op {^\perp A}\\
    &\gamma_A^L\c {^\perp A}\ot A \rarr \bot & \tau^L_A\c \top\rarr A^\perp \op A
\end{align*}
such that $(A, A^\perp, \gamma^R_A, \tau^L_A)$ and $({^\perp A}, A, \gamma^L_A, \tau^R_A)$ form complementation pairs. A SLDC has {\bf negation} if additionally $(-)^\perp = {^\perp(-)}$ and
\[ \gamma_A^L = A^\perp \ot A \xrightarrow{{\sigma_\ot}_{A^\perp, A}} A\ot A^\perp \xrightarrow{\gamma^R_A} \bot \qquad\qquad \tau^L_A = \top \xrightarrow{\tau^R_A} A\op A^\perp \xrightarrow{{\sigma_\op}_{A,A^\perp}} A^\perp\op A\]
\end{definition}

We shall only consider the symmetric case. Mirroring compact closed categories, notice that the above definition does not require $(-)^\perp$ to act on maps, nor does it ask for the many natural isomorphisms which connect the multiplicative connectives and negation together, such as $\bot^\perp\cong \top$ or $A\op B \cong (B^\perp\ot A^\perp)^\perp$. These are in fact a consequence of the above structure.

\begin{lemma}\cite[Lem 4.4]{Cockett_Seely_1997_LDC}
In a SLDC with negation, $(-)^\perp$ is a contravariant functor, which acts on maps as follows: given $f\c A\rarr B$, define $f^\perp\c B^\perp\rarr A^\perp$ by
\begin{equation*}
\begin{gathered}\xymatrixrowsep{1.75pc}\xymatrixcolsep{3.75pc}\xymatrix{
B^\perp\ar[r]^-{{u^R_\ot}_{B^\perp}}\ar[dd]_-{f^\perp} & B^\perp\ot \top \ar[r]^-{1_{B^\perp}\ot \tau^R_A} & B^\perp \ot (A\op A^\perp)\ar[d]^-{1_{B^\perp}\ot (f\op 1_{A^\perp})} \\
&& B^\perp \ot (B\op A^\perp)\ar[d]^-{\delta^L_{B^\perp, B, A^\perp}}\\
A^\perp & \bot\op A^\perp \ar[l]^{{u^L_\op}_{A^\perp}}& (B^\perp\ot B)\op A^\perp\ar[l]^-{{\gamma^L_B}\op 1_{A^\perp}} 
}
\end{gathered}\end{equation*}
\end{lemma}

Then, Cockett and Seely demonstrated that:
\begin{theorem}\cite[Thm 4.5]{Cockett_Seely_1997_LDC}\label{thm:SLDC_starauto}
The notions of SLDCs with negation and $*$-autonomous categories coincide.
\end{theorem}

\begin{example}\label{ex:LDC} Let us list three important classes of examples of LDCs.
\begin{enumerate}
    \item Every monoidal category $(\cX, \os, I)$ can be viewed as a LDC, when taking the tensor and par structures to be equal to the original monoidal structure, i.e. $\ot=\op=\os$ and $\top= \bot = I$. In this case, the linear distributors are just $\os$-associators and we call such LDCs $(\cX, \os, I, \os, I)$ {\bf degenerate}.

    \item As stated in Theorem \ref{thm:SLDC_starauto}, every $*$-autonomous category $(\cX, \os, I, *)$ is a LDC with the tensor structure given by the original monoidal structure and the par structure given by the de Morgan dual, i.e. $A \op B = (B^* \os A^*)^*$ and $\bot = I^*$. 

    \item Cockett and Seely introduced a class of examples known as {\bf shifted tensor} LDCs, the categorical analogue of shift monoids.

    \begin{definition}\cite[Sec 5.2]{Cockett_Seely_1997_LDC}
        Consider a monoidal category $(\bX, \ot, \top)$. An object $\bot\in\bX$ is said to have a {\bf tensor inverse} if there is an object $\bot^{-1}$ equipped with two isomorphisms 
        \[ s^L\c\bot \ot \bot^{-1}\rarr\top \qquad s^R\c \bot^{-1}\ot\bot\rarr\top\]
        such that 
    
    \begin{equation*}\begin{gathered}
    \xymatrixrowsep{1.25pc}\xymatrixcolsep{2.75pc}\xymatrix{
    (\bot^{-1}\ot \bot) \ot \bot^{-1} \ar[rr]^-{{\alpha_\ot}_{\bot^{-1},\bot, \bot{-1}}} \ar[d]_-{s^R\ot 1_{\bot^{-1}}}& &\bot^{-1}\ot(\bot\ot\bot^{-1}) \ar[d]^-{1_{\bot^{-1}}\ot s^L}\\
    \top\ot\bot^{-1}\ar[rd]_-{{u^L_\ot}^{-1}_{\bot^{-1}}} & & \bot^{-1}\ot \top \ar[ld]^-{{u^R_\ot}^{-1}_{\bot^{-1}}}\\
    & \bot^{-1}
    }
    \end{gathered}\end{equation*}
    
    Suppose there is an object $\bot\in\bX$ with a tensor inverse, then define a monoidal product by $A\op B = A\ot (\bot^{-1}\ot B)$, known as the {\bf $\bot$-shifted tensor}. 
    \end{definition}    

    \begin{proposition}\cite[Prop 5.3, 5.4]{Cockett_Seely_1997_LDC}\label{prop:shift_tensor}
        Consider a monoidal category $(\bX, \ot, \top)$ with an object $\bot\in\bX$ with a tensor inverse, then $(\bX, \ot, \top, \op,\bot)$, where $\op$ is the $\bot$-shifted tensor, is a LDC with invertible linear distributors. Moreover, every LDC with invertible linear distributors has a tensor inverse for $\bot$ whose shift tensor product is naturally equivalent to the par product. 
    \end{proposition}
    
\end{enumerate}
\end{example}

\subsection{Linear Functors and Transformations}  \hfill\

\vspace{0.5\baselineskip}
In the same spirit as the Fox theorem, we will need the appropriate notions of functors and transformations between LDCs. These definitions were developed by Cockett and Seely in order to model the exponential modalities $!$ and $?$ of linear logic as a paired functorial entity, while also capturing the idea of lax monoidal functors and transformations between $*$-autonomous categories. 

\begin{definition}\cite[Def 1]{Cockett_Seely_1999}
Let \bX\ and \bY\ be LDCs. A {\bf (bilax) linear functor} $F = (F_\ot, F_\op)\c\bX\rarr\bY$ consists of:
\begin{itemize}
\item a (lax) monoidal functor $(F_\ot, m_\top, m_\ot)\c(\bX,\ot,\top)\rarr(\bY,\ot,\top)$, equipped with
\[ m_\top\c\top\rarr F_\ot(\top) \qquad {m_\ot}_{A, B}\c F_\ot(A)\ot F_\ot(B)\rarr F_\ot(A\ot B)\]

\item a (colax) monoidal functor $(F_\op, n_\bot, n_\op)\c(\bX,\op,\bot)\rarr(\bY,\op,\bot)$, equipped with
\[ n_\bot\c F_\op(\bot)\rarr\bot\qquad {n_\op}_{A,B}\c F_\op(A\op B)\rarr F_\op(A)\op F_\op(B)\]

\item four natural transformations, known as {\bf linear strengths}, 
\begin{align*}
{v_\ot^R}_{A,B} \c F_\ot(A\op B)\rarr F_\op(A)\op F_\ot(B)  \qquad {v_\op^R}_{A,B} \c F_\ot(A)\ot F_\op(B)\rarr F_\op(A\ot B) \\
{v_\ot^L}_{A,B} \c F_\ot(A\op B)\rarr F_\ot(A)\op F_\op(B) \qquad {v_\op^L}_{A,B} \c F_\op(A)\ot F_\ot(B)\rarr F_\op(A\ot B)   
\end{align*}
\end{itemize}
subject to various coherence conditions between the linear strengths and
\begin{enumerate}[{\bf \footnotesize (LF.1)}]
    \item the unitors,
\begin{fullwidth}
\begin{equation*}\label{cc:linear_strength_unit}
\xymatrixrowsep{1.75pc}\xymatrixcolsep{4.25pc}\xymatrix{
F_\ot(\bot\op A)\ar[r]^-{F_\ot({u^L_\op}_A)}\ar[d]_-{{\nu^R_\ot}_{\bot, A}} & F_\ot(A) \\
F_\op(\bot) \op F_\ot(A)\ar[r]_-{n_\bot \op 1_{F_\ot(A)}} & \bot \op F_\ot(A)\ar[u]_-{{u^L_\op}_{F_\ot(A)}} 
}
\end{equation*}
\vspace{-\baselineskip}
\begin{align*}
&{\nu^R_\ot}_{\bot, A}; (n_\bot \op 1_{F_\ot(A)}); {u^L_\op}_{F_\ot(A)} = F_\ot({u^L_\op}_A) \\
&{\nu^L_\ot}_{A, \bot}; (1_{F_\ot(A)}\op n_\bot); {u^R_\op}_{F_\ot(A)} = F_\ot({u^R_\op}_A)\\
&{u^L_\ot}_{F_\op(A)};(m_\top \ot 1_{F_\op(A)}); {\nu^R_\op}_{\top, A} = F_\op({u^L_\ot}_A)\\
&{u^R_\ot}_{F_\op(A)};( 1_{F_\op(A)}\ot m_\top); {\nu^L_\op}_{A,\top} = F_\op({u^R_\ot}_A)
\end{align*}
\end{fullwidth}
    \item the associators,
\begin{fullwidth}
\begin{equation*}\label{cc:linear_strength_assoc}
\xymatrixrowsep{1.75pc}\xymatrixcolsep{7.5pc}\xymatrix{
F_\ot(A\op (B\op C))\ar[r]^-{F_\ot({\alpha_\op}_{A,B,C})}\ar[d]_-{{\nu^R_\ot}_{A, B\op C}} &F_\ot((A\op B)\op C)\ar[d]^-{{\nu^R_\ot}_{A\op B, C}} \\
F_\op(A)\op F_\ot(B\op C)\ar[d]_-{1_{F_\op(A)}\op {\nu^R_\ot}_{B, C}} & F_\op(A\op B) \op F_\ot(C)\ar[d]^-{{n_\op}_{A,B}\op 1_{F_\ot(C)}}\\
 F_\op(A) \op (F_\op(B) \op F_\ot(C))\ar[r]_-{{\alpha_\op}_{F_\op(A), F_\op(B), F_\ot(C)}} & (F_\op(A)\op F_\op(B)) \op F_\ot(C)
}
\end{equation*}
\vspace{-\baselineskip}
\begin{align*}
&{\nu^R_\ot}_{A,B\op C}; (1_{F_\op(A)}\op {\nu^R_\ot}_{B,C});{\alpha_\op}_{F_\op(A), F_\op(B), F_\ot(C)} = F_\ot({\alpha_\op}_{A,B,C}); {\nu^R_\ot}_{A\op B, C}; ({n_\op}_{A,B}\op 1_{F_\ot(C)}) \\
&{\nu^L_\ot}_{A,B\op C}; (1_{F_\ot(A)}\op {n_\op}_{B, C});{\alpha_\op}_{F_\ot(A), F_\op(B), F_\op(C)} = F_\ot({\alpha_\op}_{A,B,C}); {\nu^L_\ot}_{A\op B, C}; ({\nu^L_\op}_{A,B}\op 1_{F_\op(C)})\\
&{\alpha_\ot}_{F_\ot(A), F_\ot(B), F_\op(C)}; (1_{F_\ot(A)} \ot {\nu^R_\op}_{B,C}); {\nu^R_\op}_{A, B\ot C} = ({m_\ot}_{A,B} \ot 1_{F_\op(C)}); {\nu^R_\op}_{A\ot B, C}; F_\op({\alpha_\ot}_{A, B, C})  \\
& {\alpha_\ot}_{F_\op(A), F_\ot(B), F_\ot(C)};(1_{F_\op(A)}\ot {m_\ot}_{B,C}); {\nu^L_\op}_{A, B\ot C}= ({\nu^L_\op}_{A,B}\ot 1_{F_\ot(C)}); {\nu^L_\op}_{A\ot B, C}; F_\op({\alpha_\ot}_{A,B,C})
\end{align*}
\end{fullwidth}
    \item the linear distributors,
\begin{fullwidth}
\begin{equation*}\label{cc:linear_strength_dist}
\xymatrixrowsep{1.75pc}\xymatrixcolsep{5.75pc}\xymatrix{
F_\ot(A\op B) \ot F_\ot(C)\ar[d]_-{{\nu^R_\ot}_{A,B} \ot 1_{F_\ot(C)}}\ar[r]^-{{m_\ot}_{A\op B, C}} & F_\ot((A\op B)\ot C)\ar[d]^-{F_\ot(\delta^R_{A,B,C})} \\
(F_\op(A) \op F_\ot(B)) \ot F_\ot(C) \ar[d]_-{\delta^R_{F_\op(A), F_\ot(B), F_\ot(C)}} & F_\ot(A\op (B\ot C)) \ar[d]^-{{\nu^R_\ot}_{A, B\ot C}}\\
F_\op(A) \op (F_\ot(B) \ot F_\ot(C))\ar[r]_-{1_{F_\op(A)}\op {m_\ot}_{B,C}} & F_\op(A) \op F_\ot(B\ot C)
}
\end{equation*}
\vspace{-\baselineskip}
\begin{align*}
&{m_\ot}_{A\op B, C}; F_\ot(\delta^R_{A,B,C}); {\nu^R_\ot}_{A, B\ot C}  = ({\nu^R_\ot}_{A,B} \ot 1_{F_\ot(C)}); \delta^R_{F_\op(A), F_\ot(B), F_\ot(C)}; (1_{F_\op(A)}\op {m_\ot}_{B,C}) \\
&{m_\ot}_{A, B\op C}; F_\ot(\delta^L_{A,B,C}); {\nu^L_\ot}_{A\ot B, C}  = (1_{F_\ot(A)}\ot {\nu^L_\ot}_{B,C}); \delta^L_{F_\ot(A), F_\ot(B), F_\op(C)}; ({m_\ot}_{A,B}\op 1_{F_\op(C)}) \\
& {\nu^R_\op}_{A, B\op C}; F_\op(\delta^L_{A, B, C}); {n_\op}_{A\ot B, C} = (1_{F_\ot(A)} \ot {n_\op}_{B,C}); \delta^L_{F_\ot(A), F_\op(B), F_\op(C)}; ({\nu^R_\op}_{A,B}\op 1_{F_\op(C)}) \\
& {\nu^L_\op}_{A\op B, C}; F_\op(\delta^R_{A, B, C}); {n_\op}_{A, B\ot C} = ({n_\op}_{A,B}\op 1_{F_\ot(C)}); \delta^R_{F_\op(A), F_\op(B), F_\ot(C)}; (1_{F_\op(A)}\op {\nu^L_\op}_{B,C}) 
\end{align*}
\end{fullwidth}
    \item the other linear strengths via associators, and 
\begin{fullwidth}
\begin{equation*}\label{cc:linear_stregnth_assoc_linear_strength}
\xymatrixrowsep{1.75pc}\xymatrixcolsep{7.75pc}\xymatrix{
F_\ot(A\op (B\op C))\ar[r]^-{F_\ot({\alpha_\op}_{A,B,C})}\ar[d]_-{{\nu^R_\ot}_{A, B\op C}} & F_\ot((A\op B)\op C)\ar[d]^-{{\nu^L_\ot}_{A\op B, C}} \\
F_\op(A) \op F_\ot(B\op C)\ar[d]_-{1_{F_\op(A)}\op {\nu^L_\ot}_{B,C}} & F_\ot(A\op B)\op F_\op(C)\ar[d]^-{{\nu^R_\ot}_{A, B} \op 1_{F_\op(C)}}\\
 F_\op(A) \op (F_\ot(B) \op F_\op(C))\ar[r]_-{{\alpha_\op}_{F_\op(A), F_\ot(B), F_\op(C)}}& (F_\op(A) \op F_\ot(B))\op F_\op(C)
}
\end{equation*}
\vspace{-\baselineskip}
\begin{align*}
& {\nu^R_\ot}_{A, B\op C}; (1_{F_\op(A)}\op {\nu^L_\ot}_{B,C});{\alpha_\op}_{F_\op(A), F_\ot(B), F_\op(C)} = F_\ot({\alpha_\op}_{A,B,C});{\nu^L_\ot}_{A\op B, C}; ({\nu^R_\ot}_{A, B} \op 1_{F_\op(C)}) \\
& {\alpha_\ot}_{F_\ot(A),F_\op(B), F_\ot(C)}; (1_{F_\ot(A)}\ot {\nu^L_\op}_{B, C}) ; {\nu^R_\op}_{A, B\ot C} = ({\nu^R_\op}_{A, B} \ot 1_{F_\ot(C)}); {\nu^L_\op}_{A\ot B, C}; F_\op({\alpha_\ot}_{A, B, C})
\end{align*}
\end{fullwidth}
    \item the other linear strengths via linear distributors.
\begin{fullwidth}
\begin{equation*}\label{cc:linear_stregnth_dist_linear_strength}
\xymatrixrowsep{1.75pc}\xymatrixcolsep{5.5pc}\xymatrix{
F_\ot(A) \ot F_\ot (B\op C)\ar[r]^-{1_{F_\ot(A)}\ot {\nu^R_\ot}_{B,C}}\ar[d]_-{{m_\ot}_{A, B
\op C}} & F_\ot(A) \ot (F_\op (B) \op F_\ot(C))\ar[d]^-{\delta^L_{F_\ot(A), F_\op(B), F_\ot(C)}} \\
F_\ot(A\ot (B\op C))\ar[d]_-{F_\ot(\delta^L_{A, B, C})} & (F_\ot(A) \ot F_\op(B)) \op F_\ot(C)\ar[d]^-{{\nu^R_\op}_{A,B} \op 1_{F_\ot(C)}} \\
F_\ot((A\ot B)\op C) \ar[r]_-{{\nu^R_\ot}_{A\ot B, C}}& F_\op(A\ot B) \op F_\ot(C)
}
\end{equation*}
\vspace{-\baselineskip}
\begin{align*}
& {m_\ot}_{A, B\op C}; F_\ot(\delta^L_{A, B, C}); {\nu^R_\ot}_{A\ot B, C} = (1_{F_\ot(A)}\ot {\nu^R_\ot}_{B,C}); \delta^L_{F_\ot(A), F_\op(B), F_\ot(C)}; ({\nu^R_\op}_{A,B} \op 1_{F_\ot(C)})\\
& {m_\ot}_{A\op B, C}; F_\ot(\delta^R_{A, B, C}); {\nu^L_\ot}_{A, B\ot C} = ({\nu^L_\ot}_{A,B}\ot 1_{F_\ot(C)}); \delta^R_{F_\ot(A), F_\op(B), F_\ot(C)}; (1_{F_\ot(A)}\op {v^L_\op}_{B,C})\\
& {\nu^R_\op}_{A\op B, C}; F_\op(\delta^R_{A, B, C}); {n_\op}_{A, B\ot C} = ({\nu^R_\ot}_{A, B}\ot 1_{F_\op(C)}); \delta^R_{F_\op(A), F_\ot(B), F_\op(C)}; (1_{F_\op(A)}\op {\nu^R_\op}_{B, C}) \\
& {\nu^L_\op}_{A, B\op C}; F_\op(\delta^L_{A, B, C}); {n_\op}_{A\ot B, C} = (1_{F_\op(A)} \ot {\nu^L_\ot}_{B, C}); \delta^L_{F_\op(A), F_\ot(B), F_\op(C)}; ({\nu^L_\op}_{A, B}\op 1_{F_\op(C)})
\end{align*}
\end{fullwidth}
\end{enumerate}

If \bX\ and \bY\ are SLDCs, then a (bilax) linear functor $F=(F_\ot, F_\op)$ is {\bf symmetric} if $F_\ot$ is a symmetric (lax) monoidal functor, $F_\op$ is a symmetric (colax) monoidal functor, and the linear strengths interact coherently with the braidings:
\begin{equation}\begin{gathered}\label{cc:linear_strength_braiding}\tag{\bf\footnotesize{SLF}}
\xymatrixrowsep{1.75pc}\xymatrixcolsep{2.25pc}\xymatrix{
F_\ot(A\op B) \ar[r]^-{{v_\ot^R}_{A,B}}\ar[d]_-{F_\ot({\sigma_\op}_{A,B})} & F_\op(A)\op F_\ot(B) \ar[d]_-{{\sigma_\op}_{F_\op(A), F_\ot(B)}} \\
F_\ot (B\op A) \ar[r]_-{{v_\ot^L}_{B, A}} & F_\ot(B) \op F_\op(A)}
\qquad
\xymatrixrowsep{1.75pc}\xymatrixcolsep{2.25pc}\xymatrix{
F_\ot(A)\ot F_\op(B)\ar[r]^-{{v_\op^R}_{A,B}}\ar[d]^-{{\sigma_\ot}_{F_\ot(A),F_\op(B)}} & F_\op(A\ot B)\ar[d]^-{F_\op({\sigma_\ot}_{A,B})}\\
F_\op(B)\ot F_\ot(A)\ar[r]_-{{v_\op^L}_{B,A}} & F_\op(B\ot A)
}
\end{gathered}\end{equation}

\end{definition}

\begin{remark}
In the symmetric case, instead of requiring \mbox{\bf\footnotesize{(SLF)}} to hold, we can alternatively drop the requirement that the left (or right) linear strengths exist, drop the conditions from \mbox{\bf\footnotesize{(LF.1), (LF.2), (LF.3), (LF.5)}} involving the left (or right) linear strengths, and drop \mbox{\bf\footnotesize{(LF.4)}} entirely. This remains a (bilax) linear functor as we can define the missing pair of linear strengths via \mbox{\bf\footnotesize{(SLF)}} and subsequently prove all the missing coherence conditions. In particular, we use \mbox{\bf\footnotesize{(LF.2})} to prove \mbox{\bf\footnotesize{(LF.4)}}. Therefore, when defining symmetric (bilax) linear functors, we often only give one pair of linear strengths.
\end{remark}

\begin{definition}\cite[Def 3]{Cockett_Seely_1999}\label{def:linear_trans}
Let $F, G\c\bX\rarr\bY$ be (bilax) linear functors. A {\bf linear transformation} $\alpha\c F\Rarr G$ consists of:
\begin{itemize}
\item a monoidal transformation $\alpha_\ot\c (F_\ot, m_{\top}^{F}, m_\ot^{F}) \Rarr (G_\ot, m_{\top}^{G}, m_\ot^{G})$, and
\item a monoidal transformation $\alpha_\op\c (G_\op, n_{\bot}^{G}, n_\op^{G})\Rarr (F_\op, n_{\bot}^{F}, n_\op^{F})$,
\end{itemize}
which commute coherently with the linear strengths of $F$ and $G$:
\begin{equation}\begin{gathered}\label{cc:linear_trans_linear strength}\tag{\bf\footnotesize{LT}}
\xymatrixrowsep{1.75pc}\xymatrixcolsep{4.75pc}\xymatrix{
F_\ot(A\op B) \ar[r]^-{{\alpha_\ot}_{A\op B}} \ar[dd]_-{{v_\ot^R}^{F}_{A,B}} & G_\ot(A\op B)\ar[d]^-{{v_\ot^R}^{G}_{A,B}} \\
& G_\op (A) \op G_\ot(B) \ar[d]^-{{\alpha_\op}_{A}\op 1_{G_\ot(B)}} \\
F_\op(A)\op F_\ot(B)\ar[r]_-{1_{F_\op(A)} \op {\alpha_\ot}_{B}} & F_\op(A) \op G_\ot(B)
}
\end{gathered}\end{equation}
\vspace{-\baselineskip}
\begin{align*}
&{\alpha_\ot}_{A\op B} ; {v_\ot^R}^{G}_{A,B} ; ({\alpha_\op}_{A} \op 1_{G_\ot (B)}) = {v_\ot^R}^{F}_{A,B} ; (1_{F_\op(A)} \op {\alpha_\ot}_{B})\\
&{\alpha_\ot}_{A\op B} ; {v_\ot^L}^{G}_{A,B} ;  (1_{G_\ot(A)} \op {\alpha_\op}_{B}) = {v_\ot^L}^{F}_{A,B} ; ({\alpha_\ot}_{A}\op 1_{F_\op(B)}) \\
& ({\alpha_\ot}_{A} \ot 1_{G_\op(B)}) ; {v_\op^R}^{G}_{A,B} ; {\alpha_\op}_{A\ot B} = (1_{F_\ot(A)}\ot {\alpha_\op}_{B}) ; {v_\op^R}^{F}_{A,B} \\
& (1_{G_\op(A)}\ot {\alpha_\ot}_{B}) ; {v_\op^L}^{G}_{A,B} ; {\alpha_\op}_{A\ot B} = ({\alpha_\op}_{A} \ot 1_{F_\ot(B)}) ; {v_\op^L}^{F}_{A,B}
\end{align*}
\end{definition}

\begin{proposition}\cite[Prop 4]{Cockett_Seely_1999}\label{prop:2cat_LDC}
LDCs, bilax linear functors, and linear transformations form a 2-category, which is denoted by \LDC. Restricting to SLDCs and symmetric bilax linear functors also gives a 2-category, denoted \SLDC. 
\end{proposition}

We will also consider linear functors whose component functors are equal. These were first defined by \cite{Blute_Panangaden_Slanov_2012} and then further explored in \cite{Cockett_Comfort_Srinivasan_2021}.

\begin{definition}\cite[Def 3.1]{Cockett_Comfort_Srinivasan_2021}
Consider LDCs \bX\ and \bY, a (bilax) linear functor $F=(F_\ot,F_\op)\c \bX\rarr\bY$ is {\bf Frobenius} if
\[ F_\ot = F_\op, \qquad {\nu^R_\ot}_{A,B} = {\nu^L_\ot}_{A,B} = {n_\op}_{A,B},\quad\text{and}\quad{\nu^R_\op}_{A,B} = {\nu^L_\op}_{A,B} = {m_\ot}_{A,B}\]
\end{definition}

Given the degeneracy, we can give an alternative characterization.
\begin{proposition}\cite[Lem 3.2]{Cockett_Comfort_Srinivasan_2021}
Consider LDCs \bX\ and \bY, then the following notions coincide:
\begin{itemize}
    \item Frobenius (bilax) linear functors $F = (F_\ot, F_\op)\c \bX\rarr\bY$, and 
    \item $\ot$-(lax) monoidal and $\op$-(colax) monoidal functors $(F, m_\ot, m_\top, n_\op, n_\bot)\c \bX\rarr\bY$ satisfying
\end{itemize}
\vspace{-0.5\baselineskip}
\begin{equation}\label{cc:Frobenius_linear_functor}
\begin{gathered}\xymatrixrowsep{1.75pc}\xymatrixcolsep{4.75pc}\xymatrix{
F(A)\ot F(B\op C)\ar[r]^-{{m_\ot}_{A,B\op C}}\ar[d]_-{1_{F(A)}\ot {n_\op}_{B,C}} & F(A\ot (B\op C))\ar[d]^-{F(\delta^L_{A,B,C})} \\
F(A) \ot (F(B)\op F(C))\ar[d]_-{\delta^L_{F(A), F(B), F(C)}} & F((A\ot B)\op C)\ar[d]^-{{n_\op}_{A\ot B, C}} \\
(F(A) \ot F(B))\op F(C)\ar[r]_-{{m_\ot}_{A,B}\op 1_{F(C)}} & F(A\ot B)\op F(C)
}
\end{gathered}\end{equation}
\vspace{-\baselineskip}
\begin{align*}
& {m_\ot}_{A,B\op C}; F(\delta^L_{A,B,C}); {n_\op}_{A\ot B, C} = (1_{F(A)}\ot {n_\op}_{B,C}); \delta^L_{F(A), F(B), F(C)}; ({m_\ot}_{A,B}\op 1_{F(C)}) \\
& {m_\ot}_{A\op B, C}; F(\delta^R_{A, B, C}); {n_\op}_{A, B\ot C} = ({n_\op}_{A,B}\ot 1_{F(C)}); \delta^R_{F(A), F(B), F(C)}; (1_{F(A)}\op {m_\ot}_{B, C})
\end{align*}

\end{proposition}

If the LDCs are mix, the definition of Frobenius (bilax) linear functors can be extended to guarantee they preserve the nullary mix maps.
\begin{definition}\cite[Def 3.4]{Cockett_Comfort_Srinivasan_2021}
Consider mix LDCs \bX\ and \bY. A Frobenius (bilax) linear functor $F=(F,F)\c \bX\rarr\bY$ is {\bf mix} if
\begin{equation}\label{cc:mix_Frobenius_linear_functor}\tag{\bf\footnotesize{MixFLF}}
\begin{gathered}\xymatrixrowsep{1.75pc}\xymatrixcolsep{1.75pc}\xymatrix@L=0.3pc{
F(\bot)\ar[r]^-{n_\bot}\ar[d]_-{F(m)} & \bot\ar[d]^-{m} \\
F(\top) & \top\ar[l]^-{m_\top}
}
\end{gathered}\end{equation}
\end{definition}

If we then consider linear transformations between Frobenius linear functors, they often correspond to a natural isomorphism paired with its inverse.

\begin{lemma}\cite[Lem 3.10]{Cockett_Comfort_Srinivasan_2021}\label{lem:Frobenius_linear_trans_2}
Suppose $F,G\c \bX\rarr\bY$ are Frobenius (bilax) linear functors and $\alpha\c F\Rarr G$ is a $\ot$-monoidal and $\op$-monoidal natural isomorphism, then $(\alpha, \alpha^{-1})\c F\Rarr G$ is a linear transformation.
\end{lemma}

As Frobenius (bilax) linear functors compose and identity linear functors are Frobenius, restricting to such 1-cells determines a sub-2-category.
\begin{proposition}
LDCs, Frobenius bilax linear functors and linear transformations form a sub-2-category of \LDC, denoted by \FLDC.
\end{proposition}

\subsection{Cartesian Linearly Distributive Categories}  \hfill\

\vspace{0.5\baselineskip}
We now introduce cartesian LDCs, the subclass we are interested in characterizing via the linearly distributive Fox theorem.

\begin{definition}\cite[Sec 2]{Cockett_Seely_1997_LDC}
A {\bf cartesian linearly distributive category}, or CLDC, $(\bX, \times, \bone, +, \bzero)$ is a SLDC whose tensor monoidal structure $(\bX, \times, \bone)$ is cartesian and whose par monoidal structure $(\bX, +, \bzero)$ is cocartesian. 
\end{definition}

The development of the linearly distributive Fox theorem renewed interest in CLDCs, leading them to be further explored in their own right by the author with co-author Lemay in the article ``Cartesian linearly distributive categories: Revisited'' \cite{Kudzman-Blais_Lemay_2025}. This work examines the distinctive structures and properties of CLDCs, discusses key classes of examples, and reassesses previously proposed instances of CLDCs. We recommend reading both articles, as they complement each other, though each is written to be accessible independently. 

\begin{proposition}
There is a sub 2-category of \LDC\ consisting of CLDCs, symmetric linear functors, known as {\bf cartesian linear functors}, and linear transformations, denoted by \CLDC. Moreover, there is a sub 2-category of \FLDC\ consisting of CLDCs, symmetric Frobenius linear functors and linear transformations, denoted by \FCLDC.
\end{proposition}

\section{Duoidal Categories}\label{sec:duoidal_categories}

LDCs are of course not the only categories with two monoidal structures. Duoidal categories are another variant, characterized by two monoidal products and units that interact via additional structure maps, most notably the interchange law. Unlike LDCs, the duoidal structure arises canonically when considering finite products and coproducts. The theory of duoidal categories will play a crucial role in this paper, so we present all the necessary background in this section.

The earliest appearance of a form of the interchange law between monoidal structures is found in Joyal and Street's work on braided monoidal categories \cite{Joyal_Street_1993}. They demonstrated that a duoidal category with isomorphisms for structure maps corresponds to a braided monoidal category via a categorical Eckmann-Hilton argument. Subsequently, Balteanu, Fiedorowicz, Schwänzl, and Vogt defined 2-fold monoidal categories as monoids in the monoidal category of monoidal categories \cite{Balteanu_Fiedorowicz_Schwaenzl_Vogt_2003}. In their definition, the interchange law is not required to be an isomorphism, but the unit objects of both monoidal structures are assumed to coincide, and the monoidal structures themselves are strict. Building on this foundation, Forcey, Siehler, and Sowers generalized the definition by allowing distinct unit objects and non-strict associators, while retaining equalities for unitors and certain structure maps \cite{Forcey_Siehler_Sowers_2007}. 

The currently accepted definition of duoidal category was first introduced by Aguiar and Mahajan under the name 2-monoidal categories \cite{Aguiar_Mahajan_2010}. The background on duoidal categories presented in this section is almost entirely based on Chapter 6 of their monograph, and readers seeking additional details are encouraged to consult it. The now-standard name for duoidal categories was suggested by Street and first appeared in the literature in \cite{Batanin_Markl_2012, Booker_Street_2013}.

\subsection{Definition and Examples}  \hfill\

\begin{definition}\cite[Def 6.1, 6.3]{Aguiar_Mahajan_2010}\label{def:duoidal_cat}
A {\bf duoidal category} $(\cX, \diamond, I,\star, J)$ is category \cX\ with two monoidal structures $(\cX, \diamond, I)$ and $(\cX, \star, J)$, equipped with morphisms 
\begin{alignat*}{3}
\Delta_{I}  &\c I\rarr I\star I &\qquad  \mu_{J}& \c J\diamond J\rarr J  &\qquad  \iota  & \c I\rarr J
\end{alignat*}
and an {\bf interchange} natural transformation 
\[ \zeta_{A, B, C, D}\c (A\star B)\diamond(C\star D)\rarr(A\diamond C)\star(B\diamond D)\]
such that $(J, \mu_{J}, \iota)$ is a $\diamond$-monoid, $(I, \Delta_{I}, \iota)$ is a $\star$-comonoid, and the interchange maps interact coherently with 
\begin{enumerate}[{\bf \footnotesize (DUO.1)}]
    \item the associators, and 
\begin{fullwidth}
\begin{equation*}\label{cc:assoc_interchange}
\xymatrixrowsep{1.75pc}\xymatrixcolsep{6.75pc}\xymatrix{
((A\star B) \diamond (C\star D)) \diamond (E\star F) \ar[r]^-{{\alpha_\diamond}_{A\star B, C\star D, E\star F}} \ar[d]_-{\zeta_{A, B, C, D}\diamond 1_{E\star F}} & (A\star B) \diamond ((C\star D) \diamond (E\star F)) \ar[d]^-{1_{A\star B}\diamond \zeta_{C, D, E, F}} \\
((A\diamond C)\star (B\diamond D)) \diamond (E\star F) \ar[d]_-{\zeta_{A\diamond C, B\diamond D, E, F}} & (A\star B) \diamond ((C\diamond E)\star (D\diamond F)) \ar[d]^-{\zeta_{A, B, C\diamond E, D\diamond F}} \\
((A\diamond C)\diamond E)\star ((B\diamond D)\diamond F) \ar[r]_-{{\alpha_\diamond}_{A, C, E}\star {\alpha_\diamond}_{B, D, F}} & (A\diamond (C\diamond E))\star (B\diamond (D\diamond F))
}
\end{equation*}
\vspace{-\baselineskip}
\begin{align*}
&{\alpha_\diamond}_{A\star B, C\star D, E\star F}; (1_{A\star B}\diamond \zeta_{C, D, E, F}); \zeta_{A, B, C\diamond E, D\diamond F} = (\zeta_{A, B, C, D}\diamond 1_{E\star F}); \zeta_{A\diamond C, B\diamond D, E, F}; ({\alpha_\diamond}_{A, C, E}\star {\alpha_\diamond}_{B, D, F}) \\
&\zeta_{A\star B, C D\star E, F}; (\zeta_{A, B, D, E}\star 1_{C\diamond F}); {\alpha_\star}_{A\diamond D, B\diamond E, C\diamond F} = ({\alpha_\star}_{A,B,C}\diamond{\alpha_\star}_{D, E, F}); \zeta_{A, B\star C, D, E\star F}; (1_{A\diamond D}\star\zeta_{B, C, E, F})
\end{align*}
\end{fullwidth}

\item the unitors.
\begin{fullwidth}
\begin{equation*}\label{cc:unit_interchange}
\xymatrixrowsep{1.75pc}\xymatrixcolsep{3.75pc}\xymatrix{
A\star B\ar[r]^-{{\lambda_\diamond}_{A\star B}} \ar[d]_-{{\lambda_\diamond}_{A}\star {\lambda_\diamond}_{B}}& I \diamond (A\star B)\ar[d]^-{\Delta_I\diamond 1_{A\star B}} \\
(I\diamond A)\star (I\diamond B) & (I\star I)\diamond (A\star B)\ar[l]^-{\zeta_{I,I,A,B}}
}
\end{equation*}
\vspace{-\baselineskip}
\begin{align*}
& {\lambda_\diamond}_{A\star B}; (\Delta_I\diamond 1_{A\star B}); \zeta_{I,I,A,B} = {\lambda_\diamond}_{A}\star {\lambda_\diamond}_{B} \quad\quad\quad {\rho_\diamond}_{A\star B}; (1_{A\star B}\diamond \Delta_I); \zeta_{A,B, I, I} = {\rho_\diamond}_{A}\star {\rho_\diamond}_{B} \\
& ({\lambda_\star}_{A}\diamond{\lambda_\star}_{B}); \zeta_{J, A, J, B}; (\mu_J\star 1_{A\diamond B}) = {\lambda_\star}_{A\diamond B} \quad\quad ({\rho_\star}_{A}\diamond{\rho_\star}_{B}); \zeta_{A, J, B, J}; (1_{A\diamond B}\star \mu_J) = {\rho_\star}_{A\diamond B}
\end{align*}
\end{fullwidth}
\end{enumerate}

A duoidal category \cX\ is {\bf normal} if $\iota\c I\rarr J$ is invertible and {\bf strong} if all its structure maps are isomorphisms. 
\end{definition}

\begin{definition}\cite[Def 6.5]{Aguiar_Mahajan_2010}\label{def:braided_duoidal_cat}
A duoidal category \cX\ is {\bf braided} if $(\cX, \diamond, I)$ and $(\cX, \star, J)$ are braided monoidal categories, $(J, \mu_{J},\iota)$ is a braided $\diamond$-monoid and $(I, \Delta_{I}, \iota)$ is a cobraided $\star$-comonoid, and the interchange maps interact coherently with the braidings:
\begin{equation}\begin{gathered}\label{cc:interchange_braiding}\tag{\bf\footnotesize{BDUO}}
\xymatrixrowsep{1.75pc}\xymatrixcolsep{3.75pc}\xymatrix{
(A\star B)\diamond (C\star D)\ar[r]^-{\zeta_{A, B, C, D}}\ar[d]_-{{\sigma_\star}_{A,B} \diamond {\sigma_\star}_{C,D}} & (A\diamond C)\star (B\diamond D)\ar[d]^-{{\sigma_\star}_{A\diamond C, B\diamond D}} \\
(B\star A)\diamond (D\star C)\ar[r]_-{\zeta_{B, A, D, C}} & (B\diamond D)\star (A\diamond C)
}
\end{gathered}\end{equation}
\vspace{-\baselineskip}
\begin{align*}
& \zeta_{A, B, C, D} ; {\sigma_\star}_{A\diamond C, B\diamond D} = ({\sigma_\star}_{A,B} \diamond {\sigma_\star}_{C,D}); \zeta_{B, A, D, C}\\ 
& \zeta_{A, B, C, D} ; ({\sigma_\diamond}_{A,C} \star {\sigma_\diamond}_{B,D}) = {\sigma_\diamond}_{A\star B, C\star D};\zeta_{C,D,A,B}
\end{align*}
\end{definition}

\begin{example}\label{ex:duoidal_cat}
Every braided monoidal category $(\cX, \os, I)$ induces a duoidal category $(\cX, \os, I, \os, I)$ \cite[Prop 6.10]{Aguiar_Mahajan_2010}, with structure maps given by unitors, an identity map, and the canonical flip:
\begin{align*}
&\Delta_I = {u_\os}_I \c I\rarr I\os I \qquad \mu_I = {u_\os}^{-1}_I \c I\os I \rarr I \qquad \iota = 1_I \c I\rarr I \\
&\zeta_{A,B,C,B} = \tau_{A,B,C,D} \c (A\os B)\os (C\os D)\rarr  (A\os C)\os (B\os D)
\end{align*}
        
Here we state Joyal and Street's result about braided monoidal categories \cite{Joyal_Street_1993}, in the terminology of duoidal categories:
\begin{proposition}\cite[Prop 6.11]{Aguiar_Mahajan_2010}\label{prop:strong_duoidal}
Consider a strong duoidal category $\cX$, then $(\cX, \diamond, I)$ and $(\cX, \star, J)$ are isomorphic braided monoidal categories and the interchange arises from the canonical flips. 
\end{proposition}

\end{example}

\subsection{Duoidal Functors and Transformations}  \hfill\

\vspace{0.5\baselineskip}
Aguiar and Mahajan define different types of functors between duoidal categories: bilax duoidal functors and double (co)lax duoidal functors. We shall only need the former.

\begin{definition}\cite[Def 6.50]{Aguiar_Mahajan_2010}
Let $\cX$ and $\cY$ denote duoidal categories. A {\bf bilax duoidal functor} $(F, p_{I},p_\diamond, q_{J},q_\star) \c \cX\rarr \cY$ is a functor $F\c \cX\rarr \cY$ which is
\begin{itemize}
	\item a lax monoidal functor $(F, p_{I},p_\diamond)\c (\cX,\diamond,I)\rarr(\cY,\diamond, I)$, with 
    \[ p_{I}\c I\rarr F(I)\qquad {p_\diamond}_{A, B}\c F(A)\diamond F(B)\rarr F(A\diamond B)\]
	\item a colax monoidal functors $(F, q_{J},q_\star)\c (\cX,\star,J)\rarr(\cY,\star, J)$, with 
    \[ q_{J}\c F(J)\rarr J \qquad {q_\star}_{A, B}\c F(A\star B)\rarr F(A)\star F(B)\]
\end{itemize}
such that the monoidal structure maps interact coherently with
\begin{enumerate}[{\bf \footnotesize (DF.1)}]
    \item the interchange maps, and
\begin{fullwidth}
\begin{equation*}\begin{gathered}\label{cc:duoidalfunctor_interchange}
\xymatrixrowsep{1.75pc}\xymatrixcolsep{4.75pc}\xymatrix{
F(A\star B)\diamond F(C\star D)\ar[r]^-{{q_\star}_{A,B}\diamond {q_\star}_{C,D}}\ar[d]_-{{p_\diamond}_{A\star B, C\star D}} & (F(A)\star F(B)) \diamond (F(C)\star F(D))\ar[d]^-{\zeta_{F(A), F(B), F(C), F(D)}} \\
F((A\star B)\diamond (C\star D))\ar[d]_-{F(\zeta_{A,B,C,D})} & (F(A)\diamond F(C))\star (F(B)\diamond F(D))\ar[d]^-{{p_\diamond}_{A,C}\star {p_\diamond}_{B,D}} \\
F((A\diamond C)\star (B\diamond D))\ar[r]_-{{q_\star}_{A\diamond C, B\diamond D}} & F(A\diamond C)\star F(B\diamond D)
}
\end{gathered}\end{equation*}
\end{fullwidth}

    \item the unit (co)monoid structure maps.
\begin{equation*}\begin{gathered}\label{cc:duoidalfunctor_unit}
\xymatrixrowsep{1.75pc}\xymatrixcolsep{1.75pc}\xymatrix{
I\ar[r]^-{p_{I}}\ar[d]_-{\Delta_{I}} & F(I)\ar[r]^-{F(\Delta_{I})} & F(I\star I)\ar[d]^-{{q_\star}_{I,I}} & I\ar[r]^-{p_I}\ar[d]_-{\iota} & F(I)\ar[d]^-{F(\iota)}\\
I\star I\ar[rr]_-{p_{I}\star p_{I}} && F(I)\star F(I) & J & F(J)\ar[l]^-{q_J}\\
F(J) \diamond F(J)\ar[r]^-{{p_\diamond}_{J,J}}\ar[d]_-{q_{J}\diamond q_{J}} & F(J\diamond J)\ar[r]^-{F(\mu_{J})} & F(J)\ar[d]^-{q_J} \\
J\diamond J\ar[rr]_-{q_J} && J
}
\end{gathered}\end{equation*}
\end{enumerate}
\end{definition}

\begin{definition}\cite[Def 6.51]{Aguiar_Mahajan_2010}
A {\bf duoidal transformation} between bilax duoidal functors $\alpha\c (F, p_{I}^{F},p_\diamond^{F}, q_{J}^{F},q_\star^{F}) \Rarr (G, p_{I}^{G},p_\diamond^{G}, q_{J}^{G},q_\star^{G})$ is a natural transformation $\alpha\c F\rarr G$ such that 
\begin{itemize}
	\item $\alpha\c (F, p_{I}^{F},p_\diamond^{F})\Rarr (G,p_{I}^{G}, p_\diamond^{G})$ is a monoidal transformation, and
	\item $\alpha\c (F, q_{J}^{F}, q_\star^{F})\Rarr (G, q_{J}^{G},q_\star^{G})$ is a monoidal transformation.
\end{itemize}
\end{definition}

\begin{proposition}\cite[Prop 6.52]{Aguiar_Mahajan_2010}\label{prop:2cat_duoidalcats}
There is a 2-category of duoidal categories, bilax duoidal functors and duoidal transformations, denoted \DUO, closed under finite products. 
\end{proposition}

Recalling that a monoid in a monoidal category \cX\ can be equivalently defined as a monoidal functor from the terminal category to \cX, Aguiar and Mahajan define different types of monoid-like structures in a duoidal category based on their different functors: duoidal bimonoids and duoidal double monoids. Once more, we only introduce the former, which will be key to the linearly distributive Fox theorem.

\begin{definition}\cite[Def 6.25]{Aguiar_Mahajan_2010}
A {\bf duoidal bimonoid} in a duoidal category \cX\ is a quintuple $\langle A, \nabla_A, u_A, \Delta_A, e_A\rangle$ consisting an object $A$ equipped with four morphisms
\[ \nabla_A\c A\diamond A \rarr A \qquad u_A\c I\rarr A \qquad \Delta_A\c A\rarr A\star A \qquad e_A\c A\rarr J\]
such that $\langle A, \nabla_A, u_A\rangle$ is a $\diamond$-monoid, $\langle A, \Delta_A, e_A\rangle$ is a $\star$-comonoid, and the two structures are compatible:
\begin{equation}\begin{gathered}\label{cc:duoidal_bimonoid}\tag{\bf\footnotesize{DBM}}
\xymatrixrowsep{1.75pc}\xymatrixcolsep{1.75pc}\xymatrix{
A\diamond A\ar[r]^-{\nabla_A}\ar[d]_-{\Delta_A \diamond \Delta_A} & A\ar[r]^-{\Delta_A} & A\star A  & I\ar[r]^-{u_A}\ar[rd]_-{\iota} & A\ar[d]^-{e_A}\\
(A\star A)\diamond (A\star A)\ar[rr]_-{\zeta_{A,A,A,A}} && (A\diamond A)\star (A\diamond A)\ar[u]_-{\nabla_A \star \nabla_A} & & J \\
A\diamond A\ar[r]^-{\nabla_A} \ar[d]_-{e_A\diamond e_A}& A\ar[d]^-{e_A} & I\ar[r]^-{\Delta_I}\ar[d]_-{u_A} & I\star I\ar[d]^-{u_A\star u_A}\\
J \diamond J \ar[r]_-{\mu_J}& J & A\ar[r]_-{\Delta_A} & A\star A
}
\end{gathered}\end{equation}
A {\bf morphism of duoidal bimonoids} is a morphism of the underlying $\diamond$-monoid and $\star$-comonoid.
\end{definition}

\begin{remark}\label{rem:duoidal_bimonoid}
A duoidal bimonoid $\langle A, \nabla_A, u_A, \Delta_A, e_A\rangle$ in \cX\ can be equivalently defined as a bilax duoidal functor from the terminal category to \cX\ (with slight abuse of notation) 
\[ (A, u_A, \nabla_A, e_A, \Delta_A)\c 1\rarr \cX, \qquad *\mapsto A, \qquad 1_*\mapsto1_A\]
\vspace{-\baselineskip}
\begin{align*}
    & u_A = p_{I}\c I\rarr A && \nabla_A = {p_\diamond}_{A, B}\c A\diamond A\rarr A \\
    & e_A = q_{J}\c A\rarr J && \Delta_A = {q_\star}_{A, B}\c A\rarr A\star A
\end{align*}
and a morphism of duoidal bimonoids $f\c \langle A, \nabla_A, u_A, \Delta_A, e_A\rangle\rarr \langle B, \nabla_B, u_B, \Delta_B, e_B\rangle$ as a duoidal transformation between such functors
\begin{equation*}
\xymatrixrowsep{1.75pc}\xymatrixcolsep{7.75pc}\xymatrix{ 
1\rtwocell<5>^{(A, u_A, \nabla_A, e_A, \Delta_A)}_{(B, u_B, \nabla_B, e_B, \Delta_B)}{\,\,f} &\cX
}
\end{equation*}
\end{remark}

With the perspective outlined in the previous remark and Proposition \ref{prop:2cat_duoidalcats}, we get that: 

\begin{proposition}\cite[Cor 6.53]{Aguiar_Mahajan_2010}\label{prop:duoidal_functor_preserve}
A bilax duoidal functor preserves bimonoids and morphisms of bimonoids. Moreover, the component maps of a duoidal transformation are morphisms of the preserved duoidal bimonoids. 

\end{proposition}
\begin{proof}
Consider a bilax duoidal functor $(F, p_{I},p_\diamond, q_{J},q_\star) \c \cX\rarr \cY$ and a duoidal bimonoid $\langle A, \nabla_{A}, u_{A}, \Delta_{A}, e_{A}\rangle$ in \cX. Viewing the bimonoid as a bilax duoidal functor and post-composing with $F$, we obtain the composite bilax duoidal functor
\[ 1\xrightarrow{(A, u_A, \nabla_A, e_A, \Delta_A)}\cX \xrightarrow{(F, p_{I},p_\diamond, q_{J},q_\star)} \cY\]
As a bilax duoidal functor $1\rarr\cY$, the composite determines a duoidal bimonoid structure on $F(A)$ in \cY, given by the structure maps
\begin{align*}
&\nabla_{F(A)} = F(A)\diamond F(A) \xrightarrow{{p_\diamond}_{A,A}} F(A\diamond A) \xrightarrow{F(\nabla_{A})} F(A) \\
&u_{F(A)} =  I \xrightarrow{p_I} F(I) \xrightarrow{F(u_{A})} F(A) \\
&\Delta_{F(A)} = F(A) \xrightarrow{F(\Delta_A)} F(A\star A) \xrightarrow{{q_\star}_{A,A}} F(A)\star F(A) \\
&e_{F(A)} =F(A) \xrightarrow{F(e_{A})} F(J) \xrightarrow{q_J} J
\end{align*} 

Consider now $f\c \langle A, \nabla_{A}, u_{A}, \Delta_{A}, e_{A}\rangle\rarr \langle B, \nabla_{A}, u_{B}, \Delta_{B}, e_{B}\rangle$, a morphism of duoidal bimonoids. Viewing this as a duoidal transformation and whiskering on the right with $F$, we obtain the composite duoidal transformation
\begin{equation*}
\xymatrixrowsep{1.75pc}\xymatrixcolsep{5.75pc}\xymatrix{ 
1\rtwocell<5>^{(A, u_A, \nabla_A, e_A, \Delta_A)}_{(B, u_B, \nabla_B, e_B, \Delta_B)}{\,\,f} &\cX\ar[r]^-{(F, p_{I},p_\diamond, q_{J},q_\star)} & \cY 
}
\end{equation*}
As a duoidal transformation between bilax duoidal functors $1\rarr \cY$, $F(f)$ determines a morphism of duoidal bimonoids $F(A)\rarr F(B)$. 

Finally, consider $\alpha\c (F, p_{I}^{F},p_\diamond^{F}, q_{J}^{F},q_\star^{F}) \Rarr (G, p_{I}^{G},p_\diamond^{G}, q_{J}^{G},q_\star^{G})$. By whiskering on the left with the bilax duoidal functor $A$, we obtain the composite duoidal transformation 
\begin{equation*}
\xymatrixrowsep{1.75pc}\xymatrixcolsep{7.75pc}\xymatrix{ 
1\ar[r]^-{(A, u_A, \nabla_A, e_A \Delta_A)} & \cX\rtwocell<5>^{(F, p_{I}^{F},p_\diamond^{F}, q_{J}^{F},q_\star^{F})}_{(G, p_{I}^{G},p_\diamond^{G}, q_{J}^{G},q_\star^{G})}{\,\,\alpha} &\cY
}
\end{equation*}
As a duoidal transformation between bilax duoidal functors $1\rarr \cY$, $\alpha_A$ determines a morphism of duoidal bimonoids $F(A)\rarr G(A)$.
\end{proof}

This is precisely the reason bilax duoidal functors and transformations will be considered and a central reason why the linearly distributive Fox theorem will hold.

\subsection{Symmetry} \hfill\

\vspace{0.5\baselineskip}
The notion of symmetric duoidal categories has not yet appeared in the literature as such, although it is a trivial extension, since the motivation for duoidal structures is often braided monoidal categories. Symmetry is however essential to the Fox theorem and therefore will be a key component of this work. 

\begin{definition}
A braided duoidal category $\cX$ is {\bf symmetric} if $(\cX, \diamond, I)$ and $(\cX, \star, J)$ are symmetric.
\end{definition}

Given a symmetric duoidal category, there are two canonical flips $\tau^\diamond$ and $\tau^\star$ \mbox{\bf\footnotesize{(CF)}}. These maps ensure that much of the structures of a duoidal category become themselves bilax duoidal functors and transformations. 

\begin{lemma}
The following diagrams and equations commute in any symmetric duoidal category $\cX$:
\begin{equation}\begin{gathered}\label{cc:delta_nabla_canonicalflip}
\xymatrixrowsep{1.75pc}\xymatrixcolsep{2.75pc}\xymatrix{
(J\diamond J)\diamond(J\diamond J)\ar[r]^-{\mu_J\diamond \mu_J}\ar[d]_-{\tau^\diamond_{J,J,J,J}} & J\diamond J\ar[dd]^-{\mu_J} &I\ar[r]^-{\Delta_I}\ar[dd]_-{\Delta_I} & I\star I \ar[d]^-{\Delta_I \star \Delta_I} \\
(J\diamond J)\diamond(J\diamond J)\ar[d]_-{\mu_J} & & & (I\star I)\star (I\star I)\ar[d]^-{\tau^\star_{I,I,I,I}} \\
J\diamond J\ar[r]_-{\mu_J} & J & I\star I\ar[r]_-{\Delta_I\star \Delta_I} & (I\star I)\star (I\star I)
}
\end{gathered}\end{equation}
\vspace{-\baselineskip}
\begin{align}
\begin{split}\label{cc:medial_canonicalflip}
&\tau^\diamond_{A\star A', B\star B', C\star C', D\star D'};(\zeta_{A,A',C,C'}\diamond\zeta_{B,B',D,D'}); \zeta_{A\diamond C, A'\diamond C', B\diamond D, B'\diamond D'} \\
&\qquad=(\zeta_{A,A',B,B'}\diamond \zeta_{C,C', D, D'});\zeta_{A\diamond B, A'\diamond B', C\diamond D, C'\diamond D'}; (\tau^\diamond_{A,B,C,D}\star \tau^\diamond_{A',B',C',D'}) \\[0.5pc]
&\zeta_{A\star A', B\star B', C\star C', D\star C'}; (\zeta_{A,A',C,C'}\star \zeta_{B,B',D,D'}); \tau^\star_{A\diamond C, A'\diamond C', B\diamond D, B'\diamond D'}\\ 
&\qquad=(\tau^\star_{A,A',B,B'}\diamond \tau^\star_{C, C', D, D'}); \zeta_{A\star B, A'\star B', C\star D, C'\star D'}; (\zeta_{A,B,C,D}\star \zeta_{A',B',C',D'}) 
\end{split}
\end{align}
\end{lemma}

\begin{proposition}\label{prop:monoidal_products_bilax_duoidal_functor}
Given a symmetric duoidal category \cX, the monoidal units and products are canonically bilax duoidal functors:
\[ (I, 1_I, {\rho_\diamond}_I^{-1}, \iota, \Delta_I), (J, \iota, \mu_J, 1_J, {\rho_\star}_J)\c 1\rarr\cX\]
\[(\diamond, {\rho_\diamond}_I, \tau^\diamond, \mu_J, \zeta),(\star, \Delta_I, \zeta, {\rho_\star}_J^{-1}, \tau^\star)\c \cX\times \cX\rarr\cX \]
\end{proposition}
\begin{proof}
It is immediate that $(I, 1_I, {\rho_\diamond}_I^{-1})\c 1\rarr\cX$ is a $\diamond$-lax monoidal functor. Furthermore, $(I,\iota, \Delta_I)\c 1\rarr\cX$ is a $\star$-colax monoidal functor as it is a $\star$-comonoid, by definition of a duoidal category. 

As for $\diamond$, since $(\cX, \diamond, I)$ is a SMC, it is well-known that $(\diamond, {\rho_\diamond}_I, \tau^\diamond)\c \cX\times\cX\rarr\cX$ is a $\diamond$-lax monoidal functor. Moreover, as discussed in \cite[Prop 6.4]{Aguiar_Mahajan_2010}, the second equation in {\bf \footnotesize (DUO.1)}, and the second and fourth equations in {\bf \footnotesize (DUO.2)} ensure precisely that $(\diamond, \mu_J, \zeta)\c \cX\times\cX\rarr\cX$ is a $\star$-colax monoidal functor. It remains to show {\bf \footnotesize (DF.1-2)}. 

Now, {\bf \footnotesize (DF.1)} is the first equation in \eqref{cc:medial_canonicalflip}. The first two diagrams in {\bf \footnotesize (DF.2)} hold by the following commuting diagrams.
\begin{equation*}\begin{gathered}\label{cc:duoidalfunctor_unit}
\xymatrixrowsep{1.75pc}\xymatrixcolsep{2.75pc}\xymatrix{
I\ar[r]^-{{\rho_\diamond}_I}\ar[dd]_-{\Delta_{I}}\ar@{}[rd]|{(\nat)} & I\diamond I \ar[r]^-{\Delta_{I}\diamond \Delta_I}\ar[d]^-{\Delta_I\diamond 1_I} & (I\star I)\diamond (I\star I)\ar[dd]^-{{\zeta}_{I,I, I, I}}  & I\ar@{}[rd]|{(\nat)}\ar[dd]_-{\iota}\ar[rr]^-{{\rho_\diamond}_I} && I\diamond I \ar[dd]^-{\iota\diamond\iota}\ar[ld]_-{\iota\diamond 1_I}\\ 
& (I\star I)\diamond I\ar@{}[r]|{\bf \footnotesize (DUO.2)} & & & J\diamond I\ar@{}[rd]|{\bf \footnotesize (M.2)} & \\
I\star I\ar[rr]_-{p_{I}\star p_{I}}\ar[ru]_-{{\rho_\diamond}_{I\star I}} && (I\diamond I)\star (I\diamond I) & J\ar[ru]^-{{\rho_\diamond}_J} && J\diamond J\ar[ll]^-{\mu_J}\\
}
\end{gathered}\end{equation*}
The final diagram in {\bf \footnotesize (DF.2)} is the first diagram in \eqref{cc:delta_nabla_canonicalflip}. A similar argument holds for $J$ and $\star$ being bilax duoidal functors.
\end{proof}

\begin{proposition}\label{prop:monoidal_isomorphisms_duoidal_trans}
Given a symmetric duoidal category \cX, the associators, the unitors and the braidings are duoidal transformations between the respective canonical bilax duoidal functors. 
\end{proposition}
\begin{proof}
Consider the right $\diamond$-unitor ${\rho_\diamond}_{A}\c A\rarr A\diamond I$. It is a $\diamond$-monoidal transformation by a trivial diagram and the first diagram in \eqref{cc:canonical_flip_unitor}, and it is a $\star$-monoidal transformation by {\bf \footnotesize (M.2)} for the $\diamond$-monoid $(J, \mu_J, \iota)$ and the third diagram in {\bf \footnotesize (DUO.2)}. As such, it is a duoidal transformation.

The $\diamond$-braiding ${\sigma_\diamond}_{A,B}\c A\diamond B\rarr B\diamond A$ is $\diamond$-monoidal by compatibility with left and right unitor, as well as the equality of right and left unitors on the units, and by \eqref{cc:canonical_flip_braiding}. It is a $\star$-monoidal transformation by {\bf \footnotesize (CM)} for $(J, \mu_J, \iota)$ as a commutative $\diamond$-monoid, and {\bf \footnotesize (BDUO)}. Given that duoidal transformations compose, the left $\diamond$-unitor ${\lambda_\diamond}_{A}\c A\rarr I\diamond A$ is also a duoidal transformation.

Finally, the $\diamond$-associator ${\alpha_\diamond}_{A,B,C}\c (A\diamond B)\diamond C\rarr A\diamond (B\diamond C)$ is a $\diamond$-monoidal transformation by {\bf \footnotesize (MC.2)}, as well as the equality of right and left unitors on the units, and \eqref{cc:canonical_flip_associator}. It is a $\star$-monoidal transformation by {\bf \footnotesize (M.1)} for the $\diamond$-monoid $(J, \mu_J, \iota)$ and the first equation in {\bf \footnotesize (DUO.1)}.  Similarly, the $\star$-monoidal structure maps are duoidal transformations. 
\end{proof}

This idea extends to the structure maps of bilax duoidal functors. 

\begin{proposition}\label{prop:duoidal_functor_maps_duoidal_trans}
Consider a bilax duoidal functor $(F, p_{I},p_\diamond, q_{J},q_\star)\c\cX\rarr\cY$ between symmetric duoidal categories, then its structure maps are duoidal transformations between the respective canonical bilax duoidal functors. 
\end{proposition}
\begin{proof}
It is well-known that $p_I\c I\rarr F(I)$ is a $\diamond$-monoidal transformation, and the first and second diagrams in {\bf \footnotesize (DF.2)} mean precisely that it is a $\star$-monoidal transformation. 

Consider now ${p_\diamond}_{A,B}\c F(A)\diamond F(B)\rarr F(A\diamond B)$. It is a $\diamond$-monoidal transformation by naturality, unitality of lax monoidal functors {\bf \footnotesize (MF.2)} and its interaction with the canonical flip \eqref{cc:strong_symmetric_monoidal_flip}.  {\bf \footnotesize (DF.1)} and the third diagram of {\bf \footnotesize (DF.2)} precisely ensure it is a $\star$-monoidal transformation. A similar argument works for $q_J$ and $q_\star$. 
\end{proof}

\subsection{Cartesian and Cocartesian Duoidal Categories}  \hfill\

\vspace{0.5\baselineskip}
As previously stated, duoidal structure arises canonically whenever a monoidal category has finite products or finite coproducts. We detail the construction here.

\begin{proposition}\cite[Ex 6.19]{Aguiar_Mahajan_2010}\label{prop:cartesian_duoidal}
    \begin{enumerate}
        \item Consider a monoidal category $(\cX,\os, I)$ with finite products, then $(\cX, \os, I, \times, \bone)$ is a duoidal category with structure maps
    \[ \Delta_I = \langle 1_I, 1_I\rangle\c I\rarr I\times I \qquad \mu_{\bone} =  t_{\bone\os\bone}\c \bone\os\bone\rarr \bone \qquad  \iota = t_{I}\c I\rarr \bone\]
    \[\zeta_{A,B,C,D} = \langle \pi^0_{A,B} \os \pi^0_{C,D}, \pi^1_{A,B} \os \pi^1_{C,D}\rangle\c (A\times B)\os(C\times D)\rarr (A\os C)\times (B\os D) \]

    \item Consider a monoidal category $(\cX,\os, I)$ with finite coproducts, then $(\cX, +, \bzero, \os, I)$ is a duoidal category with structure maps
    \[  \Delta_{\bzero} = b_{\bzero\os \bzero}\c \bzero\rarr \bzero\os \bzero \qquad \mu_I =  [1_I, 1_I]\c I+I\rarr I \qquad \iota = b_{I}\c \bzero\rarr I\]
    \[\zeta_{A,B,C,D} = [\iota^0_{A,C} \os \iota^0_{B,D}, \iota^1_{A,C} \os \iota^1_{B,D}]\c (A\os B)+(C\os D)\rarr (A+C)\os (B+D) \]

    \item Consider now a category \cX\ with finite products and coproducts, then $(\cX, +, \bzero, \times, \bone)$ is of course a symmetric duoidal category with structure maps
    \[  \Delta_{\bzero} = b_{\bzero\times \bzero} = \langle 1_\bzero, 1_\bzero\rangle\c \bzero\rarr \bzero\times \bzero \qquad  \mu_{\bone} = t_{\bone+\bone} = [1_\bone, 1_\bone]\c\bone+\bone\rarr \bone\]
    \[ \iota = t_{\bzero}=b_{\bone}\c \bzero\rarr\bone \]
    \[ \zeta_{A,B,C,D} = \langle \pi^0_{A,B} + \pi^0_{C,D}, \pi^1_{A,B} +\pi^1_{C,D}\rangle = [\iota^0_{A,C} \times \iota^0_{B,D}, \iota^1_{A,C} \times \iota^1_{B,D}]\c\]
    \[ (A\times B)+(C\times D)\rarr (A+C)\times (B+D)\]
    \end{enumerate}
\end{proposition}

This idea extends to monoidal functors between monoidal categories with finite products or finite coproducts as follows. 

\begin{lemma}\label{lem:functor_comon_if_cart}
Given cartesian categories \cX\ and \cY, any functor $F\c \cX\rarr\cY$ becomes a symmetric colax monoidal functor with structure maps
\[ q_{\bone}= t_{F(\bone)}\c F(\bone)\rarr \bone\quad\quad {q_{\times}}_{A,B} = \langle F(\pi^0_{A,B}), F(\pi^1_{A, B})\rangle\c F(A\times B)\rarr F(A)\times F(B) \]
Further, any natural transformation $\alpha\c F\Rarr G$ is monoidal between the induced symmetric colax monoidal functors. 

Dually, given cocoartesian categories \cX\ and \cY, any functor $G\c \cX\rarr\cY$ becomes a symmetric lax monoidal functor with structure maps
\[ p_{\bzero}= b_{F(\bzero)}\c \bzero\rarr F(\bzero)\quad\quad {p_+}_{A,B} = [F(\iota^0_{A,B}), F(\iota^1_{A, B})]\c F(A)+ F(B)\rarr F(A+B) \]
Any natural transformation $\alpha\c F\Rarr G$ is monoidal between the induced symmetric lax monoidal functors.
\end{lemma}

\begin{proposition}\cite[Ex 6.68]{Aguiar_Mahajan_2010}\label{prop:cartesian_duoidal_functor}
\begin{enumerate}
    \item Consider duoidal categories $(\cX, \os, I, \times, \bone)$ and $(\cY, \os, I, \times, \bone)$ whose second monoidal structures are cartesian and a lax monoidal functor $(F, p_I, p_\os)\c (\cX, \os, I)\rarr (\cY, \os, I)$. By the above Lemma, $(F, q_{\bone}, q_\times)\c (\cX, \times, \bone)\rarr(\cY, \times, \bone)$ is a colax monoidal functor. This further determines a bilax duoidal functor 
\[(F, p_I, p_\os, q_{\bone}, q_\times)\c (\cX, \os, I, \times, \bone)\rarr (\cY, \os, I, \times,\bone)\]
    \item Consider duoidal categories $(\cX,+, \bzero, \os, I)$ and $(\cY,+, \bzero, \os, I)$ whose first monoidal structures are cocartesian. Then, any colax monoidal functor $(F, q_I, q_\os)\c(\cX, \os, I)\rarr (\cY, \os, I)$ is equally a lax monoidal functor $(F, p_{\bzero}, p_+)\c (\cX, +, \bzero)\rarr(\cY, +, \bzero)$ by the above Lemma. This determines a bilax duoidal functor 
\[(F, p_{\bzero}, p_+, q_I, q_\os)\c (\cX,+, \bzero, \os, I)\rarr(\cY,+, \bzero, \os, I)\]
    \item If we now consider a categories \cX\ and \cY\ with finite products and coproducts, then every functor $F\c\cX\rarr\cY$ is canonically a bilax duoidal functor.
\end{enumerate}
\end{proposition}


\section{Medial Linearly Distributive Categories}\label{sec:MLDC}

We can immediately provide a characterization of CLDCs $(\bX, \times, \bone, +, \bzero)$ by recognizing that $(\bX, \times, \bone)$ is cartesian and $(\bX, +, \bzero)$ is cocartesian, and therefore, applying Corollary \ref{cor:char_cartesian_cat} to the former and applying its dual to the latter. 

\begin{proposition}\label{prop:char_CLDC}
A SLDC \bX\ is cartesian if and only if there are natural transformations
\[ \Delta_{A}\c A\rarr A\ot A \quad\quad e_{A}\c A\rarr \top \quad\quad \nabla_{A}\c A\op A\rarr A \quad\quad u_{A}:\bot\rarr A \]
such that $\langle A, \Delta_{A}, e_{A}\rangle$ is a cocommutative $\ot$-comonoid, $\langle A, \nabla_{A}, u_{A}\rangle$ is a commutative $\op$-monoid, and 
\[\Delta_{A\ot B} = (\Delta_{A}\ot\Delta_{B}); \tau^\ot_{A, A, B, B} \qquad e_{A\ot B} = (e_{A}\ot e_{B}); {u^{R}_{\ot}}_{\top}^{-1} \]
\[\nabla_{A\op B} = \tau^\op_{A, B, A, B};(\nabla_{A}\op\nabla_{B}) \qquad u_{A\op B} = {u^R_\op}^{-1}_{\bot};(u_{A}\op u_{B}) \]
\[\Delta_{\top} = {u^{R}_{\ot}}_{\top} \qquad e_{\top} = 1_{\top} \qquad  \nabla_{\bot} = {u^{R}_{\op}}_{\bot} \qquad u_{\bot} = 1_{\bot} \]
\end{proposition}

The natural transformations of course have component morphisms for all objects. Therefore, there will be maps \[ e_\bot, u_\top\c \bot\rarr\top \qquad u_{\bot\ot\bot}, \Delta_\bot\c \bot\rarr \bot\ot \bot \qquad e_{\top\op\top}, \nabla_\top\c \top\op\top\rarr\top\] and, given $A, B \in \bX$, there will be a maps \[ \Delta_{A\op B}\c A\op B \rarr (A\op B)\ot (A\op B) \qquad \nabla_{A\ot B}\c (A\ot B)\op (A\ot B)\rarr A\ot B \] 
The conditions of Proposition \ref{prop:char_CLDC} give us additional information about these maps. Indeed, by naturality of $e$ and equation $e_\top = 1_\top$, we see that  \[ e_\bot = u_\top \qquad e_{\top\op\top} = \nabla_\top \qquad e_{A\op B} = (e_A\op e_B) ; \nabla_\top \]
Similarly, by naturality of $u$ and equation $u_\bot = 1_\bot$, 
\[ u_\top = e_\bot \qquad u_{\bot\ot\bot} = \Delta_\bot \qquad u_{A\ot B} = \Delta_\bot; (u_A \ot u_B)\]
Furthermore, we can determine a formula for the diagonal morphism $\Delta_{A\ot B}$:
\begin{equation*}\begin{gathered}
\xymatrixrowsep{2.75pc}\xymatrixcolsep{8.75pc}\xymatrix@L=0.35pc{
A\op B\ar[r]^-{\Delta_A\op \Delta_B}\ar[d]_-{\Delta_{A\op B}}\ar@{}[rd]|{(\nat)} & (A\ot A)\op (B\ot B)\ar[d]^-{\Delta_{(A\ot A)\op (B\ot B)}} \\
(A\op B)\ot (A\op B)\ar@{}[rd]|{(\mathrm{\bf CM.2})}\ar@{=}[d]\ar[r]^-{(\Delta_A\op\Delta_B)\ot (\Delta_A\op\Delta_B)} & ((A\ot A)\op (B\ot B)) \ot ((A\ot A)\op (B\ot B))\ar[d]^-{\hbox{$\scriptsize\begin{array}{@{}c@{}}((1_A\ot e_A)\op (1_B\ot e_B))\ot {}\\ ((e_A\ot 1_A)\op (e_B\ot 1_B))\end{array}$}} \\
(A\op B)\ot (A\op B) & ((A\ot\top)\op (B\ot\top))\ot ((\top\ot A)\op (\top\ot B))\ar[l]^-{({u^R_\ot}^{-1}_A \op {u^R_\op}^{-1}_B)\ot ({u^L_\ot}^{-1}_A \op {u^L_\op}^{-1}_B)}
}\end{gathered}\end{equation*}
Therefore, \[ \Delta_{A\op B} = A\op B \xrightarrow{\Delta_A\op \Delta_B} (A\ot A)\op (B\ot B) \xrightarrow{\mu^0_{A,A,B,B}} (A\op B)\ot (A\op B)\]
where $\mu^0_{A,B,C,D}\c (A\ot B)\op (C\ot D) \rarr (A\op C)\ot (B\op D)$ is the natural transformation
\[ \Delta_{(A\ot B)\op(C\ot D)}; ( ( (1_A\ot e_B); {u^R_\ot}^{-1}_A \op (1_C\ot e_D); {u^R_\ot}^{-1}_C) \ot ((e_A\ot 1_B); {u^L_\ot}^{-1}_B \op (e_C\ot 1_D); {u^L_\ot}^{-1}_D))\]
Similarly:
\begin{equation*}\begin{gathered}
\xymatrixrowsep{2.75pc}\xymatrixcolsep{8.75pc}\xymatrix@L=0.35pc{
((A\op \bot)\ot(B\op\bot))\op((\bot\op A)\ot (\bot\op B))\ar[d]_-{\hbox{$\scriptsize\begin{array}{@{}c@{}}((1_{A}\op u_{A})\ot (1_{B}\op u_{B})) {}\\ ((u_{A}\op 1_{ A})\ot (u_{B}\op 1_{B}))\end{array}$}}\ar@{}[rd]|{(\mathrm{\bf M.2})} & (A\ot B)\op (A\ot B)\ar[l]_-{({u^R_\op}_A^{-1}\ot{u^R_\op}_B^{-1})\op ({u^L_\op}_A^{-1} \ot {u^L_\op}_B^{-1})}\ar@{=}[d] \\
(A\op A)\ot (B\op B))\op ((A\op A)\ot (B\op B))\ar@{}[rd]|{(\nat)}\ar[d]_-{\nabla_{(A\op A)\ot (B\op B)}}\ar[r]_-{(\nabla_A\ot \nabla_B)\op(\nabla_A\ot \nabla_B)}& (A\ot B)\op (A\ot B)\ar[d]^-{\nabla_{A\ot B}} \\
(A\op A)\ot (B\op B))\ar[r]_-{\nabla_A\ot\nabla_B} & A\ot B
}\end{gathered}\end{equation*}
This determines that \[ \nabla_{A\ot B} = (A\ot B)\op (A\ot B) \xrightarrow{\mu^1_{A,B,A,B}} (A\op A)\ot (B\op B) \xrightarrow{\nabla_A\ot \nabla_B} A\ot B \]
where $\mu^1_{A,B,C,D}\c (A\ot B)\op (C\ot D) \rarr (A\op C)\ot (B\op D)$ is the natural transformation
\[ ( ( {u^R_\op}^{-1}_A;(1_A\ot u_B)\op {u^R_\op}^{-1}_C;(1_C\ot u_D)) \ot ( {u^L_\op}^{-1}_B; (u_A\ot 1_B) \op {u^L_\op}^{-1}_D;(u_C\ot 1_D)); \nabla_{(A\op C)\ot (B\op D)}\]
We can further demonstrate that $\mu^0_{A,B,C,D} = \mu^1_{A,B,C,D}$, though the detailed computations are omitted here due to their length.

\begin{corollary}\label{cor:char_CLDC_naive}
Consider a SLDC \bX\ satisfying the conditions of Proposition \ref{prop:char_CLDC}, then there are maps
\[ m \c \bot\rarr\top \qquad \Delta_\bot\c \bot\rarr\bot\ot\bot \qquad \nabla_\top\c \top\op\top\rarr\top\]
and a natural transformation 
\[ \mu_{A, B, C, D}\c (A\ot B)\op(C\ot D)\rarr(A\op C)\ot(B\op D) \]
such that 
\[ \Delta_{A\op B} = (\Delta_A \op\Delta_B) ; \mu_{A,A,B,B} \qquad e_{A\op B} = (e_A\op e_B);\nabla_\bot \]
\[ \nabla_{A\ot B} = \mu_{A,B,A,B}; (\nabla_A \ot \nabla_B) \qquad u_{A\ot B} = \Delta_\top; (u_A\ot u_B) \]
\[ e_\bot = u_\top = m \]
\end{corollary}

The maps indicated above will be central to our main definition and will be the focus of the following section.

\subsection{Medial Rule and Main Definition}  \hfill\

\vspace{0.5\baselineskip}
Considering SLDCs provide the categorical semantics for \MLLp, it is important to examine the logical significance of the following map: \[ (A\ot B)\op (C\ot D) \rarr (A\op C)\ot (B\op D)\] Viewed as a sequent, it is not derivable within \MLLp, as exemplified by many SLDCs which do not have such a natural transformation. However, if we consider, for example, a two-sided sequent calculus presentation of classical logic, taking $\ot$ to be classical conjunction $\wedge$ and $\op$ to be classical disjunction $\vee$, then $(A\wedge B) \vee (C\wedge D) \vdash (A\vee C)\wedge (B\vee D)$ is a provable sequent. In fact, it appears prominently as an inference rule in multiple deductive systems based on deep inference. Deep inference, more specifically the calculus of structures, is an alternative logic formalism to Gentzen's sequent calculus, which was developed by Guglielmi when studying a non-commutative extension of linear logic known as pomset logic \cite{Guglielmi_2007}. 

It was first introduced as an inference rule by Br\"{u}nnler and Tiu \cite{Brunnler_Tiu_2001} for their local system of classical logic SKS and was named the {\bf medial rule}. In this context, non-locality refers to the idea that certain rules require a global view of formulae of potentially unbounded size, like in the case of contraction within Gentzen’s system for classical logic \cite{Brunnler_2003_1}. Locality in SKS is achieved by ensuring that non-local rules, including contraction, are presented in atomic form. In order to ensure the general form of contraction remains admissible in SKS, the medial rule was introduced. The medial rule also appears in Stra\ss burger’s local system of linear logic SLLS, to ensure that contraction and cocontraction remain atomic \cite{Strassburger_2002}, and similarly in Tiu’s local system for intuitionistic logic SISgq \cite{Tiu_2006}. Following its introduction in deep inference, the medial rule has equally been studied as a term-rewriting rule, both in isolation \cite{Strassburger_2007_2} and alongside the linear distributivity and mix rules \cite{Das_Strassburger_2016, Bruscoli_Strassburger_2017}.

In categorical logic, the medial rule has been examined in the study of categorical semantics for classical logic. Unlike intuitionistic logic, which is modeled by cartesian closed categories, and linear logic, which is modeled by $*$-autonomous or LDCs, the appropriate categorical framework for classical logic remains an active area of research. Indeed, the naive definition one might consider, cartesian $*$-autonomous categories, leads only to posetal categories due to the so-called Joyal’s paradox \cite{Lambek_Scott_1988}. The medial rule plays a key role in Stra\ss burger’s Boolean categories \cite{Strassburger_2007_1} and in Lamarche’s $*$-autonomous categories with finitary medial and the absorption law \cite{Lamarche_2007}. In both cases, the authors began with $*$-autonomous categories and introduced medial maps, alongside the $\bot$-contraction map $\bot\rarr \bot\ot \bot$ and the $\top$-cocontraction map $\top\op\top\rarr\top$, in order to construct categorical models of classical logic. Furthermore, in both constructions, the nullary mix map $\bot\rarr\top$ is present. In particular, Lamarche’s categories, which will be shown to be related to our main definition in Section \ref{sec:adding_negation}, were a major inspiration for the present work.

From a categorical perspective, the medial maps can be recognized as an instance of the well-known interchange law of duoidal categories, taking $\diamond = \op$, $J=\bot$, $\star = \op$ and $I = \bot$. When combined with  $\bot$-contraction, $\top$-cocontraction and the nullary mix map, we get all the necessary structure for a duoidal category.

To establish the linearly distributive Fox theorem, we cannot begin with the 2-category of SLDCs alone. Instead, we must consider SLDCs that additionally include the medial maps, including nullary mix, $\bot$-contraction and $\top$-cocontraction. This necessity aligns with the role of the medial rule in SLLS, where it enables the atomic contraction and cocontraction rules to extend to all formulae. Therefore, we must first define what it means for a LDC to have coherent medial maps. 

\begin{definition}
A {\bf medial linearly distributive category}, or a MLDC, $(\bX,\ot,\top,\op,\bot)$ is a category \bX\ equipped with 
\begin{itemize}
\item a {\bf tensor} monoidal structure $(\bX, \ot, \top)$,
\[ {\alpha_\ot}_{A,B,C}\c (A\ot B)\ot C\rarr A\ot (B\ot C) \quad  {u_\ot^R}_{A}\c A\rarr A\ot \top \quad {u_\ot^L}_{A}\c A\rarr \top\ot A \]
\item a {\bf par} monoidal structure $(\bX, \op, \bot)$,
\[{\alpha_\op}_{A,B,C}\c A\op(B\op C)\rarr (A\op B)\op C \quad {u_\op^R}_{A}\c A\op\bot\rarr A \quad {u_\op^L}_{A}\c\bot\op A\rarr A\] 
\item {\bf $\bot$-contraction}, {\bf $\top$-cocontraction} and {\bf nullary mix} maps,
\[ \Delta_{\bot}\c\bot\rarr\bot\ot\bot \qquad \nabla_{\top}\c \top\op\top\rarr\top \qquad m\c\bot\rarr\top\]
\item a {\bf medial} natural transformation, and 
\[ \mu_{A, B, C, D}\c (A\ot B)\op(C\ot D)\rarr(A\op C)\ot(B\op D) \]
\item left and right {\bf linear distributor} natural transformations,
\begin{align*}
\delta^L_{A,B,C}\c A\ot(B\op C)\rarr (A\ot B)\op C && \delta^R_{A,B,C}\c (A\op B)\ot C\rarr A\op (B\ot C)
\end{align*}
\end{itemize}
such that 
\begin{enumerate}[(i)]
\item $(\bX,\ot,\top,\op,\bot)$ is a mix LDC,
\item $(\bX,\op,\bot,\ot,\top)$ is a duoidal category, and
\item the {\bf medial maps} interact coherently with the linear distributors:
\begin{fullwidth}
\begin{equation}\begin{gathered}\label{cc:medial_lindist}\tag{\bf\footnotesize{MLDC.1}}
\xymatrixrowsep{1.75pc}\xymatrixcolsep{5.75pc}\xymatrix{
((A\ot B)\op (C\ot D))\ot X\ar[r]^-{\mu_{A,B,C,D} \ot 1_{X}}\ar[d]_-{\delta^R_{A\ot B, C\ot D, X}} & ((A\op C)\ot(B\op D))\ot X\ar[d]^-{{\alpha_\ot}_{A\op C, B\op D, X}} \\
(A\ot B)\op ((C\ot D)\ot X)\ar[d]_-{1_{A\ot B} \op {\alpha_\ot}_{C, D, X}} &(A\op C)\ot((B\op D)\ot X) \ar[d]^-{1_{A\op C} \ot \delta^R_{B,D,X}} \\
(A\ot B)\op (C\ot(D\ot X))\ar[r]_-{\mu_{A, B, C, D\ot X}} & (A\op C) \ot (B\op (D\ot X))
}
\end{gathered}\end{equation}
\vspace{-\baselineskip}
\begin{align*}
(\mu_{A,B,C,D} \ot 1_{X}); &{\alpha_\ot}_{A\op C, B\op D, X}; (1_{A\op C} \ot \delta^R_{B,D,X})= \delta^R_{A\ot B, C\ot D, X} ; (1_{A\ot B} \op {\alpha_\ot}_{C, D, X}); \mu_{A, B, C, D\ot X} \\ 
(1_{X} \ot \mu_{A,B,C,D}); &{\alpha_\ot}^{-1}_{X, A\op C, B\op D} ; (\delta^L_{X,A,C}\ot 1_{B\op D})= \delta^L_{X,A\ot B, C\ot D} ; ({\alpha_\ot}^{-1}_{X, A, B}\op 1_{C\ot D}); \mu_{X\ot A, B, C, D} \\
(\delta^R_{X,A,B} \op 1_{C\ot D}); &{\alpha_\op}^{-1}_{X, A\ot B, C\ot D}; (1_{X} \op \mu_{A,B,C,D}) = \mu_{X\op A, B, C, D}; ({\alpha_\op}^{-1}_{X, A, C}\ot 1_{B\op D}); \delta^R_{X, A\op C, B\op D} \\
(1_{A\ot B}\op \delta^L_{C,D,X}); &{\alpha_\op}_{A\ot B, C\ot D, X}; (\mu_{A,B,C,D}\op 1_{X})  = \mu_{A,B,C,D\op X}; (1_{A\op C}\ot {\alpha_\op}_{B,D,X}); \delta^L_{A\op C, B\op D, X}
\end{align*}
\end{fullwidth}
\end{enumerate}
\end{definition}
In the above definition, $(\bX,\op,\bot,\ot,\top)$ being a duoidal category corresponds to:
\begin{enumerate}[{\bf \footnotesize (MLDC.1)}]
  \setcounter{enumi}{1}
    \item $(\bot, \Delta_{\bot}, m)$ is a $\ot$-comonoid and $(\top, \nabla_{\top}, m)$ is a $\op$-monoid,
    
\begin{fullwidth}
\begin{equation*}\label{cc:bot_top_comonoid_monoid}
\begin{aligned}
&\qquad\Delta_{\bot}; (1_{\bot}\ot \Delta_{\bot})=\Delta_{\bot}; (\Delta_{\bot}\ot 1_{\bot});{\alpha_\ot}_{\bot,\bot,\bot} \\
&\qquad(1_{\top}\op\nabla_{\top});\nabla_{\top} = {\alpha_\op}_{\top,\top,\top};(\nabla_{\top} \op 1_{\top});\nabla_{\top}\\[0.5pc]
&\Delta_{\bot};(1_{\bot}\ot m) = {u^R_\ot}_{\bot} \qquad\qquad (1_{\top}\op m);\nabla_{\top} = {u^R_\op}_{\top} \\
&\Delta_{\bot};(m\ot 1_{\bot}) = {u^L_\ot}_{\bot} \qquad\qquad (m \op1_{\top});\nabla_{\top} = {u^L_\op}_{\top} 
\end{aligned}
\end{equation*}
\end{fullwidth}

    \item the medial maps interact coherently with the associators, and 
\begin{fullwidth}
\begin{align*}
\begin{split}
&\mu_{A\ot B, C, D\ot E, F} ;(\mu_{A, B, D, E}\ot 1_{C\op F}); {\alpha_\ot}_{A\op D, B\op E, C\op F} \\
&= ({\alpha_\ot}_{A,B,C}\op{\alpha_\ot}_{D,E,F}); \mu_{A, B\ot C, D,E\ot F} ; (1_{A\op D}\ot \mu_{B,C,E,F}) \\[0.5pc]
&{\alpha_\op}_{A\ot B, C\ot D, E\ot F}; (\mu_{A,B,C,D}\op 1_{E\ot F}); \mu_{A\op C, B\op D, E, F} \\
&= (1_{A\ot B}\op \mu_{C,D,E,F}); \mu_{A,B,C\op E, D\op F}; ({\alpha_\op}_{A,C,E}\ot {\alpha_\op}_{B,D,F})\\
\end{split}
\end{align*}
\end{fullwidth}

    \item the medial maps interact coherently with the unitors.
\begin{fullwidth}
\begin{align*}
&(\Delta_{\bot}\op 1_{A\ot B}); \mu_{\bot,\bot,A,B};({u_\op^L}_{A} \ot {u_\op^L}_{B}) = {u_\op^L}_{A\ot B} \\
&(1_{A\ot B}\op\Delta_{\bot});\mu_{A,B,\bot,\bot}; ({u_\op^R}_{A} \ot {u_\op^R}_{B}) = {u_\op^R}_{A\ot B}\\
&({u_\ot^L}_{A} \op {u_\ot^L}_{B}); \mu_{\top,A,\top,B}; (\nabla_{\top}\ot 1_{A\op B}) = {u_\ot^L}_{A\op B} \\
&({u_\ot^R}_{A} \op {u_\ot^R}_{B}); \mu_{A,\top,B,\top}; (1_{A\op B}\ot \nabla_{\top}) = {u_\ot^R}_{A\op B} 
\end{align*}
\end{fullwidth}
\end{enumerate}

As discussed, MLDCs provide a suitable framework for categories that are both independently and coherently linearly distributive and duoidal. We use the term independently to emphasize that, while certain duoidal categories can inherit a linearly distributors from their interchange maps (this construction will be reviewed in Section \ref{sec:properties}), this is not what we are currently discussing. We are concerned with categories where the duoidal and linearly distributive structures do not induce one another, yet nonetheless interact coherently.  This perspective will be particularly useful as it allows us to leverage the well-developed theory of duoidal categories.

Alternatively, MLDCs can be viewed as potentially providing categorical semantics for \MLLp\ with the medial rule (\MLLp+Medial), an extension of \MLLp+MIX. This variant may be incompatible with standard sequent calculus, as the medial rule is rather unconventional: its premise is a disjunction of formulae, while its conclusion is a conjunction. However, \MLLp+Medial can likely be described within the calculus of structures, where linear distributivity and the medial rule already play central roles in several systems within the deep inference paradigm. This question was similarly posed by Lamarche for \MLL+Medial at the end of \cite[Sec 2]{Lamarche_2007}. Overall, the logical study of \MLLp+Medial and \MLL+Medial remains an open area for future research. 

\subsection{Symmetry}\label{sec:symmetry}  \hfill\

\vspace{0.5\baselineskip}
If we revisit Proposition \ref{prop:char_CLDC}, which provides a preliminary characterization of CLDCs via the traditional Fox theorem, we observe that, for any objects $A$ and $B$ in a CLDC, the diagonal $\Delta_{A\ot B}$ and the multiplication $\nabla_{A\op B}$ are constructed using the canonical flip maps $\tau^\ot$ and $\tau^\op$ \mbox{\bf\footnotesize{(CF)}}. By the same reasoning that led us to define MLDCs to ensure the presence of the necessary medial maps, we now define symmetric MLDCs to guarantee the existence of these canonical flip maps.

\begin{definition}\label{def:symm_MLDC}
A MLDC $\bX$ is {\bf symmetric}, or a SMLDC, if it is symmetric as a LDC and symmetric as a duoidal category, in other words:
\begin{enumerate}[(i)]
    \item $(\bX,\ot,\top)$ and $(\bX,\op,\bot)$ are SMCs,
    \item $(\bot,\Delta_\bot, m)$ is a cocommutative $\ot$-comonoid, $(\top, \nabla_\top, m)$ is a commutative $\op$-monoid, and 
    \item \mbox{\bf\footnotesize{(SLDC)}} and \mbox{\bf\footnotesize{(BDUO)}} hold.
\end{enumerate}
\end{definition}

There is however an equivalent definition for SMLDCs which does not introduce the nullary mix map $m\c\bot\rarr\top$, nor does it require \mbox{\bf\footnotesize{(MIX)}} to hold. Indeed, in the presence of the braidings, we can define it in terms of the other structure maps and get the mix structure as a consequence.

It is well-known that given a duoidal category, the nullary mix map is defined by the rest of the structure.

\begin{proposition}\cite[Prop 6.9]{Aguiar_Mahajan_2010}\label{prop:iota_defined_by_zeta}
The structure map $m\c\bot\rarr\top$ in a duoidal category $(\bX, \op,\bot,\ot,\top)$ is equal to following equivalent composites
\begin{equation}\label{eqn:nullary_mix_map_composites}
\begin{gathered}\xymatrixrowsep{1.75pc}\xymatrixcolsep{1.75pc}\xymatrix@L=0.5pc{
&&(\bot\ot\top)\op(\top\ot\bot)\ar[r]^-{\mu_{\bot,\top,\top\,\bot}} & (\bot\op\top)\ot(\top\op\bot)\ar[rd]^(0.65){{u^L_\op}_\top \ot {u^R_\op}_\top} \\
\bot\ar[r]^-{{u_\op}_\bot^{-1}} & \bot\op\bot\ar[ru]^(0.35){{u^R_\ot}_\bot \op {u^L_\ot}_\bot}\ar[rd]_(0.35){{u^L_\ot}_\bot \op {u^R_\ot}_\bot} & & & \top\ot\top\ar[r]^-{{u_\ot}^{-1}_\top} & \top \\
&& (\top\ot\bot)\op(\bot\ot\top)\ar[r]_*+<0.25em>{_{\mu_{\top,\bot,\bot,\top}}} & (\top\op\bot)\ot(\bot\op\top)\ar[ur]_(0.65){{u^R_\op}_\top \ot {u^L_\op}_\top}
}
\end{gathered}\end{equation}
\end{proposition}

The proof that both these composites are equivalent, and themselves equal to $m$ follows from the conditions that $(\top, \nabla_\top, m)$ is a $\op$-monoid and $(\bot, \Delta_\bot, m)$ is a $\ot$-comonoid. Therefore, the definition of an arbitrary duoidal category must include the nullary mix map. 

However, in the presence of braidings, it can be excluded. It was first noticed by Lamarche \cite[Prop 2.11]{Lamarche_2007} that, given a category \bX\ with two symmetric monoidal structures $(\bX, \ot, \top)$ and $(\bX, \op, \bot)$, and a coherent medial natural transformation $\mu_{A,B,C,D}$, i.e. a category where \mbox{\bf\footnotesize{(BDUO)}} holds, the composites in \eqref{eqn:nullary_mix_map_composites} are equal. 

Furthermore, in the presence of $\bot$-contraction $\Delta_\bot$ such that $\bot$ is a cocommutative $\ot$-cosemigroup and $\top$-cocontraction map $\nabla_\top$ such that $\top$ is commutative $\op$-semigroup, which interact coherently with the medial maps, i.e. ({\footnotesize{\bf DUO.2}}) holds, then the equivalent composites in \eqref{eqn:nullary_mix_map_composites} further endow $\bot$ and $\top$ with $\ot$-comonoid and $\op$-monoid structures respectively, as shown in \cite[Prop 2.19]{Lamarche_2007}. 

Then, as $\bot$ is endowed with a $\ot$-comonoid structure and $\top$ is endowed with a $\op$-monoid structure in a SLDC, Fuhrmann and Pym's result \cite[Thm 3.11]{Fuhrmann_Pym_2007} implies their counit and unit, which are equal in this context, induce a mix structure.  

Altogether, these observations allows us to give an alternative definition of SMLDCs which excludes the nullary mix, defining it entirely based on the rest of the structure.

\begin{proposition}\label{prop:alternative_SMLDC}
A SMLDC $\bX$ is equivalently defined as a SLDC equipped with a $\bot$-contraction map, a $\top$-cocontraction map and a medial natural transformation
\[ \Delta_{\bot}\c\bot\rarr\bot\ot\bot \qquad \nabla_{\top}\c \top\op\top\rarr\top\]
\[ \mu_{A, B, C, D}\c (A\ot B)\op(C\ot D)\rarr(A\op C)\ot(B\op D) \]
such that 
\begin{enumerate}[(i)]
\item $(\bot, \Delta_\bot)$ is a cocommutative $\ot$-cosemigroup and $(\top, \nabla_\top)$ is commutative $\op$-semigroup,
\begin{fullwidth}
\begin{align}\label{cc:bot_top_comm_semigroups}
\begin{split}
&\Delta_{\bot} = \Delta_{\bot}; {\sigma_\ot}_{\bot,\bot} \qquad \Delta_{\bot}; (1_{\bot}\ot \Delta_{\bot}) = \Delta_{\bot}; (\Delta_{\bot}\ot 1_{\bot}); {\alpha_\ot}_{\bot,\bot,\bot} \\
&\nabla_{\top} = {\sigma_\op}_{\top,\top}; \nabla_{\top} \qquad (1_{\top}\op\nabla_{\top});\nabla_{\top} = {\alpha_\op}_{\top,\top,\top};(\nabla_{\top} \op 1_{\top});\nabla_{\top} 
\end{split}
\end{align}
\end{fullwidth}

\item the medial maps interact coherently with the braidings, the associators and with the unitors,
\begin{fullwidth}
\begin{align}\label{cc:sym_medial_braiding}
\begin{split}
& \mu_{A, B, C, D} ; {\sigma_\ot}_{A\op C, B\op D} = ({\sigma_\ot}_{A,B} \op {\sigma_\ot}_{C,D}); \mu_{B, A, D, C}\\ 
& \mu_{A, B, C, D} ; ({\sigma_\op}_{A,C} \ot {\sigma_\op}_{B,D}) = {\sigma_\op}_{A\ot B, C\ot D};\mu_{C,D,A,B}
\end{split}
\end{align}

\begin{align}\label{cc:sym_medial_assoc}
\begin{split}
&\mu_{A\ot B, C, D\ot E, F} ; (\mu_{A, B, D, E}\ot 1_{C\op F}); {\alpha_\ot}_{A\op D, B\op E, C\op F} = \\
&({\alpha_\ot}_{A,B,C}\op{\alpha_\ot}_{D,E,F}); \mu_{A, B\ot C, D,E\ot F} ; (1_{A\op D}\ot \mu_{B,C,E,F}) \\[0.5pc]
&{\alpha_\op}_{A\ot B, C\ot D, E\ot F};(\mu_{A,B,C,D}\op 1_{E\ot F}); \mu_{A\op C, B\op D, E, F} \\
&= (1_{A\ot B}\op \mu_{C,D,E,F}); \mu_{A,B,C\op E, D\op F}; ({\alpha_\op}_{A,C,E}\ot {\alpha_\op}_{B,D,F})
\end{split}
\end{align}

\begin{align}\label{cc:sym_absorption_law}
\begin{split}
&(1_{A\ot B}\op\Delta_{\bot});\mu_{A,B,\bot,\bot}; ({u_\op^R}_{A} \ot {u_\op^R}_{B}) = {u_\op^R}_{A\ot B}\\
&({u_\ot^R}_{A} \op {u_\ot^R}_{B}); \mu_{A,\top,B,\top}; (1_{A\op B}\ot \nabla_{\top}) = {u_\ot^R}_{A\op B} 
\end{split}
\end{align}
\end{fullwidth}

\item the medial maps interact coherently with the linear distributors, and
\begin{fullwidth}
\begin{align}\label{cc:sym_medial_lindist}
\begin{split}
&(1_{X} \ot \mu_{A,B,C,D}); {\alpha_\ot}^{-1}_{X, A\op C, B\op D} ; (\delta^L_{X,A,C}\ot 1_{B\op D})= \\
&\delta^L_{X,A\ot B, C\ot D} ; ({\alpha_\ot}^{-1}_{X, A, B}\op 1_{C\ot D}); \mu_{X\ot A, B, C, D} \\[0.5pc]
&(1_{A\ot B}\op \delta^L_{C,D,X}); {\alpha_\op}_{A\ot B, C\ot D, X}; (\mu_{A,B,C,D}\op 1_{X})  \\
&= \mu_{A,B,C,D\op X}; (1_{A\op C}\ot {\alpha_\op}_{B,D,X}); \delta^L_{A\op C, B\op D, X}
\end{split}
\end{align} 
\end{fullwidth}

\item the medial maps, the linear distributors and the braidings interact coherently together.
\begin{fullwidth}
\begin{align}\label{cc:sym_medial_lindist_braiding}
\begin{split}
&(1_{X} \ot \mu_{A,B,C,D}); {\alpha_\ot}^{-1}_{X, A\op C, B\op D} ; ({\sigma_\ot}_{X,A\op C}\ot 1_{B\op D}); {\alpha_\ot}_{A\op C, X, B\op D}; (1_{A\op C}\ot \delta^L_{X,B,D}) =\\
&\delta^L_{X,A\ot B, C\ot D} ; ({\alpha_\ot}^{-1}_{X, A, B}\op 1_{C\ot D}); (({\sigma_\ot}_{X,A}\ot 1_B)\op 1_{C\ot D}); ({\alpha_\ot}_{A,X,B}\op 1_{C\ot D}); \mu_{A, X\ot B, C, D} \\[0.5pc]
&(\delta^L_{A,B,X}\op 1_{C\ot D}); {\alpha_\op}^{-1}_{A\ot B, X, C\ot D}; (1_{A\ot B}\op {\sigma_\op}_{X, C\ot D}); {\alpha_\op}_{A\ot B, C\ot D, X}; (\mu_{A,B,C,D}\op 1_X)= \\
&\mu_{A,B\op X, C, D}; (1_{A\op C}\ot {\alpha_\op}^{-1}_{B,X,D}); (1_{A\op C}\ot (1_B\op {\sigma_\op}_{X,D})); (1_{A\op C}\ot {\alpha_\op}_{B,D,X}); \delta^L_{A\op C, B\op D, X}
\end{split}
\end{align}  
\end{fullwidth}
\end{enumerate}

\end{proposition}
\begin{proof}
It is immediate that Definition \ref{def:symm_MLDC} implies the definition presented in the Proposition, since the former is a version of the latter with more structures and conditions, noting that \eqref{cc:sym_medial_lindist_braiding} is obtained by replacing the right linear distributor in \mbox{\bf\footnotesize{(MLDC.1)}} by the composite in \mbox{\bf\footnotesize{(SLDC)}}, and simplifying using {\bf\footnotesize{(BMC.1)}} and \eqref{cc:sym_medial_braiding}.

Now, consider a SLDC $\bX$ as outlined in the Proposition, then 
\begin{equation*}
\begin{gathered}\xymatrixrowsep{1.75pc}\xymatrixcolsep{3.75pc}\xymatrix@L=0.5pc{
&(\bot\ot\top)\op(\top\ot\bot)\ar[r]^-{\mu_{\bot,\top,\top\,\bot}}\ar[dd]^-{{\sigma_\ot}_{\bot,\top}\op {\sigma_\ot}_{\top,\bot}}\ar@{}[rdd]|{\eqref{cc:sym_medial_braiding}}  & (\bot\op\top)\ot(\top\op\bot)\ar[rd]^(0.65){{u^L_\op}_\top \ot {u^R_\op}_\top}\ar[dd]_-{{\sigma_\ot}_{\bot\op\top, \top\op\bot}} \\
\bot\op\bot\ar[ru]^(0.35){{u^R_\ot}_\bot \op {u^L_\ot}_\bot}\ar[rd]_(0.35){{u^L_\ot}_\bot \op {u^R_\ot}_\bot}\ar@{}[r]|{\qquad(\mathrm{\bf BMC.2})} & & \ar@{}[r]|{(\mathrm{\bf BMC.2})\qquad} & \top\ot\top \\
& (\top\ot\bot)\op(\bot\ot\top)\ar[r]_*+<0.25em>{_{\mu_{\top,\bot,\bot,\top}}} & (\top\op\bot)\ot(\bot\op\top)\ar[ur]_(0.65){{u^R_\op}_\top \ot {u^L_\op}_\top}
}
\end{gathered}\end{equation*}
As such, we can define $m\c\bot\rarr\top$ to be the equivalent composites in \eqref{eqn:nullary_mix_map_composites}. 

$(\bot, \Delta_\bot, m)$ is a cocommutative $\ot$-comonoid by Figure \ref{fig:bot_cocomm_ot_comonoid} in Appendix \ref{app:commuting_diagrams}. The upper composite in the diagram is equal to $\Delta_{\bot};(1_{\bot}\ot m)$, while the lower composite
\[ {u^R_\ot}_{\bot}; {u^R_\op}^{-1}_{\bot\ot\top}; (1_{\bot\ot \top}\op\Delta_{\bot});\mu_{\bot,\top,\bot,\bot}; ({u_\op^R}_{\bot} \ot {u_\op^R}_{\top}) = {u^R_\ot}_{\bot}; {u^R_\op}^{-1}_{\bot\ot\top};{u_\op^R}_{\bot\ot \top} = {u^R_\ot}_{\bot} \] by \eqref{cc:sym_absorption_law}. Thus,  $\Delta_{\bot};(1_{\bot}\ot m) = {u^R_\ot}_{\bot}$. Similarly, we can show that $m\c\bot\rarr\top$ endows $\top$ with a commutative $\op$-monoid structure. By Theorem \ref{thm:prove_mix}, this implies the SLDC is mix. 

Finally, it straightforward to prove the symmetric versions of \eqref{cc:sym_medial_assoc},  \eqref{cc:sym_absorption_law} using the braidings. The only technical computation is showing that the coherence conditions \eqref{cc:sym_medial_braiding}, \eqref{cc:sym_medial_lindist} and \eqref{cc:sym_medial_lindist_braiding} imply the right linear distributor, defined by \mbox{\bf\footnotesize{(SLDC)}}, interacts coherently with the medial maps, in other words \mbox{\bf\footnotesize{(MLDC.1)}} holds.
\end{proof}

This alternative definition is of particular use because it involves fewer coherence conditions and will be helpful to show that SMLDCs with negation are equivalent to Lamarche's definition.

Furthermore, an important consequence of the SMLDC coherence conditions is they ensure the left linear distributor interacts coherently with the duoidal structure maps.

\begin{proposition}\label{prop:interaction_canonical_flip_medial_maps_linear_dist}
The following holds in any SMLDC:
\begin{equation}\begin{gathered}\label{eqn:interaction_linear_dist_duoidal_trans_1}
\xymatrixrowsep{1.75pc}\xymatrixcolsep{2.75pc}\xymatrix{
\top\ot(\top\op\top)\ar[rr]^-{\delta^L_{\top,\top,\top}}\ar[d]_-{1_\top\ot\nabla_\top} && (\top\ot\top)\op\top\ar[d]^-{{u_\ot}^{-1}_\top\op1_\top} \\
\top\ot\top\ar[r]_-{{u_\ot}^{-1}_\top} & \top & \top\op\top\ar[l]^-{\nabla_\top}
}
\end{gathered}\end{equation}
\vspace{-\baselineskip}
\begin{align*}
& \delta^L_{\top,\top,\top};({u_\ot}^{-1}_\top\op1_\top);\nabla_\top = 1_\top\ot\nabla_\top;{u_\ot}^{-1}_\top   \\
&\Delta_\bot; (1_\bot\ot{u_\op}_\bot^{-1}); \delta^L_{\bot,\bot,\bot} = {u_\op}_\bot^{-1}; (\Delta_\bot\op 1_\bot)
\end{align*}
\begin{equation}\begin{gathered}\label{cc:interaction_canonical_flip_medial_maps_linear_dist}
\xymatrixrowsep{1.75pc}\xymatrixcolsep{5.75pc}\xymatrix{
(A\ot B)\ot((C\ot D)\op (E\ot F))\ar[r]^-{\delta^L_{A\ot B, C\ot D, E\ot F}}\ar[d]_-{1_{A\ot B}\ot \mu_{C,D,E,F}} & ((A\ot B)\ot (C\ot D))\op (E\ot F)\ar[d]^-{\tau^\ot_{A,B,C,D}\op 1_{E\ot F}} \\
(A\ot B)\ot ((C\op E)\ot(D\op F))\ar[d]_-{\tau^\ot_{A,B,C\op E, D\op F}} & ((A\ot C)\ot(B\ot D))\op (E\ot F)\ar[d]^-{\mu_{A\ot C, B\ot D, E,F} } \\
(A\ot (C\op E))\ot (B\ot (D\op F))\ar[r]_-{\delta^L_{A,C,E}\ot \delta^L_{B,D,F}} & ((A\ot C)\op E) \ot ((B\ot D)\op F)
}
\end{gathered}\end{equation}
\vspace{-\baselineskip}
\begin{align*}
&\delta^L_{A\ot B, C\ot D, E\ot F};(\tau^\ot_{A,B,C,D}\op 1_{E\ot F});\mu_{A\ot C, B\ot D, E,F} \\
&= (1_{A\ot B}\ot \mu_{C,D,E,F}); \tau^\ot_{A,B,C\op E, D\op F}; (\delta^L_{A,C,E}\ot \delta^L_{B,D,F}) \\[0.5pc]
&\mu_{A, B\op C, D, E\op F}; (1_{A\op D}\ot \tau^\op_{B,C,E,F}); \delta^L_{A\op D, B\op E, C\op F}= \\
&(\delta^L_{A,B,C}\op \delta^L_{D, E, F}); \tau^\op_{A\ot B, C, D\ot E, F}; (\mu_{A,B,D,E}\op 1_{C\op F})
\end{align*}
As such, the left linear distributor is a duoidal transformation between bilax duoidal functors.
\end{proposition}
\begin{proof}
The first equality in \eqref{eqn:interaction_linear_dist_duoidal_trans_1} holds by the following commuting diagram, with the second equality holding similarly. 
\begin{equation*}
\xymatrixrowsep{1.75pc}\xymatrixcolsep{2.75pc}\xymatrix{
\top\ot(\top\op\top)\ar[r]^-{\delta^L_{\top,\top,\top}}\ar[dd]_-{1_\top\ot\nabla_\top}\ar[rd]_-{{u^L_\ot}^{-1}_{\top\op\top}}^-{\quad(\mathrm{\bf LDC.1})} & (\top\ot\top)\op\top\ar[d]^-{{u_\ot}^{-1}_\top\op1_\top} \\
\ar@{}[rd]|{(\nat)}& \top\op\top\ar[d]^-{\nabla_\top} \\
\top\ot\top\ar[r]_-{{u_\ot}^{-1}_\top} & \top
}
\end{equation*}
The proof of the first equality in \eqref{cc:interaction_canonical_flip_medial_maps_linear_dist} is given by Figure \ref{fig:interaction_canonical_flip_medial_maps_linear_dist} in Appendix \ref{app:commuting_diagrams} and the other equality follows similarly. 

Now, the monoidal products $\ot$ and $\op$ are bilax duoidal functors by Proposition \ref{prop:monoidal_products_bilax_duoidal_functor}. As such, the left linear distributors is a natural transformation between functors which are canonically bilax duoidal:
\[ \delta^L\c (1_{\bbX}\times\op);\ot\Rarr (\ot\times1_{\bbX}); \op\]

The first equalities in \eqref{eqn:interaction_linear_dist_duoidal_trans_1} and \eqref{cc:interaction_canonical_flip_medial_maps_linear_dist} are precisely the conditions such that $\delta^L$ is a $\ot$-monoidal transformation, while the second equalities in \eqref{eqn:interaction_linear_dist_duoidal_trans_1} and \eqref{cc:interaction_canonical_flip_medial_maps_linear_dist} are precisely the conditions such that $\delta^L$ is a $\op$-monoidal transformation. Therefore, $\delta^L$ is a duoidal transformation.
\end{proof}

\subsection{Adding Negation}\label{sec:adding_negation}  \hfill\

\vspace{0.5\baselineskip}
As linear negation is central to linear logic and a key attribute of many models of \MLL, it is important we describe how the duoidal structure interacts with it. This will recover Lamarche's definition of $*$-autonomous categories with finitary (binary and nullary) medial and with the absorption law \cite{Lamarche_2007}. As such, we can say SMLDCs are the ``negation-free'' version of the former.

\begin{definition}
A SMLDC $\bX$ {\bf has negation} if it has negation as a SLDC and satisfies the following conditions
\begin{enumerate}[(i)]
    \item $\bot$-contraction $\Delta_\bot$ and $\top$-cocontraction $\nabla_\top$ are dual to one another, and 
\begin{fullwidth}
\begin{equation*}\begin{gathered}\label{cc:delta_nabla_dual}
\xymatrixrowsep{1.75pc}\xymatrixcolsep{1.75pc}\xymatrix{
(\bot\ot\bot)^\perp\ar[d]_-{(\Delta_\bot)^\perp}\ar[r]^-{\sim} & \top\op\top\ar[d]^-{\nabla_\top} \\
\bot^\perp\ar[r]_-{\sim} & \top
}
\end{gathered}\end{equation*}
\end{fullwidth}
    \item the medial maps are self-dual.
\begin{fullwidth}
\begin{equation*}\begin{gathered}\label{cc:medial_selfdual}
\xymatrixrowsep{1.75pc}\xymatrixcolsep{2.75pc}\xymatrix{
(A^\perp \ot B^\perp) \op (C^\perp \op D^\perp)\ar[d]_-{\mu_{A^\perp, B^\perp, C^\perp, D^\perp}} \ar[r]^-{\sim} & ((D\op C) \ot (B\op A))^\perp \ar[d]^-{\mu_{D,B,C,A}^\perp} \\
(A^\perp \op C^\perp) \ot (B^\perp \op D^\perp) \ar[r]_-{\sim}&  ((D\ot B)\op (C\ot A))^\perp\\
}
\end{gathered}\end{equation*}
\end{fullwidth}
\end{enumerate}   
\end{definition}

\begin{lemma}
A SMLDC with negation $\bX$ has a self-dual nullary mix map:
\begin{equation*}\begin{gathered}\label{cc:m_dual}
\xymatrixrowsep{1.75pc}\xymatrixcolsep{1.75pc}\xymatrix{
\bot\ar[d]_-{m}\ar[r]^-{\sim} & \top^\perp\ar[d]^-{m^\perp} \\
\top\ar[r]_-{\sim} & \bot^\perp
}
\end{gathered}\end{equation*}
\end{lemma}
\begin{proof}
This follows from the definition of $m\c \bot\rarr\top$ as a composite of unitors and the medial map in Proposition \ref{prop:iota_defined_by_zeta}, and the fact that the medial maps are self-dual themselves.  
\end{proof}

By Theorem \ref{thm:SLDC_starauto}, which states the correspondence between SLDCs with negation and $*$-autonomous categories, and by comparing Lamarche's definition to the alternative definition of a SMLDC in Proposition \ref{prop:alternative_SMLDC}, it is immediate that: 

\begin{proposition}
The notions of SMLDCs with negation and of $*$-autonomous categories with finitary medial, and with the absorption law, as defined in \cite{Lamarche_2007}, coincide.
\end{proposition}

We will however be particularly invested in SMLDCs which do not have a notion of negation. This is due to the ultimate goal of providing a 2-functor which maps SMLDCs to CLDCs as part of the linearly distributive Fox theorem. If this construction preserves the negation of SMLDCs, then by Joyal's paradox, the resulting CLDCs would necessarily be posetal. This will be discussed in greater detail in Section \ref{sec:medial_bimonoids}. 

\subsection{Examples}\label{sec:examples}  \hfill\ 

\vspace{0.5\baselineskip}
The most important class of MLDC examples is of course the CLDCs, as the former were defined precisely to be an appropriate generalization of the latter. The canonical duoidal structure of a CLDC was discussed by this author with co-author Lemay \cite[Sec 3.5]{Kudzman-Blais_Lemay_2025}. It is precisely the one for cartesian and cocartesian categories, as described in Proposition \ref{prop:cartesian_duoidal}. In particular, \cite[Prop 3.12]{Kudzman-Blais_Lemay_2025} details how this duoidal structure interacts with the linear distributors, which implies that:

\begin{theorem}\label{thm:CLDC_is_SMLDC}
Every CLDC \bX\ is canonically a SMLDC, with structure maps
\[ \Delta_{\bzero} = b_{\bzero\times \bzero} = \langle 1_\bzero, 1_\bzero\rangle\c \bzero\rarr \bzero\times \bzero \qquad  \nabla_{\bone} = t_{\bone+\bone} = [1_\bone, 1_\bone]\c\bone+\bone\rarr \bone\]
\[ m = t_{\bzero}=b_{\bone}\c \bzero\rarr\bone \]
\[ \mu_{A,B,C,D} = \langle \pi^0_{A,B} + \pi^0_{C,D}, \pi^1_{A,B} +\pi^1_{C,D}\rangle = [\iota^0_{A,C} \times \iota^0_{B,D}, \iota^1_{A,C} \times \iota^1_{B,D}]\c\]
\[ (A\times B)+(C\times D)\rarr (A+C)\times (B+D)\]
\end{theorem}

Another immediate class of examples are the braided monoidal categories, as every monoidal category is a degenerate LDC and every braided monoidal category is a strong duoidal category.

\begin{proposition}\label{prop:braided_monoidal_MLDC}
Consider a braided monoidal category $(\cX, \os, I)$, then it is a MLDC $(\cX, \os, I, \os, I)$ with structure maps
\[ \nabla_I = {u_\os}^{-1}_I \c I\os I \rarr I \qquad\qquad \Delta_I = {u_\os}_I \c I\rarr I\os I \qquad \qquad m = 1_I\c I\rarr I\]
\[ \mu_{A,B,C,D} = \tau^\os_{A,B,C,D}\c (A\os B)\os (C\os D) \rarr (A\os C)\os (B\os D)\]
\[ \delta^L_{A,B,C} = {\alpha_\os^{-1}}_{A,B,C}\c A\os (B\os C) \rarr (A\os B)\os C\] \[\delta^R_{A,B,C} = {\alpha_\os}_{A,B,C}\c (A\os B)\os C \rarr A\os (B\os C) \]
\end{proposition}

Our final source of examples stems from generalizing $Q$-coherences to the linearly distributive context, as first introduced by Lamarche \cite{Lamarche_1995} to capture Girard's coherence spaces and Ehrhard's hypercoherences, and later investigated in the context of $*$-autonomous categories with finitary medial, and with the absorption law \cite{Lamarche_2007}.  

Firstly, let $(\bP, \ot, \top, \op, \bot)$ be a small posetal SMLDC. In other words, $\bP$ is a poset with two commutative, associative and monotone binary operations $\ot$ and $\op$, equipped with units $\top$ and $\bot$ respectively, such that the following inequalities hold $\forall a, b, c, d\in \bP$:
\[ \top\op\top\leq \top \qquad \bot\leq \bot\ot\bot \qquad \bot\leq \top \]
\[ (a\ot b)\op (c\ot d) \leq (a\op c)\ot (b\op d) \]
\[ (a\op b)\ot c \leq a\op (b\ot c) \qquad a\ot (b\op c) \leq (a\ot b)\op c\]

While this may seem at first glance to be an involved definition, we can quickly see that every bounded distributive lattice is such a poset with $\ot = \wedge$ and $\op = \vee$. Indeed, as discussed in \cite{Kudzman-Blais_Lemay_2025}, a small CLDC is posetal if and only if it is a bounded distributive lattice, and a posetal CLDC is a posetal SMLDC. 

\begin{definition}\cite[Def 4.5]{Lamarche_2007} 
A {\bf \bP-coherence} $A = (|A|, \rho_A)$ consists of a poset \mbox{$(|A|, \sqsubseteq)$} and a symmetric monotone function $\rho_A\c |A|\times|A|\rarr \bP$. A {\bf \bP-coherence map} $f\c A\rarr B$ is a relation $f\c |A|\nrightarrow |B|$ which is 
\begin{enumerate}[(i)]
    \item down-closed in the source: $(a,b)\in f \wedge a' \sqsubseteq a \implies (a',b)\in f$,
    \item up-closed in the target: $(a,b)\in f \wedge b\sqsubseteq b' \implies (a,b')\in f$, and 
    \item compatible with $\rho$: $(a,b)\in f \wedge (a',b')\in f \Rarr \rho_A(a,a')\leq \rho_B(b,b')$
\end{enumerate}
\end{definition}

\begin{definition}\cite[Def 4.5]{Lamarche_2007}\label{def:PCOH}
Let \PCOH\ denote the category of \bP-coherences and \bP-coherence maps. Define additional structure on \PCOH\ as follows.
\begin{itemize}
    \item Tensor product: $A\ot B = (|A|\times |B|, \rho_{A\ot B})$ where $\rho_{A\ot B}\c (|A|\times |B|)\times ((|A|\times |B|)\rarr \bP$ is given by 
    \[ \rho_{A\ot B}((a,b),(a',b')) = \rho_A(a,a')\ot \rho_B(b,b')\]
    \item Tensor unit: $(\{*\}, \rho_\top)$ where $\rho_\top(*,*) = \top$
    \item Par product: $A\op B = (|A|\times |B|, \rho_{A\op B})$ where $\rho_{A\op B}\c (|A|\times |B|)\times ((|A|\times |B|)\rarr \bP$ is given by 
    \[ \rho_{A\op B}((a,b),(a',b')) = \rho_A(a,a')\op \rho_B(b,b')\]
    \item Par unit: $(\{*\}, \rho_\bot)$ where $\rho_\bot(*,*) = \bot$
     \item Medial maps: $\mu_{A,B,C,D}\c (A\ot B)\op (C\ot D)\rarr (A\op C)\ot (B\op D)$ is given by
    \[ (((a,b),(c,d)), ((a',c'),(b',d')))\in \mu_{A,B,C,D} \iff a\sqsubseteq a' \wedge b\sqsubseteq b' \wedge c\sqsubseteq c' \wedge d\sqsubseteq d'\]
    \item Left linear distributor: $\delta^L_{A,B,C}\c A\ot (B\op C)\rarr (A\ot B)\op C$ is given by 
    \[ (((a,(b,c)), ((a',b'),c'))\in \delta^L_{A,B,C} \iff a\sqsubseteq a' \wedge b\sqsubseteq b' \wedge c\sqsubseteq c'\]
\end{itemize}
\end{definition}

\begin{proposition}
\PCOH\ is a SMLDC.
\end{proposition}
As the above construction is a trivial generalization of Lamarche's work, we shall not include the proof in this paper.

\subsection{Properties}\label{sec:properties} \hfill\

\vspace{0.5\baselineskip}
A key aspect of the investigation into CLDCs in \cite{Kudzman-Blais_Lemay_2025} was the development of new collapse results for CLDCs. In particular, a CLDC has invertible linear distributors if and only if it is isomix if and only if it is compact \cite[Thm 5.12]{Kudzman-Blais_Lemay_2025}. In this subsection, we shall see that these results generalize to the MLDC setting.

If we consider a LDC with invertible linear distributors, the par product must be equivalent to the $\bot$-shifted tensor, as stated in Proposition \ref{prop:shift_tensor}. In the presence of medial structure, we can say even more.

\begin{proposition}
Consider a MLDC $\bX$ with invertible linear distributors, then its linearly distributive structure is compact.    
\end{proposition}
\begin{proof}
Consider a MLDC $\bX$ whose linear distributors $\delta^L_{A,B,C}$ and $\delta^R_{A,B,C}$ have inverses. By Proposition \ref{prop:shift_tensor} and the details of the proof in \cite{Cockett_Seely_1997_LDC}, $\top\op\top$ is the tensor inverse of $\bot$ via isomorphisms
\[ s^L= \bot\ot(\top\op\top)\xrightarrow{\delta^L_{\bot,\top,\top}} (\bot\ot\top)\op\top \xrightarrow{{u^R_\ot}^{-1}_\bot\op1_\top} \bot\op\top \xrightarrow{{u^L_\op}_\top} \top \]  
\[s^R=(\top\op\top)\ot\bot\xrightarrow{\delta^R_{\top, \top,\bot}} \top\op(\top\ot\bot) \xrightarrow{1_\top\op{u^L_\ot}^{-1}_\bot} \top\op\bot \xrightarrow{{u^R_\op}_\top} \top \]
This implies that the MLDC is isomix as $m\c \bot\rarr\top$ has an inverse 
\[ m^{-1} = \top \xrightarrow{{s^L}^{-1}} \bot\ot(\top\op\top) \xrightarrow{1_\bot\ot \nabla_\top} \bot\ot \top \xrightarrow{{u^R_\ot}_\bot^{-1}} \bot\]
$m^{-1}; m = 1_\top$ holds by the following commuting diagram
\begin{equation*}\begin{gathered}
\xymatrixrowsep{2.25pc}\xymatrixcolsep{5.75pc}\xymatrix{
\top\ar[rr]^-{{u^L_\op}^{-1}_\top}\ar@{=}[dddd]\ar@{}[rd]|{(\mathrm{\bf MLDC.2})} && \bot\op \top\ar[r]^-{{u^R_\ot}_\bot}\ar[ld]_-{m\op 1_\top}\ar@{}[d]|{(\nat)} & (\bot\ot\top)\op\top\ar[dd]^-{{\delta^L}_{\bot,\top,\top}^{-1}}\ar[dl]_-{(m\ot 1_\top)\op 1_\top\quad} \\
& \top\op\top\ar[dddl]_-{\nabla_\top}\ar[r]^-{{u_\ot}_\top \op 1_\top}\ar[rd]_-{{u^L_\ot}_{\top\op\top}}^{\qquad(\mathrm{\bf LDC.1})}& (\top\ot\top)\op \top\ar[d]^-{{\delta^L}_{\top,\top,\top}^{-1}}\ar@{}[rd]|{(\nat)} & \\
&\ar@{}[rd]|{(\nat)} & \top\ot(\top\op\top)\ar[d]^-{1_\top\ot \nabla_\top}\ar@{}[rd]|{(\exch)}  & \bot\ot(\top\op\top)\ar[dd]^-{1_\bot\ot \nabla_\top}\ar[l]_-{m\ot 1_{\top\op\top}}\\
&&\top\ot \top\ar[lld]_-{{u_\ot}_\top^{-1}}\ar@{}[d]|{(\nat)} &\\
\top && \bot\ar[ll]^-{m} & \bot\ot\top\ar[l]^-{{u^R_\ot}_\bot^{-1}}\ar[ul]_-{m\ot 1_\top}
}
\end{gathered}\end{equation*}
Similarly for $m; m^{-1} = 1_\bot$. As the linear distributors and the nullary mix map are isomorphisms, $\mix_{A,B}$ is an isomorphism.
\end{proof}

A similar result can be proved for isomix MLDCs. Indeed, given a normal duoidal category, there is a known construction providing a linearly distributive structure: let $(\cX, \diamond, I, \star, J)$ be a duoidal category where $\iota\c I\rarr J$ is an isomorphism. We can define the following linear distributors, as in \cite{Garner_LopezFranco_2016}:
\begin{align*}
&\partial^L_{A,B,C}= A\diamond (B\star C) \cong (A\star J)\diamond (B\star C) \xrightarrow{\zeta_{A,J,B,C}} (A\diamond B)\star (J\diamond C) \cong (A\diamond B)\star C \\
&\partial^R_{A,B,C} = (A\star B)\diamond C \cong (A\star B)\diamond (J\star C) \xrightarrow{\zeta_{A,B,J,C}} (A\diamond J) \star (B\diamond C) \cong A\star (B\diamond C)
\end{align*}
Then, $(\cX, \diamond, I, \star, J)$ is an isomix LDC \cite{Spivak_Srinivasan_2024}. Consider now an isomix MLDC $(\bX, \ot, \top,\op,\bot)$ with linear distributors
\[ \delta^L_{A,B,C}\c A\ot (B\op C)\rarr (A\ot B)\op C \qquad\qquad \delta^R_{A,B,C}\c (A\op B)\ot C\rarr A\op (B\op C)\]
Then, it is in particular a normal duoidal category $(\bX, \op, \bot,\ot,\top)$ and thus an isomix LDC with linear distributors 
\[ \partial^L_{A,B,C}\c A\op(B\ot C) \rarr (A\op B)\ot C\ \qquad\qquad \partial^R_{A,B,C}\c (A\ot B)\op C\rarr 
A\ot (B\op C) \]
Notice the direction of the latter linear distributors is opposite to the former. It turns out they are inverses of each other.

\begin{proposition}\label{prop:isomix_MLDC}
Consider an isomix MLDC $\bX$, then its linearly distributive structure is compact. 
\end{proposition}
\begin{proof}
Consider an isomix MLDC $\bX$ and define the natural transformations $\partial^L$ and $\partial^R$ as follows.
\begin{align*}
\partial^L_{A,B,C}= &A\op (B\ot C) \xrightarrow{{u^R_\ot}_A\op 1_{B\ot C}} (A\ot \top)\ot (B\ot C) \xrightarrow{\mu_{A,\top,B,C}} (A\op B)\ot (\top\op C) \\
&\xrightarrow{1_{A\op B}\ot (m^{-1}\op 1_C)} (A\op B)\ot (\bot \op C) \xrightarrow{1_{A\op B}\ot {u^L_\op}_C} (A\op B)\ot C \\
\partial^R_{A,B,C} =& (A\ot B)\op C \xrightarrow{1_{A\ot B}\op {u^L_\ot}_C} (A\ot B)\op (\top\ot C) \xrightarrow{\mu_{A,B,\top,C}} (A\op \top) \ot (B\op C)\\
&\xrightarrow{(1_A\op m^{-1})\ot 1_{B\op C}} (A\op \bot)\ot (B\op C) \xrightarrow{{u^R_\op}_A\ot 1_{B\op C}} A\ot (B\op C)
\end{align*}
Then, $\partial^R_{A,B,C}$ is the inverse of $\delta^L_{A,B,C}$ and $\partial^L_{A,B,C}$ is the inverse of $\delta^R_{A,B,C}$. The proof that $\delta^L_{A,B,C}; \partial^R_{A,B,C} = 1_{A\ot (B\op C)}$ follows from Figure \ref{fig:partial_R_inverse_delta_L} in Appendix \ref{app:commuting_diagrams}. The other equalities follow similarly. As the linear distributors and the nullary mix map are isomorphisms, $\mix_{A,B}$ is an isomorphism.
\end{proof}

Let us now consider the symmetric context. Recall that in a SMLDC, we have the canonical flip maps $\tau^\ot$ and $\tau^\op$ \mbox{\bf\footnotesize{(CF)}}. Proposition \ref{prop:interaction_canonical_flip_medial_maps_linear_dist} allows us to show that the medial maps are in some sense the canonical flip maps ``modulo'' the mix maps, as was shown to be the case for CLDCs \cite[Prop 3.14]{Kudzman-Blais_Lemay_2025}.

\begin{proposition}\label{prop:interchange_canonicalflip_mix}
Given a SMLDC, the following diagrams commute
\begin{equation}\label{diag:interchange_canonicalflip_mix}
\begin{gathered}\xymatrixrowsep{1.75pc}\xymatrixcolsep{3.75pc}\xymatrix{
(A\ot B)\op(C\ot D)\ar[r]^-{\mu_{A,B,C,D}}\ar[d]_-{\mix_{A,B}\op\mix_{C,D}} & (A\op C)\ot (B\op D)\ar[d]^-{\mix_{A\op C, B\op D}} \\
(A\op B)\op(C\op D)\ar[r]_-{\tau^\op_{A, B, C, D}} & (A\op C)\op(B\op D) \\
(A\ot B)\ot (C\ot D)\ar[r]^-{\tau^\ot_{A,B,C,D}}\ar[d]_-{\mix_{A\ot B, C\ot D}} & (A\ot C)\ot (B\ot D)\ar[d]^-{\mix_{A,C}\ot \mix_{B\ot D}} \\
(A\ot B)\op(C\ot D)\ar[r]_-{\mu_{A,B,C,D}} & (A\op C)\ot(B\op D) 
}
\end{gathered}\end{equation}  
\end{proposition}
\begin{proof}
The proof of the first equality is given by Figure \ref{fig:interaction_canonical_flip_medial_maps_mix} in Appendix \ref{app:commuting_diagrams} and the other equality follows similarly. 
\end{proof}

Now, if the linearly distributive structure of a SMLDC is compact, then by \eqref{diag:interchange_canonicalflip_mix}, we know the medial maps are isomorphisms as the mix maps and the canonical flip maps are. Moreover, $\Delta_\bot$ and $\nabla_\top$ are invertible by ({\footnotesize{\bf MLDC.2}}), since the nullary mix map is invertible. Thus, the duoidal structure is strong. Altogether, we get the following, which generalizes the CLDC result:

\begin{theorem}
Given a SMLDC, the following are equivalent:
\begin{enumerate}[(i)]
    \item the linearly distributive structure is compact and the duoidal structure is strong,
    \begin{itemize}
        \item it is isomix,
        \item mix maps $\mix_{A,B}\c A\ot B\rarr A\op B$ are isomorphisms, 
        \item linear distributors are associators, modulo the mix maps (see \cite[Lem 2.7]{Kudzman-Blais_Lemay_2025}),
        \item $\bot$-contraction/$\top$-cocontraction are unitors, modulo the nullary mix map,  and
        \item medial maps are canonical flip maps, modulo the mix maps.
    \end{itemize}
    \item the linear distributors are isomorphisms, and
    \item the linearly distributive structure is isomix.
\end{enumerate}
\end{theorem}

The above results demonstrate why the landscape of examples of MLDCs and SMLDCs is currently restrained. We hope however that more examples will surface as the theory of MLDCs and CLDCs continues to be explored.

\section{Medial Linear Functors and Transformations}\label{sec:medial_linear_functor}

With the definition of SMLDCs established, we now turn to the appropriate functors and transformations between these categories. Similar to the 2-category \CLDC, the corresponding 2-category should form a sub-2-category of \SLDC. 

As such, these ``medial linear functors'' are linear functors which also interact coherently with the duoidal structure. Given that there are two established definitions of duoidal functors, there are naturally at least two distinct types of ``medial linear functors''. We define only the appropriate class required for the linearly distributive Fox theorem. 

\begin{definition}
    Let $\bX$ and $\bY$ be SMLDCs. A {\bf symmetric medial linear functor} $F = (F_\ot, F_\op)\c \bX\rarr\bY$ consists of:
    
\begin{itemize}
    \item a functor $F_\ot\c \bX\rarr \bY$, equipped with maps and natural transformations
    \begin{align*}
    &m_{\top}\c\top\xrightarrow{\sim} F_\ot(\top) &&\quad m_{\bot}\c\bot\rarr F_\ot(\bot)\\
    &{m_\ot}_{A, B}\c F_\ot(A)\ot F_\ot(B)\xrightarrow{\sim} F_\ot(A\ot B) &&\quad {m_\op}_{A, B}\c F_\ot(A)\op F_\ot(B)\rarr F_\ot(A\op B)
    \end{align*}
    such that
     \begin{enumerate}[(i)]
        \item $(F_\ot, m_{\top},m_\ot)\c (\bX, \ot, \top)\rarr(\bY,\ot,\top)$ is a symmetric monoidal functor,
         \item $(F_\ot, m_{\bot}, m_\op)\c(\bX, \op, \bot)\rarr(\bY,\op,\bot)$ is a symmetric lax monoidal functor, 
     \end{enumerate}
        
    \item a functor $F_\op\c \bX\rarr \bY$, equipped with maps and natural transformations
    \begin{align*}
    &n_{\bot}\c F_\op(\bot)\xrightarrow{\sim} \bot &&\quad n_{\top}\c F_\op(\top)\rarr\top\\
    &{n_\op}_{A,B}\c F_\op(A\op B)\xrightarrow{\sim} F_\op(A)\op F_\op(B) &&\quad {n_\ot}_{A,B}\c F_\op(A\ot B)\rarr F_\op(A)\ot F_\op(B)
    \end{align*}
    such that 
    \begin{enumerate}[(i)]
         \item \((F_\op, n_{\bot}, n_\op)\c(\bX, \op, \bot)\rarr(\bY,\op,\bot)\) is a symmetric monoidal functor,
         \item \((F_\op, n_{\top}, n_\ot)\c(\bX, \op, \bot)\rarr(\bY,\op,\bot)\) is a symmetric colax monoidal functor,
    \end{enumerate}
    
    \item linear strength natural transformations
    \begin{align*}
        &{v_\ot^R}_{A,B} \c F_\ot(A\op B)\rarr F_\op(A)\op F_\ot(B) & {v_\op^R}_{A,B} \c F_\ot(A)\ot F_\op(B)\rarr F_\op(A\ot B)
    \end{align*}
\end{itemize}   
satisfying the coherence conditions ensuring that
\begin{enumerate}[{\bf \footnotesize (SMLF.1)}]
    \item $(F_\ot, m_{\bot},m_\op, m_{\top}^{-1}, m_\ot^{-1})$ is a bilax duoidal functor,

\begin{fullwidth}
\begin{equation*}\label{cc:F_ot_bilax}
\begin{gathered}\begin{aligned}
&\Delta_{\bot};(m_{\bot}\ot m_{\bot});{m_\ot}_{\bot,\bot} = m_{\bot};F_\ot(\Delta_{\bot}) \\
&(m_{\top}\op m_{\top}); {m_\op}_{\top,\top}; F_{\ot}(\nabla_{\top}) = \nabla_{\top};m_{\top} \\ 
& m;m_{\top} = m_{\bot};F_\ot(m)\\[0.5pc]
&({m_\ot}_{A,B}\op {m_\ot}_{C,D}); {m_\op}_{A\ot B, C\ot D}; F_\ot(\mu_{A,B,C,D}) \\
&= \mu_{F_\ot(A), F_\ot(B), F_\ot(C), F_\ot(D)}; ({m_\op}_{A,C}\ot {m_\op}_{B, D}); {m_\ot}_{A\op C,B\op D}
\end{aligned}\end{gathered}
\end{equation*}    
\end{fullwidth}
    
    \item $(F_\op, n_{\bot}^{-1}, n_\op^{-1}, n_{\top}, n_\ot)$ is a bilax duoidal functor,

\begin{fullwidth}
\begin{equation*}\label{cc:F_op_bilax}
\begin{gathered}\begin{aligned}
&F_\op(\Delta_{\bot}); {n_\ot}_{\bot,\bot}; (n_{\bot}\ot n_{\bot}) = n_{\bot};\Delta_{\bot} \\
& {n_\op}_{\top,\top}; (n_{\top}\op n_{\top}); \nabla_{\top} = F_\op(\nabla_{\top}); n_{\top} \\
& n_{\bot};m=F_\op(m);n_{\top} \\[0.5pc]
&F_\op(\mu_{A,B,C,D}); {n_\ot}_{A\op C, B\op D};({n_\op}_{A,C}\ot {n_\op}_{B,D})\\
&= {n_\op}_{A\ot B, C\ot D};({n_\ot}_{A,B}\op {n_\ot}_{C,D}); \mu_{F_\op(A), F_\op(B), F_\op(C), F_\op(D)}
\end{aligned}\end{gathered}
\end{equation*}
\end{fullwidth}
    
    \item $F  = (F_\ot, F_\op)\c\bX\rarr\bY$ is a symmetric linear functor, i.e. {\bf \footnotesize (LF.1-5)} and {\bf \footnotesize (SLF)},
    \item the linear strengths interact coherently with $\top$-cocontraction and $\bot$-contraction,
\begin{fullwidth}
\begin{equation*}\begin{gathered}\label{cc:linear_strength_nabla_delta}
\xymatrixrowsep{1.75pc}\xymatrixcolsep{1pc}\xymatrix{
F_\ot(\top\op\top)\ar[d]_-{F_\op(\nabla_\top)}\ar[rr]^-{{\nu^R_\ot}_{\top,\top}} & & F_\op(\top) \op F_\ot(\top)\ar[d]^-{n_\top \op m_\top^{-1}} & F_\op(\bot)\ar[r]^-{n_\bot}\ar[d]_-{F_\op(\Delta_\bot)} & \bot\ar[r]^-{\Delta_\bot} & \bot\ot\bot\ar[d]^-{m_\bot\ot n_\bot^{-1}} \\
F_\ot(\top) & \top \ar[l]^-{m_\top}& \top \op \top \ar[l]^-{\nabla_\top} & F_\op(\bot\ot\bot) & & F_\ot(\bot)\ot F_\op(\bot) \ar[ll]^-{{\nu^R_\op}_{\bot,\bot}}
}
\end{gathered}\end{equation*}   
\end{fullwidth}
    
    \item  the linear strengths interact coherently with the medial maps, and 
\begin{fullwidth}
\begin{equation*}\begin{gathered}\label{cc:linear_strength_medial}
\xymatrixrowsep{1.75pc}\xymatrixcolsep{3.75pc}\xymatrix{
F_\ot((A\ot B)\op (C\ot D))\ar[r]^-{{\nu_\ot^R}_{A\ot B, C\ot D}}\ar[d]_-{F_\ot(\mu_{A,B,C,D})} & F_\op(A\ot B)\op F_\ot(C\ot D)\ar[d]^-{{n_\ot}_{A,B}\op {m_\ot}^{-1}_{C,D}} \\
F_\ot ((A\op C)\ot (B\op D))\ar[d]_-{{m_\ot}^{-1}_{A\op C, B\op D}}  & (F_\op(A)\ot F_\op(B))\op (F_\ot(C)\ot F_\ot(D))\ar[d]^-{\mu_{F_\op(A), F_\op(B), F_\ot(C), F_\ot(D)}}  \\
F_\ot(A\op C)\ot F_\ot(B\op D)\ar[r]_-{{\nu_\ot^R}_{A,C}\ot {\nu_\ot^R}_{B,D}} & (F_\op(A) \op F_\ot(C)) \ot (F_\op(B)\op F_\ot(D)) \\
(F_\ot(A) \ot F_\op(B)) \op (F_\ot(C)\ot F_\op(D))\ar[r]^-{{\nu_\op^R}_{A,B}\op{\nu_\op^R}_{C,D}}\ar[d]_-{\mu_{F_\ot(A), F_\op(B), F_\ot(C), F_\op(D)}} & F_\op(A\ot B) \op F_\op (C\ot D)\ar[d]^-{{n_\op}^{-1}_{A\ot B, C\ot D}} \\
(F_\ot(A) \op F_\ot(C))\op (F_\op(B)\op F_\op(D)\ar[d]_-{{m_\op}_{A,C}\ot {n_\op}^{-1}_{B,D}}  & F_\op((A\ot B)\op(C\ot D))\ar[d]^-{F_\op(\mu_{A,B,C,D})}\\
F_\ot(A\op C) \ot F_\op(B\op D)\ar[r]_-{{\nu^R_\op}_{A\op C, B\op D} } & F_\op((A\op C)\ot (B\op D))
}
\end{gathered}\end{equation*}   
\end{fullwidth}
    
    \item the linear strengths interact coherently with $m_\op$ and $n_\ot$.
\begin{fullwidth}
\begin{equation*}\label{cc:linear_strength_m_op_n_ot}
\begin{gathered}\xymatrixrowsep{1.75pc}\xymatrixcolsep{6.75pc}\xymatrix{
F_\ot(A\op B)\op F_\ot(C)\ar[r]^-{{\nu^R_\ot}_{A, B}\op 1_{F_\ot(C)}}\ar[d]_-{{m_\op}_{A\op B, C}} & (F_\op(A) \op F_\ot(B))\op F_\ot(C) \ar[d]^-{ {\alpha_\op}^{-1}_{F_\op(A), F_\ot(B), F_\ot(C)}}\\
F_\ot((A\op B)\op C)\ar[d]_-{F_\ot({\alpha_\op}^{-1}_{A, B, C})} & F_\op(A) \op (F_\ot(B) \op F_\ot(C))\ar[d]^-{1_{F_\op(A)}\op {m_\op}_{B, C}}\\
F_\ot(A\op(B\op C))\ar[r]_-{{\nu^R_\ot}_{A,B\op C}} & F_\op(A) \op F_\ot(B\op C) \\
F_\ot(A) \ot F_\op(B\ot C)\ar[r]^-{{\nu^R_\op}_{A, B\ot C}} \ar[d]_-{1_{F_\ot(A)}\ot {n_\ot}_{B,C}} & F_\op(A \ot(B\ot C))\ar[d]^-{F_\op({\alpha_\ot}^{-1}_{A, B, C})} \\
F_\ot(A) \ot (F_\op(B)\ot F_\ot(C))\ar[d]_-{{\alpha_\ot}^{-1}_{F_\ot(A), F_\op(B), F_\op(C)}} & F_\op((A\ot B)\ot C) \ar[d]^-{{n_\ot}_{A\ot B, C} } \\
(F_\ot(A) \ot F_\op(B)) \ot F_\op(C)\ar[r]_-{{\nu^R_\op}_{A,B}\ot 1_{F_\op(C)}} & F_\op(A\ot B)\ot F_\op(C)
}
\end{gathered}\end{equation*}
\end{fullwidth}
\end{enumerate}
\end{definition}

The coherence conditions above ensure that the linear strengths interact coherently with the duoidal functor structures.

\begin{proposition}\label{prop:interaction_canonical_flip_medial_linear_functor}
Given a symmetric medial linear functor $F=(F_\ot, F_\op)\c\bX\rarr\bY$, the following diagrams and equations hold. 
\begin{equation}\begin{gathered}\label{eqn:interaction_linear_strength_duoidal_trans_1}
\xymatrixrowsep{1.75pc}\xymatrixcolsep{2.75pc}\xymatrix{
F_\ot(\bot)\ar[d]_-{F_\ot({u_\op}^{-1}_\bot)} & \bot\ar[l]_-{m_\bot}\ar[r]^-{{u_\op}^{-1}_\bot} & \bot\op\bot\ar[d]^-{n_\bot^{-1}\op m_\bot} \\
F_\ot(\bot\op\bot)\ar[rr]_-{{\nu^R_\ot}_{\bot,\bot}} && F_\ot(\bot)
}
\end{gathered}\end{equation}
\vspace{-\baselineskip}
\begin{align*}
& m_\bot; F_\ot({u_\op}^{-1}_\bot); {\nu^R_\ot}_{\bot,\bot} = {u_\op}^{-1}_\bot; n_\bot^{-1}\op m_\bot \\
&{\nu^R_\op}_{\top,\top}; F_\op({u_\ot}^{-1}_\top); n_\top = (m_\top^{-1}\ot n_\top);{u_\ot}^{-1}_\top
\end{align*}
\begin{equation}\begin{gathered}\label{eqn:interaction_canonical_flip_medial_linear_functor}
\xymatrixrowsep{1.75pc}\xymatrixcolsep{4.75pc}\xymatrix{
F_\ot(A\op B) \op F_\ot(C\op D)\ar[r]^-{{\nu^R_\ot}_{A,B}\op{\nu^R_\ot}_{C,D}}\ar[d]_-{{m_\op}_{A\op B, C\op D}} & (F_\op(A)\op F_\ot(B)) \op (F_\op(C)\op F_\ot(D))\ar[d]^-{\tau^\op_{F_\op(A), F_\ot(B), F_\op(C), F_\ot(D)}}\\
F_\ot(A\op B) \op (C\op D))\ar[d]_-{F_\ot(\tau^\op_{A,B,C,D})} & (F_\op (A)\op F_\op(C)) \op (F_\ot(B)\op F_\ot(D))\ar[d]^-{{n_\op}_{A,C}^{-1} \op {m_\op}_{B,D}} \\
F_\op((A\op C)\op (B\op D))\ar[r]_-{{\nu^R_\ot}_{A\op C, B\op D}} &F_\op(A\op C) \op F_\ot(B\op D)
}
\end{gathered}\end{equation}
\vspace{-\baselineskip}
\begin{align*}
& ({\nu^R_\ot}_{A,B}\op{\nu^R_\ot}_{C,D}); \tau^\op_{F_\op(A), F_\ot(B), F_\op(C), F_\ot(D)}; ({n_\op}^{-1}_{A,C} \op {m_\op}_{B,D})\\
&= {m_\op}_{A\op B, C\op D}; F_\ot(\tau^\op_{A,B,C,D}); {\nu^R_\ot}_{A\op C, B\op D} \\[0.5pc]
& {\nu^R_\op}_{A\ot B, C\ot D}; F_\op(\tau^\ot_{A,B,C,D}); {n_\ot}_{A\ot C, B\ot D}\\
&= ({m_\ot}^{-1}_{A,B}\ot {n_\ot}_{C,D}); \tau^\ot_{F_\ot(A), F_\ot(B), F_\op(C), F_\op(D)};({\nu^R_\op}_{A,C} \ot {\nu^R_\op}_{B,D}) 
\end{align*}
As such, the linear strengths are duoidal transformations between bilax duoidal functors.
\end{proposition}
\begin{proof}
The first equality in \eqref{eqn:interaction_linear_strength_duoidal_trans_1} holds by the following commuting diagram, with the second equality holding similarly. 
\begin{equation*}
\xymatrixrowsep{1.75pc}\xymatrixcolsep{2.75pc}\xymatrix{
\bot\ar@{=}[rr]\ar[d]_-{m_\bot}\ar@{}[rrd]|{(\nat)} && \bot\ar[d]^-{{u_\op}^{-1}_\bot}  \\
F_\ot(\bot)\ar[d]_-{F_\ot({u_\op}^{-1}_\bot)}\ar[r]^-{{u^L_\op}^{-1}_{F_\ot(\bot)}}\ar@{}[rrd]|{(\mathrm{\bf LF.1})} & \bot\op F_\ot(\bot) & \bot\op\bot\ar[d]^-{n_\bot^{-1}\op m_\bot}\ar[l]_-{1_\bot\op m_\bot} \\
F_\ot(\bot\op\bot)\ar[rr]_-{{\nu^R_\ot}_{\bot,\bot}} && F_\ot(\bot)
}
\end{equation*}
The proof for the first equality in \eqref{eqn:interaction_canonical_flip_medial_linear_functor} is given by Figure \ref{fig:proof_interaction_canonical_flip_medial_linear_functor} of Appendix \ref{app:commuting_diagrams} and the second equality is proved similarly.

Now,the monoidal products $\ot$ and $\op$ are bilax duoidal functors by Proposition \ref{prop:monoidal_products_bilax_duoidal_functor}, and $F_\ot$ and $F_\op$ are bilax duoidal functors by definition. As such, the linear strengths are natural transformations between functors which are canonically bilax duoidal:
\begin{align*}
    &{v_\ot^R} \c \op; F_\ot \Rarr (F_\op\times F_\ot);\op  & {v_\op^R}\c (F_\ot\times F_\op);\ot\Rarr \ot;F_\op
\end{align*}

The first equality in \eqref{eqn:interaction_linear_strength_duoidal_trans_1} and the first equality in \eqref{eqn:interaction_canonical_flip_medial_linear_functor} are precisely the conditions such that ${v_\ot^R}$ is a $\op$-monoidal transformation, while the first diagrams of {\bf \footnotesize (SMLF.4)} and {\bf \footnotesize (SMLF.5)} are the conditions such that ${v_\ot^R}$ is a $\ot$-monoidal transformation. Therefore, ${v_\ot^R}$ is a duoidal transformation. Similarly, ${v_\op^R}$ is a duoidal transformation. 
\end{proof}

As there are several potential notions of medial linear functors, there are multiple definitions of ``medial linear transformations''. However, we will introduce only the specific definition relevant to our treatment of MLDCs.

\begin{definition}
Let $F, G\c\bX\rarr\bY$ be symmetric medial linear functors. A {\bf medial linear transformation} $\alpha = (\alpha_\ot, \alpha_\op)\c F\Rarr G$ consists of a pair of natural transformations 
\[ {\alpha_\ot}_A\c F_\ot(A) \rarr G_\ot(A) \qquad {\alpha_\op}_A\c G_\op(A)\rarr F_\op(A)\]
such that 
\begin{enumerate}[(i)]
    \item $\alpha_\ot\c (F_\ot, m_{\top}^{F}, m_\ot^{F}) \Rarr (G_\ot, m_{\top}^{G}, m_\ot^{G})$ is a monoidal transformation,
    \item  $\alpha_\ot\c (F_\ot, m_{\bot}^{F},m_\op^{F})\Rarr (G_\ot,m_{\bot}^{G}, m_\op^{G})$ is a monoidal transformation,
    \item $\alpha_\op\c(G_\op, n_{\bot}^{G}, n_\op^{G}) \Rarr (F_\op, n_{\bot}^{F}, n_\op^{F})$ is a monoidal transformation,
    \item $\alpha_\op\c (G_\op, {n_{\top}}^{G}, {n_\ot}^{G})\Rarr(F_\op, {n_\top^F}, {n_\ot}^{F})$ is a monoidal transformation, and 
    \item $\alpha = (\alpha_\ot, \alpha_\op)$ is a linear transformation, in other words \eqref{cc:linear_trans_linear strength} holds.
\end{enumerate}	
\end{definition}

\begin{remark}\label{rem:medial_linear_trans}
Notice that conditions \textit{(i-ii)} in the above definition are equivalent to asking that
\[\alpha_\ot\c (F_\ot, m_{\bot}^{F}, {m_\op}^{F}, {m_{\top}^{-1}}^{F},  {m_\ot^{-1}}^{F})\Rarr(G_\ot, m_{\bot}^{G}, {m_\op}^{G}, {m_{\top}^{-1}}^{G},{m_\ot^{-1}}^{G})\]
is a bilax duoidal transformation and, similarly, conditions \textit{(iii-iv)} are equivalent to
\[\alpha_\op\c (G_\op, {n_{\bot}^{-1}}^{G},{n_\op^{-1}}^{G}, {n_\op}^{G},{n_{\top}}^{G})\Rarr(F_\op,{n_{\bot}^{-1}}^{F},{n_\op^{-1}}^{F}, {n_\op}^{F}, {n_{\top}}^{F})\]
being a bilax duoidal transformation.
\end{remark}

\subsection{2-Category \SMLDC} \hfill\

\vspace{0.5\baselineskip}
These medial functors and transformations, alongside SMLDCs, give us a 2-category:
\begin{theorem}
SMLDCs, symmetric medial linear functors, and medial linear transformations form a 2-category, denoted by \SMLDC.  
\end{theorem}
\begin{proof}
Let $F=(F_\ot,F_\op)\c \bX\rarr\bY$ and $G=(G_\ot, G_\op)\c \bY\rarr\bZ$ be symmetric medial linear functors. Consider their horizontal composite $F; G = (F_\ot; G_\ot, F_\op; G_\op)\c \bX\rarr \bY$, with the obvious structure maps. It is immediate that $F;G\c\bX\rarr\bZ$ is a symmetric linear functor, and that $F_\ot; G_\ot \c \bX\rarr \bZ$ and $F_\op; G_\op\c \bX\rarr \bZ$ are bilax duoidal functors, as the composition is inherited from \LDC\ and \DUO\ in the appropriate ways. It remains only to show the additional coherence conditions also hold for $F;G\c\bX\rarr\bZ$. 

        
    

The second diagram of {\bf \footnotesize (SMLF.4)} holds by the following commuting diagram.
\begin{fullwidth}
\begin{equation*}\begin{gathered}
\xymatrixrowsep{2.75pc}\xymatrixcolsep{2pc}\xymatrix{
G_\op(F_\op(\bot))\ar[r]^-{G_\op(n^F_\bot)}\ar[dd]_-{G_\op(F_\op(\Delta_\bot))} & G_\op(\bot)\ar[d]_-{G_\op(\Delta_\bot)}\ar[r]^-{n^G_\bot} & \bot\ar[r]^-{\Delta_\bot} & \bot\ot \bot\ar[d]^-{m^G_\bot \ot {n^G_\bot}^{-1}} \ar@{}[dll]|{(\mathrm{\bf SMLF.4})_G}\\
\ar@{}[r]|{(\mathrm{\bf SMLF.4})_F} & G_\op(\bot\ot\bot)\ar[d]_-{G_\op(m^F_\bot \ot {n^F_\bot}^{-1})} && G_\ot(\bot)\ot G_\op(\bot)\ar[d]^-{G_\ot(m^F_\bot)\ot G_\op({n^F_\bot}^{-1})} \ar[ll]_-{{\nu^R_\op}^G_{\bot,\bot}} \\
G_\op(F_\op(\bot\ot\bot)) & G_\op(F_\ot(\bot)\ot F_\op(\bot))\ar[l]^*+<0.3em>{^{G_\op({\nu^R_\op}^F_{A,B})}}\ar@{}[rru]|{(\nat)} && G_\ot(F_\ot(\bot))\ot G_\op(F_\op(\bot)) \ar[ll]^-{{\nu^R_\op}^G_{F_\ot(\bot),F_\op(\bot)}}
}
\end{gathered}\end{equation*}
\end{fullwidth}
The first diagram {\bf \footnotesize (SMLF.5)} holds by the following commuting diagram.
\begin{equation*}
\resizebox{\linewidth}{!}{
\xymatrixrowsep{3.75pc}\xymatrixcolsep{5.75pc}\xymatrix{
{\begin{array}{@{}c@{}}G_\ot(F_\ot((A\ot B)\op{}\\ (C\ot D)))\end{array}}\ar[d]_-{G_\ot(F_\ot(\mu_{A,B,C,D}))}\ar[r]_-{G_\ot({\nu^R_\ot}^F_{A\ot B, C\ot D})} & {\begin{array}{@{}c@{}}G_\ot(F_\op(A\ot B) \op{}\\ F_\ot(C\ot D))\end{array}}\ar[r]_-{{\nu^R_\ot}^G_{F_\op(A\ot B), F_\ot(C\ot D)}}\ar[d]_-{G_\ot({n^F_\ot}_{A,B}\op {m^F_\ot}^{-1}_{C,D})}^-{\qquad\qquad\qquad(\nat)} & {\begin{array}{@{}c@{}}G_\op(F_\op(A\ot B)) \op{}\\ G_\ot(F_\ot(C\ot D))\end{array}}\ar[d]_-{G_\ot({n^F_\ot}_{A,B}) \op G_\ot({m^F_\ot}^{-1}_{C,D})}\\
{\begin{array}{@{}c@{}}G_\ot(F_\ot((A\op C)\ot{}\\ (B\op D)))\end{array}} \ar[d]_-{G_\ot({m^F_\ot}^{-1}_{A\op C, B\op D})}\ar@{}[r]|{(\mathrm{\bf SMLF.5})_F} & {\begin{array}{@{}c@{}} G_\ot((F_\op(A)\ot F_\op(B))\op{}\\ (F_\ot(C)\ot F_\ot(D))) \end{array}}\ar[r]_-{{\nu^R_\ot}^G_{F_\op(A)\ot F_\op(B), F_\ot(C)\ot F_\ot(D)}} \ar[d]_-{G_\ot(\mu_{F_\op(A), F_\op(B), F_\ot(C), F_\ot(D)})} & {\begin{array}{@{}c@{}}G_\op(F_\op(A)\ot F_\op(B)) \op{}\\  G_\ot(F_\ot(C)\ot F_\ot(D))\end{array}} \ar[d]_-{{n^G_\ot}_{F_\op(A),F_\ot(B)} \op{m^G_\ot}^{-1}_{F_\ot(C),F_\op(D)}} \\
{\begin{array}{@{}c@{}}G_\ot(F_\ot(A\op C) \ot{}\\F_\ot(B\op D))\end{array}} \ar[d]_-{{m^G_\ot}^{-1}_{F_\ot(A\op C), F_\ot(B\op D)}}^-{\qquad(\nat)} \ar[r]_-{G_\ot({\nu^R_\ot}^F_{A,C} \ot{\nu^R_\ot}^F_{B,D})} & {\begin{array}{@{}c@{}} G_\ot((F_\op(A)\op F_\ot(C))\ot{}\\ (F_\op(B)\op F_\ot(D))) \end{array}} \ar[d]_-{{m^G_\ot}^{-1}_{F_\op(A)\op F_\ot(B), F_\ot(C)\ot F_\ot(D)}}\ar@{}[r]|{(\mathrm{\bf SMLF.5})_G} & {\begin{array}{@{}c@{}}(G_\op(F_\op(A)) \ot G_\op(F_\op(B))) \op{}\\ (G_\ot(F_\ot(C))\ot G_\ot(F_\ot(D)))\end{array}}\ar[d]_-{\mu_{G_\op(F_\op(A)), G_\op(F_\op(B)), G_\ot(F_\ot(C)), G_\ot(F_\ot(D))}} \\
{\begin{array}{@{}c@{}}G_\ot(F_\ot(A\op C)) \ot{}\\ G_\ot(F_\ot(B\op D))\end{array}} \ar[r]_*+<1.5em>{_{G_\ot({\nu^R_\ot}^F_{A,C}) \ot G_\ot({\nu^R_\ot}^{F}_{B,D})}} & {\begin{array}{@{}c@{}} G_\ot(F_\op(A)\op F_\ot(C))\ot{}\\ G_\ot(F_\op(B) \op F_\ot(D))\end{array}} \ar[r]_*+<1.5em>{_{{\nu^R_\ot}^G_{F_\op(A),F_\ot(C)} \op {\nu^R_\ot}^{G}_{F_\op(B),F_\op(D)}}}  & {\begin{array}{@{}c@{}} (G_\op(F_\op(A)) \op G_\ot(F_\ot(C))) \ot{}\\ (G_\op(F_\op(B))\op G_\ot(F_\ot(D)))\end{array}}
}}
\end{equation*}
The first diagram of {\bf \footnotesize (SMLF.6)} holds by the following commuting diagram.
\begin{equation*}
\resizebox{\linewidth}{!}{
\xymatrixrowsep{3.75pc}\xymatrixcolsep{3.75pc}\xymatrix@L=0.4pc{
{\begin{array}{@{}c@{}}G_\ot(F_\ot(A\op B)) \op{}\\ G_\ot(F_\ot(C))\end{array}}\ar[d]_-{{m^G_\op}_{F_\ot(A\op B), F_\ot(C)}}^-{\qquad\qquad(\nat)}\ar[r]_-{G_\ot({\nu^R_\ot}^F_{A,B})\op 1_{G_\ot(F_\ot(C))}} & {\begin{array}{@{}c@{}}G_\ot(F_\op(A)\op F_\ot(B)) \op{}\\ G_\ot(F_\ot(C))\end{array}}\ar[d]_-{{m^G_\op}_{F_\op(A)\op F_\ot(B), F_\ot(C)}}\ar[r]_-{{\nu^R_\ot}^G_{F_\op(A), F_\ot(B)}\op1_{G_\ot(F_\ot(C))}} & {\begin{array}{@{}c@{}}(G_\op(F_\op(A)) \op G_\ot(F_\ot(B)))\op{}\\  G_\ot(F_\ot(C))\end{array}}\ar[d]_-{{\alpha_\op}^{-1}_{G_\op(F_\op(A)), G_\ot(F_\ot(B)), G_\ot(F_\ot(C)}} \\
G_\ot(F_\ot(A\op B)\op F_\ot(C))\ar[d]_-{G_\ot({m^F_\op}_{A\op B, C})}\ar[r]^-{G_\ot({\nu^R_\ot}^F_{A,B}\op 1_{F_\ot(C)})} & G_\ot((F_\op(A)\op F_\ot(B)) \op F_\ot(C))\ar[d]_-{G_\ot({\alpha_\op}^{-1}_{F_\op(A), F_\ot(B), F_\ot(C)})}\ar@{}[r]|{\quad(\mathrm{\bf SMLF.6})_G}& {\begin{array}{@{}c@{}}G_\op(F_\op(A)) \op{}\\  (G_\ot(F_\ot(B))\op G_\ot(F_\ot(C))) \end{array}} \ar[d]_-{1_{G_\op(F_\op(A))} \op {m_\op^G}_{F_\ot(B), F_\ot(C)}} \\
G_\ot(F_\ot((A\op B)\op C))\ar[d]_-{G_\ot(F_\ot({\alpha_\op}^{-1}_{A,B,C}))}\ar@{}[r]|{(\mathrm{\bf SMLF.6})_F\qquad} &  G_\ot(F_\op(A) \op (F_\ot(B) \op F_\ot(C)))\ar[d]_-{G_\ot(1_{F_\op(A)}\op {m^F_\op}_{B,C})}^-{\qquad\qquad(\nat)}\ar[r]^-{{\nu^R_\ot}^G_{F_\op(A), F_\ot(B)\op F_\ot(C)}} & {\begin{array}{@{}c@{}}G_\op(F_\op(A)) \op{}\\ G_\ot(F_\ot(B)\op F_\ot(C)) \end{array}} \ar[d]_-{1_{G_\op(F_\op(A))} \op G_\ot({m^F_\op}_{B,C})}\\
G_\ot(F_\ot(A\op (B\op C)))\ar[r]_-{G_\ot({\nu^R_\ot}^F_{A,B\op C})} & G_\ot(F_\op(A) \op F_\ot(B\op C))\ar[r]_-{{\nu^R_\ot}^G_{F_\op(A), F_\ot(B\op C)}} & {\begin{array}{@{}c@{}}G_\op(F_\op(A)) \op{}\\ G_\ot(F_\ot(B\op C)) \end{array}} 
}}
\end{equation*}
The other diagrams of {\bf \footnotesize (SMLF.4-6)} hold similarly. 

Finally, it is immediate to show that the identity linear functors and the identity linear transformations are inherited by \SMLDC\ and that medial linear transformations are closed under the standard vertical and horizontal composition of linear transformations.
\end{proof}

\begin{proposition}\label{prop:inclusion_CLDC_SMLDC}
The inclusion map $\inc\c \CLDC\rarr\SMLDC$ determines a 2-functor.
\end{proposition}
\begin{proof}
By Theorem \ref{thm:CLDC_is_SMLDC}, every CLDC has a canonical SMLDC structure.

Consider a cartesian linear functor
\[ F=(F_\times, F_+)\c (\bX, \times, \bone, +, \bzero)\rarr (\bY, \times, \bone, +, \bzero) \]  Then, in particular, $(F_\times, m_\bone^{-1}, m^{-1}_\times)\c (\bX,\times, \bone)\rarr (\bY, \times, \bone)$ and $(F_+, n^{-1}_\bzero, n^{-1}_+)\c (\bX, +,\bzero)\rarr (\bY, +, \bzero)$ are respectively colax monoidal and lax monoidal functors. By Proposition \ref{prop:cartesian_duoidal_functor}, the following are bilax duoidal functors, \[ (F_\times, p_\bzero, p_+, m_\bone^{-1}, m^{-1}_\times)\c (\bX, +,\bzero, \times, \bone)\rarr (\bY, +,\bzero, \times, \bone)\]\[ (F_+, n_\bzero^{-1}, n_+^{-1}, q_\bone, q_\times)\c (\bX, +,\bzero, \times, \bone)\rarr (\bY, +,\bzero, \times, \bone)\] where $p_\bzero, p_\times, q_\bone$ and $q_\times$ are defined in Lemma \ref{lem:functor_comon_if_cart}. In order to prove $F=(F_\times, F_+)$ is a symmetric medial linear functor, it remains to show the additional coherence conditions. 

{\bf \footnotesize (SMLF.4)} holds by the universal property of initial and terminal objects:
\begin{align*}
F_\times(\bone+\bone) \xrightarrow{t_{F_\times(\bone+\bone)}} \bone&=F_\times(\bone+\bone) \xrightarrow{{\nu_\times^R}_{\bone,\bone}} F_+(\bone)+F_\times(\bone) \xrightarrow{q_{\bone}+ m^{-1}_{\bone}} \bone+\bone \xrightarrow{\nabla_\bone} \bone \\
&=F_\times(\bone+\bone) \xrightarrow{F_\times(\nabla_{\bone})} F_\times(\bone)\xrightarrow{m_{\bone}^{-1}} \bone \\ \\
\bzero \xrightarrow{b_{F_+(\bzero\times\bzero)}} F_+(\bzero\times\bzero) &=\bzero \xrightarrow{\Delta_\bzero} \bzero\times\bzero \xrightarrow{p_{\bzero}\times n_{\bzero}^{-1}} F_\times(\bzero)\times F_\times(\bzero) \xrightarrow{{\nu_+^R}_{\bzero,\bzero} } F_+(\bzero\times\bzero) \\
& = \bzero \xrightarrow{n_\bzero^{-1}} F_+(\bzero) \xrightarrow{F_+(\Delta_\bzero)} F_+(\bzero\times\bzero) 
\end{align*}

The first coherence condition of {\bf \footnotesize (SMLF.5)} is 
\begin{align*}
&F_\times(\mu_{A,B,C,D}); {m_\times}^{-1}_{A+ C, B+ D}; ({\nu_\times^R}_{A,C}\times {\nu_\times^R}_{B,D})  \\
&={\nu_\times^R}_{A\times B, C\times D}; ({q_\times}_{A,B}+ {m_\times}^{-1}_{C,D}); \mu_{F_+(A), F_+(B), F_\times(C), F_\times(D)}
\end{align*}
and it holds by the universal property of products: consider the left-hand and right-hand sides post composed with the relevant projections
\begin{align*}
&F_\times(\mu_{A,B,C,D}); {m_\times}^{-1}_{A+ C, B+ D}; ({\nu_\times^R}_{A,C}\times {\nu_\times^R}_{B,D}); \pi^0_{F_+(A)+F_\times(C), F_+(B)+ F_\times(D)} \\
&= F_\times(\mu_{A,B,C,D}); {m_\times^{-1}}_{A+ C, B+ D}; \pi^0_{F_\times(A+C), F_\times(B+D)}; {\nu_\times^R}_{A,C} \\
&= F_\times(\mu_{A,B,C,D}); F_\times(\pi^0_{A+C, B+D}); {\nu_\times^R}_{A,C} = F_\times(\pi^0_{A,B}+\pi^0_{C,D}); {\nu_\times^R}_{A,C} \\
\\
&{\nu_\times^R}_{A\times B, C\times D}; ({q_\times}_{A,B}+ {m_\times}^{-1}_{C,D}); \mu_{F_+(A), F_+(B), F_\times(C), F_\times(D)}; \pi^0_{F_+(A)+F_\times(C), F_+(B)+ F_\times(D)} \\
&= {\nu_\times^R}_{A\times B, C\times D}; ({q_\times}_{A,B}+ {m^{-1}_\times}_{C,D}); (\pi^0_{F_+(A), F+(B)}+\pi^0_{F_\times(C), F_\times(D)}) \\
&=  {\nu_\times^R}_{A\times B, C\times D}; (F_+(\pi^0_{A, B})+F_\times(\pi^0_{C,D})) = F_\times(\pi^0_{A,B}+\pi^0_{C,D}); {\nu_\times^R}_{A,C}
\end{align*}
\begin{align*}
&F_\times(\mu_{A,B,C,D}); {m_\times}^{-1}_{A+ C, B+ D}; ({\nu_\times^R}_{A,C}\times {\nu_\times^R}_{B,D}); \pi^1_{F_+(A)+F_\times(C), F_+(B)+ F_\times(D)} \\
& = F_\times(\pi^1_{A,B}+\pi^1_{C,D}); {\nu_\times^R}_{B,D} \\
\\
&{\nu_\times^R}_{A\times B, C\times D}; ({q_\times}_{A,B}+ {m_\times}^{-1}_{C,D}); \mu_{F_+(A), F_+(B), F_\times(C), F_\times(D)}; \pi^1_{F_+(A)+F_\times(C), F_+(B)+ F_\times(D)} \\
& = F_\times(\pi^1_{A,B}+\pi^1_{C,D}); {\nu_\times^R}_{B,D}
\end{align*}
The second condition holds similarly, by the universal properties of coproducts. 

The second condition of {\bf \footnotesize (SMLF.6)} is
\begin{align*}
&{\nu_+^R}_{A, B\times C}; F_+({\alpha_\times}^{-1}_{A,B,C}); {q_\times}_{A\times B, C}\\
&= (1_{F_\times(A)}\times {q_\times}_{B,C}); {\alpha_\times}^{-1}_{F_\times(A), F_+(B), F+(C)}; ({\nu^R_+}_{A, B}\times 1_{F_+(C)})
\end{align*}
The left-hand side and right-hand sides post-composed with the relevant projections are given by
\begin{align*}
&{\nu_+^R}_{A, B\times C};F_+({\alpha_\times}^{-1}_{A,B,C}); {q_\times}_{A\times B, C}; \pi^0_{F_+(A\times B), F_+(C)} = {\nu_+^R}_{A, B\times C}; F_+({\alpha_\times}^{-1}_{A,B,C}); F_+(\pi^0_{A\times B, C}) \\
&= {\nu_+^R}_{A, B\times C}; F_+(1_{A}\times \pi^0_{B,C}) = (1_{F_\times(A)}\times F_+(\pi^0_{B,C})); {\nu_+^R}_{A, B} \\
\\
&(1_{F_\times(A)}\times {q_\times}_{B,C}); {\alpha_\times}^{-1}_{F_\times(A), F_+(B), F+(C)}; ({\nu^R_+}_{A, B}\times 1_{F_+(C)}); \pi^0_{F_+(A\times B), F_+(C)} \\
&= (1_{F_\times(A)}\times {q_\times}_{B,C}); {\alpha_\times}^{-1}_{F_\times(A), F_+(B), F+(C)}; \pi^0_{F_\times(A)\times F_+(B), F_+(C)}; {\nu^R_+}_{A, B} \\
&= (1_{F_\times(A)}\times {q_\times}_{B,C}); (1_{F_\times(A)}\times \pi^0_{F_+(B), F_+(C)}); {\nu^R_+}_{A, B} = (1_{F_\times(A)}\times F_+(\pi^0_{B,C})); {\nu_+^R}_{A, B} \\
\\
&{\nu_+^R}_{A, B\times C}; F_+({\alpha_\times}^{-1}_{A,B,C}); {q_\times}_{A\times B, C}; \pi^1_{F_+(A\times B), F_+(C)} = {\nu_+^R}_{A, B\times C}; F_+({\alpha_\times}^{-1}_{A,B,C}); F_+(\pi^1_{A\times B, C}) \\
&= {\nu_+^R}_{A, B\times C}; F_+(\pi^1_{A, B\times C}); F_+(\pi^1_{B,C}) = {\nu_+^R}_{A, B\times C}; F_+(t_{A}\times 1_{B\times C}); F_+({u^L_\times}_{B\times C}^{-1});F_+(\pi^1_{B,C}) \\
&=(F_\times(t_{A})\times F_+(1_{B\times C})); {\nu_+^R}_{\bone, B\times C};F_+({u^L_\times}_{B\times C}^{-1});F_+(\pi^1_{B,C})\\
&=(F_\times(t_{A})\times F_+(1_{B\times C})); (m_{\bone}^{-1}\times 1_{F_+(B\times C)}); {u^L_\times}_{F_+(B\times C)}^{-1}; F_+(\pi^1_{B,C}) =\pi^1_{F_\times(A), F_+(B\times C)}; F_+(\pi^1_{B,C})\\
\\
&(1_{F_\times(A)}\times {q_\times}_{B,C}); {\alpha_\times}^{-1}_{F_\times(A), F_+(B), F+(C)}; ({\nu^R_+}_{A, B}\times 1_{F_+(C)}); \pi^1_{F_+(A\times B), F_+(C)} \\
&= (1_{F_\times(A)}\times {q_\times}_{B,C}); {\alpha_\times}^{-1}_{F_\times(A), F_+(B), F+(C)}; \pi^{1}_{F_\times(A)\times F_+(B), F_+(C)} \\
&= (1_{F_\times(A)}\times {q_\times}_{B,C}); \pi^{1}_{F_\times(A), F+(B)\times F_+(C)}; \pi^1_{F_+(B), F_+(C)} \\
&=  \pi^1_{F_\times(A), F_+(B\times C)}; {q_\times}_{B,C}; \pi^1_{F_+(B), F_+(C)} = \pi^1_{F_\times(A), F_+(B\times C)}; F_+(\pi^1_{B,C}) 
\end{align*}
Then, by the universal property of products, the second condition of {\bf \footnotesize (SMLF.6)}  holds. The first condition of {\bf \footnotesize (SMLF.6)} follows by the universal properties of coproducts for the same reasons as its counterpart.

Consider a linear transformation $\alpha = (\alpha_\times,\alpha_+)\c F\Rarr G\c\bX\rarr\bY$ between cartesian linear functors. It consists of a monoidal transformations $\alpha_\times\colon (F_\times, m^{F}_{\bone}, m^{F}_\times)\Rarr (G_\times, m^{G}_{\bone}, m^{G}_\times)$ and $\alpha_+\colon (F_+, n^{F}_{\bzero}, n^{F}_+)\Rarr (G_+, n^{G}_{\bzero}, n^{G}_+)$. By Proposition \ref{lem:functor_comon_if_cart}, $\alpha_\times\colon(F, p^{F}_{\bzero}, p^{F}_+)\Rarr(G, p^{G}_{\bzero}, p^{G}_+)$ and $\alpha_+\colon (F_+, q^{F}_{\bone}, q^{F}_\times)\Rarr (G_+, q^{G}_{\bone}, q^{G}_\times)$ are monoidal as well. As such, $\alpha$ is a medial linear transformation.

Therefore, the inclusion map $\inc\c\CLDC\rarr\SMLDC$ is well-defined. The proof that it is moreover a 2-functor follows by examining that the canonical medial linear functor structure of a cartesian linear functor described above respects the composition and the identities. 
\end{proof}

\subsection{Frobenius Medial Linear Functors and Transformations} \hfill\

\vspace{0.5\baselineskip}
As in the case of Frobenius linear functors, we can consider symmetric medial linear functors whose pair of functors $F_\ot$ and $F_\op$ are equal.
\begin{definition}
A symmetric medial linear functor $F=(F_\ot,F_\op)\c \bX\rarr\bY$ is {\bf Frobenius} if
\[ F_\ot = F_\op, \qquad {\nu^R_\ot}_{A,B} = {n_\op}_{A,B}, \qquad {\nu^R_\op}_{A,B} = {m_\ot}_{A,B}, \]
\[ m_\bot = n_\bot^{-1}, \qquad {m_\op}_{A,B} = {n_\op}^{-1}_{A,B}, \qquad n_\top = m_\top^{-1}, \quad\text{and}\quad {n_\ot}_{A,B} = {m_\ot}^{-1}_{A,B}\]
\end{definition}

\begin{lemma}
Every Frobenius symmetric medial linear functor $F=(F,F)\c \bX\rarr\bY$ is mix.
\end{lemma}
\begin{proof}
In the context of Frobenius medial linear functors, the third condition of \mbox{\bf\footnotesize{(SMLF.1)}}, or of \mbox{\bf\footnotesize{(SMLF.2)}}, coincides with the additional condition \mbox{\bf\footnotesize{(MixFLF)}}.
\end{proof}

Given the degeneracy, we can once again give an alternative characterization.
\begin{proposition}
Let $\bX$ and $\bY$ be SMLDCs, then the following notions coincide:
\begin{itemize}
    \item Frobenius symmetric medial linear functors $F = (F_\ot, F_\op)\c \bX\rarr\bY$, and 
    \item symmetric $\ot$- and $\op$-monoidal functors
    $(F, m_\top,m_\ot, n_\bot, n_\op)\c \bX\rarr\bY$ satisfying \eqref{cc:Frobenius_linear_functor}, \mbox{\bf\footnotesize{(MixFLF)}} and the following diagrams.
\end{itemize}
\begin{equation}\begin{gathered}\label{cc:Frobenius_medial_linear_functor_nabla_delta}
\xymatrixrowsep{1.75pc}\xymatrixcolsep{1.25pc}\xymatrix{
F(\bot)\ar[r]^-{n_\bot}\ar[d]_-{F(\Delta_\bot)} & \bot\ar[r]^-{\Delta_\bot} & \bot\ot\bot\ar[d]^-{n_\bot^{-1}\ot n_\bot^{-1}} & F(\top\op\top)\ar[d]_-{F(\nabla_\top)}\ar[rr]^-{{n_\op}_{\top,\top}} & & F(\top) \op F(\top)\ar[d]^-{m_\top^{-1} \op m_\top^{-1}}\\
F(\bot\ot\bot) & & F(\bot)\ot F(\bot) \ar[ll]^-{{m_\ot}_{\bot,\bot}} & F(\top) & \top \ar[l]^-{m_\top}& \top \op \top \ar[l]^-{\nabla_\top} 
}\\
\xymatrixrowsep{1.75pc}\xymatrixcolsep{5.25pc}\xymatrix{
F((A\ot B)\op (C\ot D))\ar[r]^-{F(\mu_{A,B,C,D})} \ar[d]_-{{n_\op}_{A\ot B, C\ot D}} & F ((A\op C)\ot (B\op D)) \ar[d]^-{{m_\ot}^{-1}_{A\op C, B\op D}} \\
F(A\ot B)\op F(C\ot D)\ar[d]_-{{m_\ot}^{-1}_{A,B}\op {m^{-1}_\ot}_{C,D}} & F(A\op C)\ot F(B\op D)\ar[d]^-{{n_\op}_{A,C}\ot {n_\op}_{B,D}} \\
(F(A)\ot F(B))\op (F(C)\ot F(D)) \ar[r]_*+<0.75em>{_{\mu_{F(A), F(B), F(C), F(D)}}} &  (F(A)\op F(C))\ot (F(B)\op F(D)) \\
}
\end{gathered}\end{equation}
\end{proposition}

If we then consider medial linear transformations between such linear functors, we can see that the paired transformations must in fact be a section-retraction pair.

\begin{proposition}
Consider a medial linear transformation $\alpha = (\alpha_\ot, \alpha_\op)\c F\Rarr G$ between Frobenius symmetric medial linear functors $F, G\c \bX\rarr\bY$, then 
\begin{enumerate}[(i)]
    \item $\alpha_\ot$ and $\alpha_\op$ are both $\ot$-monoidal and $\op$-monoidal transformations, and 
    \item ${\alpha_\ot}_A\c F(A)\rarr G(A)$ is a section of ${\alpha_\op}_A\c G(A)\rarr F(A)$, i.e. ${\alpha_\ot}; \alpha_\op = 1_{F}$.
\end{enumerate}
\end{proposition} 
\begin{proof}
Let $\alpha = (\alpha_\ot, \alpha_\op)\c F\Rarr G$ be a medial linear transformation between Frobenius symmetric medial linear functors. This means in particular that $\alpha_\ot\c (F, m_\top^F, m_\ot^F)\Rarr (G, m_\top^G, m_\ot^G)$ and $\alpha_\ot\c (F, {n_\bot^F}^{-1}, {n_\op^F}^{-1})\Rarr (G, {n_\bot^G}^{-1}, {n_\op^G}^{-1})$ are monoidal transformations. The latter implies $\alpha_\ot\c (F, n_\bot^F, n_\op^F)\Rarr (G, n_\bot^G, n_\op^G)$ is a monoidal transformation. Similarly, for $\alpha_\op$. 

Further, as $\alpha=(\alpha_\ot, \alpha_\op)\c F\Rarr G$ is a linear transformation, \mbox{\bf\footnotesize{(LT)}} holds, the first of which is as follows in the context of Frobenius linear functors:
\[ {\alpha_\ot}_{A\op B}; {n_\op^G}_{A,B}; ({\alpha_\op}_A\op 1_{G(B)}) = {n_\op^F}_{A,B}; (1_{F(A)}\op {\alpha_\ot}_B)\qquad (*)\]
Letting $B=\bot$, we get the following commuting diagram.
\begin{equation*}
\xymatrixrowsep{2.75pc}\xymatrixcolsep{3.25pc}\xymatrix{
& G(A)\ar[rd]^-{G({u^R_\op}_A^{-1})}\ar@{}[d]|{(\nat)} & & \\
F(A)\ar[ru]^-{{\alpha_\ot}_A}\ar[r]_-{F({u^R_\op}_A^{-1})}\ar[dd]_-{{u^R_\op}^{-1}_{F(A)}}\ar@{}[rdd]|{(\mathrm{\bf MT})} & F(A\op \bot)\ar@{}[rdd]|{(*)}\ar[dd]_-{{n_\op^F}_{A,\bot}}\ar[r]_-{{\alpha_\ot}_{A\op \bot}} & G(A\op \bot)\ar[r]^-{G({u^R_\op}_A)}\ar[d]_-{{n_\op^G}_{A,\bot}}\ar@{}[rd]|{(\mathrm{\bf MT})} & G(A)\ar[rd]^-{{\alpha_\op}_A}\ar[d]_-{{u^R_\op}^{-1}_{G(A)}} \\
& & G(A)\op G(\bot)\ar[d]_-{{\alpha_\op}_A\op 1_{G(\bot)}}\ar[r]^-{1_{G(A)}\op {n^G_\bot}^{-1}}\ar@{}[rd]|{(\exch)} & G(A)\op\bot\ar[d]_-{{\alpha_\op}_A\op 1_\bot}\ar@{}[r]|{(\nat)} & F(A) \\
F(A)\op\bot\ar[r]_-{1_{F(A)}\op {n_\bot^F}^{-1}}&F(A)\op F(\bot)\ar[r]_-{1_{F(A)}\op {\alpha_\ot}_\bot} & F(A)\op G(\bot)\ar[r]_-{1_{F(A)}\op n^G_\bot} & F(A)\op\bot\ar[ru]_-{{u^R_\op}_{F(A)}}
}
\end{equation*}
The top composite is ${\alpha_\ot}_A; {\alpha_\op}_A$, while the bottom composite is $1_{F(A)}$ since $\alpha_\ot$ is $\op$-monoidal. As such, ${\alpha_\ot}; \alpha_\op = 1_{F}$.
\end{proof}

An immediate situation where one gets such a linear transformation is when considering an invertible $\ot$- and $\op$-monoidal natural transformation.
\begin{proposition}
Given a $\ot$- and $\op$-monoidal isomorphism $\alpha\c F\Rarr G$ between Frobenius symmetric medial linear functors, the pair $(\alpha, \alpha^{-1})\c F\Rarr G$ is a medial linear transformation. 
\end{proposition}
\begin{proof}
Consider an invertible $\ot$-monoidal and $\op$-monoidal natural transformation $\alpha\c F\Rarr G$. By Lemma \ref{lem:Frobenius_linear_trans_2}, $(\alpha, \alpha^{-1})\c F\Rarr G$ is a linear transformation. To be medial, it remains to show $\alpha\c (F_\ot, m_\bot^F, m_\op^F)\Rarr (G_\ot, m_\bot^G, m_\op^G)$ and $\alpha^{-1}\c (G_\op, n_\top^G, n_\ot^G)\Rarr (F_\op, n_\top^F, n_\ot^F)$ are monoidal. Since $F$ and $G$ are Frobenius, this is in fact requiring that $\alpha\c (F, {n_\bot^F}^{-1}, {n_\op^F}^{-1})\Rarr (G, {n_\bot^G}^{-1}, {n_\op^G}^{-1})$ and $\alpha^{-1}\c (G, {m_\top^G}^{-1}, {m_\ot^G}^{-1})\Rarr(F, {m_\top^F}^{-1}, {m_\ot^F}^{-1})$ are monoidal. This is immediate as $\alpha\c (F, n_\bot^F, n_\op^F)\Rarr (G, n_\bot^G, n_\op^G)$ and $\alpha\c (F, m_\top^F, m_\ot^F)\Rarr (G, m_\top^G, m_\ot^G)$ are monoidal.
\end{proof}

Consider a Frobenius symmetric medial linear functor, in the sense of the above Proposition. Recall that \eqref{cc:Frobenius_linear_functor} is equivalent to requiring that $(F, m_\ot, m_\top, n_\op, n_\bot)\c \bX\rarr\bY$ is a Frobenius linear functor. Moreover, \mbox{\bf\footnotesize{(MixFLF)}} and \eqref{cc:Frobenius_medial_linear_functor_nabla_delta} together are equivalent to asking that $(F, n_\bot^{-1}, n_\op^{-1}, m_\top^{-1}, m_\ot^{-1})\c (\bX, \op, \bot, \ot, \top)\rarr(\bY, \op, \bot, \ot, \top)$ is a bilax duoidal functor.  As Frobenius linear functors and bilax duoidal functors compose, we know that Frobenius medial linear functors are closed under composition. It is thus fairly immediate to see that we can restrict to Frobenius linear functors.

\begin{corollary}
SMLDCs, Frobenius symmetric medial linear functors and medial linear transformations form a sub 2-category of \SMLDC, denoted by $\mathbf{FSMLDC}$. Further, the inclusion map $\inc\c \CLDC\rarr\SMLDC$ restricts to an inclusion 2-functor $\inc\c\mathbf{FCLDC}\rarr\mathbf{FSMLDC}$
\end{corollary}
\begin{proof}
The only technicality to verify is that given a Frobenius cartesian linear functor $F=(F_\times,F_+)\c \bX\rarr \bY$, the induced medial linear functor is Frobenius, in other words that the canonical lax monoidal structure $(F_\times,p_\bzero, p_+)\c (\bX, \bzero, +)\rarr (\bY, \bzero, +)$ corresponds to $(F_+, n_\bzero^{-1}, n_+^{-1})\c (\bX, \bzero, +)\rarr (\bY, \bzero, +)$ and, similarly, for the canonical colax momonoidal structure on $F_+$. This is immediate however as $F_+ = F_\times$. Indeed, while the structure maps $n_\bzero$ and ${n_+}_{A,B}$ are not guaranteed to have inverses, if they do, they must be the canonical ones $p_\bzero$ and ${p_+}_{A,B}$. Similarly, $q_\bone = {m_\bone}^{-1}$ and ${q_\times}_{A,B} = {m_\times^{-1}}_{A,B}$.
\end{proof}

\section{Linearly Distributive Fox Theorem}\label{sec:LD_Fox}

It is now time to discuss the construction mapping a SMLDC to a CLDC, mirroring the category of cocommutative comonoids of a SMC in the standard Fox theorem. 

If we review our earlier characterization of CLDCs in Proposition \ref{prop:char_CLDC}, we see that a CLDC \bX\ is in particular a SLDC, where each object $A\in\bX$ is equipped with maps
\[\Delta_A\c A\rarr A\ot A \qquad e_A\c A\rarr\top \qquad \nabla_A\c A\op A\rarr A \qquad u_A\c \bot\rarr A\]
such that $\langle A, \Delta_A, e_A\rangle$ is a cocommutative $\ot$-comonoid and $\langle A, \nabla_A, u_A\rangle$ is a commutative $\op$-monoid. To someone well-versed in duoidal categories, this structure will be familiar. Indeed, these are the same maps and structures of a duoidal bimonoid. This is the key to the theorem and why the duoidal structure of the MLDCs is essential.

\subsection{Medial Bimonoids}\label{sec:medial_bimonoids}  \hfill\

\begin{definition}
Let $\bX$ be a SMLDC. A {\bf medial bimonoid} in \bX\ is a quintuple \\$\langle A, \Delta_{A}, u_{A}, \nabla_{A}, e_{ A}\rangle$ consisting of an object $A$ and four maps in \bX,
\[\Delta_A\c A\rarr A\ot A \qquad e_A\c A\rarr\top \qquad \nabla_A\c A\op A\rarr A \qquad u_A\c \bot\rarr A\] such that $\langle A, \Delta_{A}, e_{A}\rangle$ is a $\ot$-comonoid, $\langle A, \nabla_{A}, u_{A}\rangle$ is a $\op$-monoid, and 
satisfying
\begin{equation}\begin{gathered}\label{cc:medial_bimonoid}\tag{\bf\footnotesize{MB}}
\xymatrixrowsep{1.75pc}\xymatrixcolsep{1.75pc}\xymatrix@L=0.5pc{
A\op A\ar[r]^-{\nabla_A}\ar[d]_-{\Delta_A\op\Delta_A} & A\ar[r]^-{\Delta_A} & A\ot A  & \bot\ar[r]^-{m}\ar[d]_-{u_A} & \top\\
(A\ot A)\op (A\ot A)\ar[rr]_-{\mu_{A,A,A,A}} && (A\op A)\ot (A\op A)\ar[u]_-{\nabla_A\ot\nabla_A} & A\ar[ru]_-{e_A} \\
A\op A\ar[r]^{\nabla_A}\ar[d]_-{e_A\op e_A} & A\ar[d]^-{e_A} & \bot\ar[r]^-{u_A}\ar[d]_-{\Delta_\bot} & A\ar[d]^-{\Delta_A} \\
\top \op \top \ar[r]_-{\nabla_\top}& \top & \bot\ot \bot\ar[r]_-{u_A\ot u_A} & A\ot A
}
\end{gathered}\end{equation}
A medial bimonoid is {\bf bicommutative} if $\langle A, \Delta_{A}, e_{A}\rangle$ is cocommutative $\ot$-comonoid and $\langle A, \nabla_{A}, u_{A}\rangle$ is a commutative $\op$-monoid. 

A {\bf medial bimonoid morphism} $f\c\langle A, \Delta_{A}, e_{A}, \nabla_{A}, u_{A}\rangle\rarr\langle B, \Delta_{B}, e_{B}, \nabla_{B}, u_{B}\rangle$ is a morphism $f\c A\rarr B$ such that $f\c\langle A, \Delta_{A}, e_{A}\rangle\rarr\langle B, \Delta_{B}, e_{B}\rangle$ is a $\ot$-comonoid morphism and $f\c\langle A, \nabla_{A}, u_{A}\rangle\rarr\langle B, \nabla_{B}, u_{B}\rangle$ is a $\op$-monoid morphism.
\end{definition}

\begin{remark}\label{rem:medial_bimonoid_bilax_functors}
Let $(\bX, \ot, \top, \op, \bot)$ be a SMLDC. Recall that $(\bX, \op, \bot, \ot, \top)$ is a duoidal category and Remark \ref{rem:duoidal_bimonoid}. It is then immediate that the following are equivalent:
\begin{itemize}
    \item a medial bimonoid $\langle A, \Delta_{A}, e_{A}, \nabla_{A}, u_{A}\rangle$ in \bX,
    \item a duoidal bimonoid $\langle A, \nabla_{A}, u_{A}, \Delta_{A}, e_{A}\rangle$ in \bX, and 
    \item a bilax duoidal functor from the terminal category \[ (A, u_A, \nabla_A, e_A, \Delta_A)\c 1\rarr \bX \]
\end{itemize}
Moreover, the following are also equivalent:
\begin{itemize}
    \item a medial bimonoid morphism $f\c\langle A, \Delta_{A}, e_{A}, \nabla_{A}, u_{A}\rangle\rarr\langle B, \Delta_{B}, e_{B}, \nabla_{B}, u_{B}\rangle$ in \bX, 
    \item a morphism of duoidal bimonoids $f\c\langle A,\nabla_{A}, u_{A}, \Delta_{A}, e_{A}\rangle\rarr\langle B, \nabla_{B}, u_{B}, \Delta_{B}, e_{B}\rangle$ in \bX, and 
    \item a duoidal transformation between the bilax duoidal functors
\[ f\c (A, u_A, \nabla_A, e_A, \Delta_A)\Rarr(B, u_B, \nabla_B, e_B, \Delta_B) \c 1\rarr \bX\]
\end{itemize}
We will frequently move between these different perspectives and rely upon the 2-categorical nature of \DUO\ in the proofs that follow, in order to avoid displaying large amounts of commuting diagrams. 

Note, however, that all of the arguments could instead be proved by straightforward diagram chasing, from the various coherence conditions for SMLDCs, bicommutative medial bimonoids, and their morphisms.
\end{remark}

It is immediate that the $\ot$ and $\op$ units in a SMLDC have canonical medial bimonoid structures, since, by Proposition \ref{prop:monoidal_products_bilax_duoidal_functor}, the units of a symmetric duoidal category are bilax duoidal functors from the terminal category.

\begin{proposition}\label{prop:top_bottom_medial_bimonoids}
Given a SMLDC \bX, the following structures
\[\langle \top, {u_\ot}_\top\c \top\rarr\top\ot\top, 1_\top\c\top\rarr\top, \nabla_\top\c\top\op\top\rarr\top, m\c \bot\rarr\top \rangle\]
\[ \langle \bot, \Delta_\bot\c\bot\rarr\bot\ot\bot, m\c\bot\rarr\top, {u_\op}_\bot\c\bot\op\bot\rarr\bot, 1_\bot\c\bot\rarr\bot\rangle\]
are bicommutative medial bimonoids. 
\end{proposition}

Furthermore, medial bimonoids combine together via the tensor and par monoidal products to produce new medial bimonoids. 

\begin{proposition}\label{prop:constructing_medial_bimonoids}
Consider a pair of bicommutative medial bimonoids $\langle A, \Delta_A, e_A, \nabla_A, u_A\rangle $ and $\langle B, \Delta_B, e_B, \nabla_B, u_B\rangle$ in a SMLDC \bX. Then, the object $A\ot B$ equipped with morphisms 
\begin{align*}
& \Delta_{A\ot B} = A\ot B \xrightarrow{\Delta_A\ot\Delta_B} (A\ot A)\ot (B\ot B) \xrightarrow{\tau^\ot_{A,A,B,B}} (A\ot B)\ot(A\ot B) \\
&t_{A\ot B} = A\ot B\xrightarrow{e_A\ot e_B} \top\ot \top \xrightarrow{{u_\ot}_\top^{-1}} \top \\
&\nabla_{A \ot B} = (A\ot B)\op(A\ot B)\xrightarrow{\mu_{A,B,A,B}} (A\op A)\ot (B\op B) \xrightarrow{\nabla_A\ot \nabla_B} A\ot B \\
&u_{A\ot B}= \bot\xrightarrow{\Delta_\bot} \bot\ot\bot \xrightarrow{u_A\ot u_B} A\ot B
\end{align*}
and the object $A\op B$ equipped with morphisms 
\begin{align*}
& \Delta_{A\op B} = A\op B \xrightarrow{\Delta_A\op\Delta_B} (A\ot A)\op (B\ot B) \xrightarrow{\mu_{A,A,B,B}} (A\op B)\ot(A\op B) \\
&t_{A\op B} = A\op B\xrightarrow{e_A\op e_B} \top\op\top \xrightarrow{\nabla_\top} \top \\
&\nabla_{A\op B} = (A\op B)\op(A\op B)\xrightarrow{\tau^\op_{A,B,A,B}} (A\op A)\op (B\op B) \xrightarrow{\nabla_A\op \nabla_B} A\op B \\
&u_{A\op B}= \bot\xrightarrow{{u_\op}_\bot^{-1}} \bot\op\bot \xrightarrow{u_A\op u_B} A\op B
\end{align*}
are bicommutative medial bimonoids. 
\end{proposition}
\begin{proof}
Given $\langle A, \Delta_A, e_A, \nabla_A, u_A\rangle $ and $\langle B, \Delta_B, e_B, \nabla_B, u_B\rangle$, these structures paired together determine a duoidal bimonoid in the duoidal category $\bX\times\bX$.  Recall that $\ot$ and $\op$ are bilax duoidal functors $\bX\times\bX\rarr\bX$ by Proposition \ref{prop:monoidal_products_bilax_duoidal_functor}. As such, by Proposition \ref{prop:duoidal_functor_preserve}, when $\ot$ and $\op$ are applied to $(A,B)$ in $\bX\times\bX$, the paired duoidal structure is preserved and we get the following duoidal bimonoids in \bX:
\[\langle A\ot B, \Delta_{A\ot B}, e_{A\ot B}, \nabla_{A\ot B}, u_{A\ot B}\rangle\qquad \langle A\op B, \Delta_{A\op B}, e_{A\op B}, \nabla_{A\op B}, u_{A\op B}\rangle \]

It remains only to demonstrate they are bicommutative. It follows immediately, by Lemma \ref{lem:constructing_comonoids} that $\langle A\ot B, \Delta_{A\ot B}, e_{A\ot B}\rangle$ is cocommutative. $\langle A\ot B, \nabla_{A\ot B}, u_{A\ot B}\rangle$ is a commutative $\op$-monoid by the commuting diagram below.
\begin{equation*}
\xymatrixrowsep{2.75pc}\xymatrixcolsep{5pc}\xymatrix{
(A\ot B)\op (A\ot B)\ar[r]^-{\mu_{A,B,A,B}}\ar[d]_-{{\sigma_\op}_{A\ot B, A\ot B}}\ar@{}[rd]|{(\mathrm{\bf BDUO})} & (A\op A)\ot (B\op B)\ar[r]^-{\nabla_A\ot \nabla_B}\ar[d]_-{{\sigma_\op}_{A,A}\ot {\sigma_\op}_{B,B}} \ar@{}[rd]|{(\mathrm{\bf CM})_{A,B}} & A\ot B\\
(A\ot B)\op (A\ot B)\ar[r]_-{\mu_{A,B,A,B}} & (A\op A)\ot (B\op B)\ar@/_2pc/[ru]_-{\quad\nabla_A\ot \nabla_B} &
}
\end{equation*}
The proof that $A\op B$ is bicommutative follows similarly. 
\end{proof}

Let us for a moment consider our most important class of SMLDCs, the CLDCs, and discuss the medial bimonoids present in these cases.

\begin{proposition}\label{prop:canonical_medial_bimonoids_in_CLDC}
Every object $A$ in a CLDC \bX\ has a unique canonical bicommutative medial bimonoid structure 
\[ \langle A, \langle 1_A, 1_A\rangle\c A\rarr A\times A, t_A\c A\rarr \bone, [1_A, 1_A]\c A+A\rarr A, b_A\c \bzero\rarr A\rangle \]
and every arrow $f\c A\rarr B$ in \bX\ is a medial bimonoid morphism with respect to them.
\end{proposition}
\begin{proof}
By Lemma \ref{lem:unique_comonoid_cartesian_cat}, it is clear $\langle A, \langle 1_A, 1_A\rangle, t_A\rangle$ is the unique cocommutative $\times$-comonoid structure on $A$ and, dually, $\langle A, [1_A, 1_A], b_A\rangle$ is the unique commutative $+$-monoid structure on $A$. The proof that the structure maps satisfy \mbox{\bf\footnotesize{(MB)}} is a simple exercise using the canonical medial structure of a CLDC and universal properties. Once more, by Lemma \ref{lem:unique_comonoid_cartesian_cat} and its dual, every arrow is a medial bimonoid morphism with respect to these structures. 
\end{proof}

\begin{definition}
Let \bX\ be a SMLDC. Define $\B[\bX]$ to be the category of bicommutative medial bimonoids and medial bimonoid morphisms in \bX. 
\end{definition}

The fact that this is in fact a category is immediate, as the category of commutative $\op$-monoids and of cocommutative $\ot$-comonoids in \bX\ are themselves be categories. However, $\B[\bX]$ is more than a simple category. It is a SLDC.

\begin{lemma}
Let \bX\ be a SMLDC, then $\B[\bX]$ is a SLDC, with $\ot$ and $\op$ monoidal structures as defined in Propositions \ref{prop:top_bottom_medial_bimonoids} and \ref{prop:constructing_medial_bimonoids}.
\end{lemma}
\begin{proof}
Firstly, consider bicommutative medial bimonoids structures on objects $A, B, A', B'$ in \bX, and medial bimonoid morphisms $f\c A\rarr B$ and $g\c A'\rarr B'$. Viewed as morphisms of duoidal bimonoids, the pair $(f,g)$ is once again such a morphism in the duoidal category $\bX\times\bX$. Once more, by Proposition \ref{prop:monoidal_products_bilax_duoidal_functor} and Proposition \ref{prop:duoidal_functor_preserve}, we see that the bilax duoidal functors $\ot$ and $\op$ map $(f,g)$ to morphisms of duoidal bimonoids
\[ f\ot g\c A\ot B\rarr A'\ot B' \qquad f\op g\c A\op B\rarr A'\op B'\]
Therefore, $\ot$ and $\op$ are inherited in $\B[\bX]$ from \bX. 

Now, consider bicommutative medial bimonoid structures on $A, B, C$ in \bX. By Proposition \ref{prop:interaction_canonical_flip_medial_maps_linear_dist}, the left linear distributor is a duoidal transformation. Now, once more by Proposition \ref{prop:duoidal_functor_preserve}, we see that the component map $\delta^L_{A,B,C}$ is a morphism of duoidal bimonoids. As such, $\B[\bX]$ inherits $\delta^L$ from \bX.  

The proofs that the other structure maps, the associators, the unitors and the braidings, inherited by \bX, follows by the same argument, recalling Proposition \ref{prop:monoidal_isomorphisms_duoidal_trans} which proves these maps are duoidal transformations.
\end{proof}

\begin{theorem}\label{thm:BX_CLDC}
Let \bX\ be a SMLDC then $\B[\bX]$ is a CLDC.    
\end{theorem}
\begin{proof}
$\B[\bX]$ satisfies the characterization in Proposition \ref{prop:char_CLDC}. Given a bicommutative medial bimonoid $\langle A, \Delta_A, e_A, \nabla_A, u_A\rangle$, then the four maps $\Delta_A$, $e_A$, $\nabla_A$ and $u_A$ are themselves medial bimonoid morphisms as follows. By viewing the medial bimonoid as a bilax duoidal functor $(A, u_A, \nabla_A, e_A, \Delta_A)\c 1\rarr\bX$, the structure maps of such functors are duoidal transformations by Proposition \ref{prop:duoidal_functor_maps_duoidal_trans}, and the component maps of duoidal transformations are morphisms of duoidal bimonoids by Proposition \ref{prop:duoidal_functor_preserve}. These are precisely the component maps of the four natural transformations required in Proposition \ref{prop:char_CLDC}.

The conditions that the component maps form comonoid and monoid structures is immediate by the definition of a bicommutative medial bimonoid and the remaining coherence conditions are satisfied as they correspond to precisely to the construction of the tensor and par monoidal structures in $\B[\bX]$ as outlined in Propositions \ref{prop:top_bottom_medial_bimonoids} and \ref{prop:constructing_medial_bimonoids}.
\end{proof}

\begin{example}
\begin{enumerate}[(i)]
    \item Consider a SMC $(\cX, \os, I)$, viewed as a SMLDC $(\cX, \os, I, \os, I)$ by Proposition \ref{prop:braided_monoidal_MLDC}. Then, it is immediate that a bicommutative medial bimonoid in $\cX$ is a bicommutative bimonoid, or bialgebra, in the traditional algebraic sense. As such, $\B[\cX]$ is the category of bicommutative bimonoids and bimonoid morphisms in \cX, where the inherited monoidal product is the biproduct, and therefore an example of a compact CLDC.

    \item Consider \PCOH, as defined in Definition \ref{def:PCOH}. It is fairly straightforward to write what conditions must be satisfied in \PCOH\ to obtain bicommutative medial bimonoids and medial bimonoid maps, although it is not particularly instructive. As it stands, $\B[\PCOH]$ does not seem to provide a nice or recognizable CLDC, even for simple choices of $\mathbb{P}$. In fact, much of the work done in \cite[Sec 4.12]{Lamarche_2007} is to build particular classes of bimonoids in \PCOH, for $\mathbb{P}=\{0,1\}$. Note that the classes of bimonoids studied are somewhat different than our bicommutative medial bimonoids. Nevertheless, we can adapt some of the results to demonstrate that $\B[\PCOH]$ is non-trivial. 

    Firstly, $\B[\PCOH]$ is far from non-empty, since by \cite[Prop 4.13]{Lamarche_2007}, the following are bicommutative medial bimonoids in \PCOH:
    \begin{definition}
    Consider a poset $A = (|A|, \sqsubseteq)$, equipped with either the {\bf discrete structure}, $\rho_A(a,a') = \bot$, or the {\bf codiscrete structure}, $\rho_A(a,a') = \top$. Define the following $\mathbb{P}$-coherence maps 
    \begin{align*}
        &\Delta_A\c A\rarr A\ot A & (a_0,(a_1,a_2))\in\Delta_A \iff a_0\sqsubseteq a_1\wedge a_2 \\
        &e_A\c A\rarr\top & (a,*)\in e_A \quad\forall a\in|A| \\
        &\nabla_A\c A\op A\rarr A & ((a_0,a_1),a_2)\in\nabla_A \iff a_0\vee a_1\sqsubseteq a_2 \\
        &u_A\c \bot\rarr A & (*,a)\in u_A \quad\forall a\in|A|       
    \end{align*}
    This is a bicommutative medial bimonoid, known as the {\bf diagonal-is-inf model}.
    \end{definition}
    
    Secondly, if $\mathbb{P}$ is non-trivial, $\B[\PCOH]$ is not isomix as the tensor and par units are not isomorphic. A complete understanding of $\B[\PCOH]$ remains an open question, particularly whether it is posetal for all $\mathbb{P}$. 
\end{enumerate}   
\end{example}

It is also important to note that negation is preserved by the $\B[-]$ construction. 

\begin{proposition}
Let \bX\ be a SMLDC with negation, then $\B[\bX]$ is a CLDC with negation, and therefore a posetal distributive category, in the sense of \cite{Kudzman-Blais_Lemay_2025}.
\end{proposition}
\begin{proof}
It is easiest in this case to view the SMLDC with negation \bX\ as a $*$-autonomous category with a full and faithful functor $(-)^\perp\c \bX^{op}\rarr \bX$ satisfying $\bX(A\ot B, C^\perp)\cong \bX(A, (B\ot C)^\perp)$. We can show that $(-)^\perp$ extends to a functor on $\B[\bX]$. Consider a bicommutative medial bimonoid $\langle A, \Delta_A, e_A, \nabla_A, u_A\rangle$, define $\langle A^\perp, \Delta_{A^\perp}, e_{A^\perp}, \nabla_{A^\perp}, u_{A^\perp}\rangle$ by 
\begin{align*}
&\Delta_{A^\perp} = A^\perp \xrightarrow{(\nabla_A)^\perp} (A\op A)^\perp \xrightarrow{\sim} A^\perp \ot A^\perp & e_{A^\perp} = A^\perp \xrightarrow{(u_A)^\perp} \bot^\perp \xrightarrow{\sim} \top \\
& \nabla_{A^\perp} = A^\perp\op A^\perp \xrightarrow{\sim} (A\ot A)^\perp \xrightarrow{(\Delta_A)^\perp} A^\perp & u_{A^\perp} = \bot \xrightarrow{\sim} \top^\perp \xrightarrow{(e_A)^\perp} A^\perp 
\end{align*}
This is a bicommutative medial bimonoid. Moreover, given any medial bimonoid morphism $f\c \langle A, \Delta_A, e_A, \nabla_A, u_A\rangle \rarr \langle B, \Delta_B, e_B, \nabla_B, u_B\rangle$, then $f^\perp\c B^\perp\rarr A^\perp$ is a medial bimonoid morphism. The proof of the necessary coherence conditions is straightforward for anyone well-versed in categorical linear logic, but involves many technical details concerning linear negation in SLDCs/$*$-autonomous categories which have not been the focus of this article, and therefore has been excluded for readability. 

As $(-)^\perp$ extends to $\B[\bX]$ and $\B[\bX]$ inherits its monoidal structures from \bX, it is immediate that $\B[\bX]$ is a CLDC with negation. By Joyal's paradox, detailed in \cite[Sec 4.2]{Kudzman-Blais_Lemay_2025} for CLDCs, $\B[\bX]$ is a posetal distributive category.
\end{proof}

With Theorem \ref{thm:BX_CLDC}, the first part of the linearly distributive Fox theorem has been established. It remains only to show that the $\B[-]$ construction extends canonically to symmetric medial linear functors and transformations, giving a right adjoint to the inclusion functor outlined in Proposition \ref{prop:inclusion_CLDC_SMLDC}.

\subsection{The Right Adjoint 2-Functor {\ensuremath {\rm B}}[-]}\label{sec:right_adjoint}  \hfill\

\begin{lemma}\label{lem:medial_linear_functor_cartesian_linear_functor}
Consider a symmetric medial linear functor $F=(F_\ot, F_\op)\c \bX\rarr\bY$ between SMLDCs, then it canonically extends to a cartesian linear functor \[ \B[F]=(\B[F]_\ot, \B[F]_\op)\c \B[\bX]\rarr\B[\bY]\]
\end{lemma}

\begin{proof}
Consider a symmetric medial linear functor $F=(F_\ot, F_\op)\c\bX\rarr\bY$. Recall that 
\[(F_\ot,  m_{ \bot},m_\op, m_{ \top}^{-1}, m_\ot^{-1}),\,\, (F_\op, n_{\bot}^{-1}, n_\op^{-1}, n_{\top}, n_\op)\c(\bX, \op,\bot,\ot, \top)\rarr(\bY,\op,\bot,\ot,\top)\]
are bilax duoidal functors and therefore preserve medial bimonoids and medial bimonoid morphisms by Proposition \ref{prop:duoidal_functor_preserve}. Thus, given $\langle A, \Delta_{A}, e_{A}, \nabla_{A}, u_{A}\rangle$ in \bX, the object $F_\ot(A)$ equipped with
\begin{align*}
&\Delta_{F_\ot(A)} = F_\ot(A) \xrightarrow{F_\ot(\Delta_{A})} F_\ot(A\ot A) \xrightarrow{{m_\ot}^{-1}_{A,A}} F_\ot(A)\ot F_\ot(A) \\
&e_{F_\ot(A)} =F_\ot(A) \xrightarrow{F_\ot(e_{A})} F_\ot(\top) \xrightarrow{m_{\top}^{-1}} \top \\
&\nabla_{F_\ot(A)} = F_\ot(A)\op F_\ot(A) \xrightarrow{{m_\op}_{A,A}} F_\ot(A\op A) \xrightarrow{F_\ot(\nabla_{A})} F_\ot(A) \\
&u_{F_\ot(A)} = \bot \xrightarrow{m_{\bot}} F_\ot(\bot) \xrightarrow{F_\ot(u_{A})} F_\ot(A)
\end{align*}
and the object $F_\op(A)$ equipped with 
\begin{align*}
&\Delta_{F_\op(A)} = F_\op(A) \xrightarrow{F_\op(\Delta_{A})} F_\op(A\ot A)\xrightarrow{{n_\ot}_{A,A}} F_\op(A)\ot F_\op(A) \\
& e_{F_\op(A)} = F_\op(A) \xrightarrow{F_\op(e_{A})} F_\op(\top) \xrightarrow{n_{\top}} \top\\
&\nabla_{F_\op(A)} = F_\op(A) \op F_\op(A) \xrightarrow{{n_\op}^{-1}_{A,A}} F_\op(A\op A) \xrightarrow{F_\op(\nabla_{A})} F_\op(A) \\
&u_{F_\op(A)} = \bot \xrightarrow{n^{-1}_{\bot}} F_\op(\bot) \xrightarrow{F_\op(u_{A})} F_\op(A)
\end{align*}
are bicommutative medial bimonoid in \bY. Moreover, given a medial bimonoid morphism $f$ in \bX, then $F_\ot(f)$ and $F_\op(f)$ are medial bimonoid morphisms in \bY. As such, let $\B[F]_\ot, \B[F]_\op\c\B[\bX] \rarr \B[\bY]$ denote the functors mapping a bicommutative medial bimonoid on $A$ to the bicommutative medial bimonoid on $F_\ot(A)$ and $F_\op(A)$ respectively, and mapping medial bimonoid morphism $f$ to $F_\ot(f)$ and $F_\op(f)$ respectively. 

Now, as the structure maps of bilax duoidal functors $m_{ \top}^{-1}, m_{ \top}^{-1}, n_{\bot}^{-1}$ and $n_\op^{-1}$ are bilax duoidal transformations by Proposition \ref{prop:duoidal_functor_maps_duoidal_trans}, so are their inverses
\begin{align*}
& m_{\top}\c \top\rarr F_\ot(\top) & {m_\ot}_{A,B}\c F_\ot(A)\ot F_\ot(B)\rarr F_\ot(A\ot B) \\
& n_{\bot}\c F_\op(\bot)\rarr\bot & {n_\op}_{A,B}\c F_\op(A\op B)\rarr F_\op(A)\op F_\op(B)
\end{align*}
By Proposition \ref{prop:duoidal_functor_preserve}, the component maps are morphisms of the relevant duoidal bimonoids. Therefore, $\B[F]_\ot$ can inherit $m_\top$ and $m_\ot$ from $F_\ot$, making it a $\ot$-monoidal functor, and $\B[F]_\op$ can inherit $n_\bot$ and $n_\op$ from $F_\op$, making it a $\op$-monoidal functor.

Finally, by Proposition \ref{prop:interaction_canonical_flip_medial_linear_functor}, the linear strengths 
\[ {v^R_\ot}_{A,B}\c F_\ot(A\op B)\rarr F_\op(A)\op F_\ot(B) \qquad {v^R_\op}_{A,B}\c F_\ot(A)\ot F_\op(B)\rarr F_\op(A\ot B) \]
are duoidal transformations. By the same argument as above, $(\B[F]_\ot, \B[F]_\op)\c \bX\rarr\bY$ inherits the linear strengths and is a well-defined cartesian linear functor.
\end{proof}

\begin{lemma}\label{lem:medial_linear_trans_is_cartesian_linear_trans}
Given a medial linear transformations $\alpha=(\alpha_\ot, \alpha_\op)\c F\Rarr G$ between medial linear functors, it canonically extends to a linear transformation between cartesian linear functors \[ \B[\alpha] = (\B[\alpha]_\ot, \B[\alpha]_\op)\c \B[F]\Rarr\B[G]\] 
\end{lemma}
\begin{proof}
This is immediate by Remark \ref{rem:medial_linear_trans}, which states that $\alpha_\ot$ and $\alpha_\op$ are appropriately duoidal transformations, and Proposition \ref{prop:duoidal_functor_preserve}, which states that the components maps of duoidal transformations are morphisms of duoidal bimonoids. 
\end{proof}

With the previous lemmas, we can state that the medial bimonoid construction is a 2-functor between the 2-category of SMLDCs and the 2-category of CLDCs as desired.

\begin{proposition}
$\B[-]\c\SMLDC\rarr\CLDC$ is a 2-functor.
\end{proposition}
\begin{proof}
By Theorem \ref{thm:BX_CLDC}, Lemma \ref{lem:medial_linear_functor_cartesian_linear_functor} and Lemma \ref{lem:medial_linear_trans_is_cartesian_linear_trans}, $\B[-]$ provides well-defined mappings, for 0-cells, 1-cells and 2-cells respectively, from \SMLDC\ to \CLDC. The proof that this is 2-functorial follows immediately since composition in the 2-categories are defined similarly and $\B[-]_{\bbX,\bbY}\c \SMLDC(\bX, \bY)\rarr\CLDC(\bX, \bY)$ is a forgetful functor. 
\end{proof}

It is fairly immediate that given a Frobenius medial linear functor $F = (F_\ot, F_\op)\c \bX\rarr\bY$, the induced cartesian linear functor $\B[F]\c \B[\bX]\rarr \B[\bY]$ is itself Frobenius as the medial bimonoid structure maps of $F_\ot(A)$ are equal to the ones of $F_\op(A)$ in the Frobenius context. Therefore:
\begin{corollary}
$\B[-]$ restricts to a 2-functor $\B[-]\c\mathbf{SFMLDC}\rarr\mathbf{FCLDC}$.  
\end{corollary}

We have finally arrived at the main result of our paper: $\B[-]$ is the right adjoint to the inclusion 2-functor $\inc \vdash \B[-]\colon\CLDC\rarr\SMLDC$.  

\begin{theorem}\label{thm:Foxthm}\emph{({\bf Linearly Distributive Fox Theorem})}
$\inc \vdash \B[-]\colon\CLDC\rarr\SMLDC$ is a 2-adjunction.
\end{theorem}

\begin{proof}
To prove this 2-adjunction, we will provide the unit and counit 2-transformations. 

Firstly, let us define the unit $\boldsymbol{\eta}\c 1_{\CLDC}\Rarr\inc;\B[-]$. Let $(\bX, \times, \bone, +, \bzero)$ be a CLDC. Recall Proposition \ref{prop:canonical_medial_bimonoids_in_CLDC}, which states that every object $A\in \bX$ has a canonical bicommutative medial bimonoid structure 
\[ \langle A, \langle 1_A, 1_A\rangle\c A\rarr A\times A, t_A\c A\rarr \bone, [1_A, 1_A]\c A+A\rarr A, b_A\c \bzero\rarr A\rangle \]
and every arrow $f\c A\rarr B \in \bX$ is a medial bimonoid morphism with respect to these structures. Then, let $\eta_{\bbX}\c\bX\rarr\B[\bX]$ be the functor mapping objects in \bX\ to their canonical medial bimonoid structures and which is identity on the arrows. 

$\eta_{\bbX}$ is a symmetric monoidal functor between $(\bX,\times, \bone)\rarr (\B[\bX], \times, \langle \bone \rangle)$ when equipped with 
\[ m_\bone = 1_\bone\c \bone\rarr \bone \qquad {m_\times}_{A, B} = 1_{A\times B}\c A\times B \rarr A\times B\]
as the canonical medial bimonoid structure on $\bone$ described above is the same as as the medial bimonoid structure on $\bone$ when \bX\ is viewed as MLDC,
    \[\langle \bone, {u_\times}_\bone, 1_\bone, \nabla_\bone, m\rangle = \langle \bone, \langle 1_\bone, 1_\bone\rangle, t_{\bone}, [1_\bone, 1_\bone], b_{\bone}\rangle\]
and the canonical medial bimonoid structure on $A\times B$ is the same as the medial bimonoid structure of the tensor product $\times$ in $\B[\bX]$ between the canonical structures on $A$ and $B$, 
\begin{align*}
\langle A\times B, & ([1_A, 1_A]\times [1_B, 1_B]); \tau^\times_{A,A,B,B}, (t_A\times t_B); {u_\times}_\bone^{-1}, \\&\mu_{A,B,A,B}; ([1_A, 1_A]\times [1_B, 1_B]), \Delta_\bot; (b_A\times b_B)\rangle \\
&= \langle A\times B, \langle 1_{A\times B}, 1_{A\times B}\rangle, t_{A\times B}, [1_{A\times B}, 1_{A\times B}], b_{A\times B}\rangle
\end{align*}

Similarly, $\eta_X$ is a symmetric monoidal functor between $(\bX,+, \bzero)\rarr (\B[\bX], +, \langle \bzero \rangle)$ when equipped with
\[ n_\bzero =1_\bzero\c \bzero\rarr \bzero \qquad  {n_+}_{A, B} = 1_{A+ B}\c A+B\rarr A+B\]
Then, $\boldsymbol{\eta}_{\bbX} = (\eta_{\bbX}, \eta_{\bbX})\c \bX\rarr \B[\bX]$ is a Frobenius cartesian linear functor, when equipped with identity linear strenghts. 

Let $\boldsymbol{\eta}\c 1_{\CLDC}\Rarr\inc;\B[-]$ be defined to be the family of maps $\boldsymbol{\eta}_{\bbX} = (\eta_{\bbX}, \eta_{\bbX})\c \bX\rarr \B[\bX]$. This is a 2-transformation as follows.
\begin{itemize}
    \item $\boldsymbol{\eta}$ satisfies 1-cell naturality: Let $F = (F_\times, F_+)\c \bX\rarr\bX'$ be a cartesian linear functor between CLDCs \bX\ and \bX'. 1-cell naturality in this case means 
    \[ \eta_{\bbX}; \B[\inc(F)]_\times = F_\times;\eta_{\bbX} \qquad \eta_{\bbX}; \B[\inc(F)]_+ = F_+; \eta_{\bbX}\]
    Now, the functors in the first equation are defined as follows
    \begin{align*}
    \eta_{\bbX}; &\B[\inc(F)]_\times\c \\
    & A\mapsto \langle F_\times(A), F_\times(\langle 1_A, 1_A\rangle); {m_\times}^{-1}_{A,A}, F_\times(t_A); {m_\bone^{-1}}, {p_+}_{A,A}; F_\times([1_A, 1_A]), {p_\bzero}; F_\times(b_A)\rangle \\
    & f \mapsto F_\times(f)\\ 
    \\
    F_\times; &\eta_{\bbX}\c\\
    &A\mapsto \langle F_\times(A), \langle 1_{F_\times(A)}, 1_{F_\times(A)}\rangle, t_{F_\times(A)}, [1_{F_\times(A)}, 1_{F_\times(A)}], b_{F_\times(A)}\rangle \\
    & f \mapsto F_\times(f)
    \end{align*}
    which are equal as the medial bimonoid structure on $F_\times(A)$ is unique, since \bX' is a CLDC. One shows that the other equality holds as both sides are functors mapping to the unique medial bimonoid structure on $F_+(A)$.
    
    \item $\boldsymbol{\eta}$ satisfies 2-cell naturality: Let $F = (F_\times, F_+), G=(G_\times, G_+)\c \bX\rarr\bX'$ be cartesian linear functor between CLDCs \bX\ and \bX', and $\alpha = (\alpha_\times, \alpha_+)\c F\Rarr G$ be a linear transformation between them. 2-cell naturality in this context means 
    \[ 1_{\eta_{\bbX}}; \B[\inc(\alpha)]_\times = \alpha_\times; 1_{\eta_{\bbX}} \qquad 1_{\eta_{\bbX}}; \B[\inc(\alpha)]_+ = \alpha_+; 1_{\eta_{\bbX}}\]

    Now, looking at the component morphisms for each object $A\in\bX$, these are medial bimonoid morphisms given by 
    \[ (1_{\eta_{\bbX}}; \B[\inc(\alpha)]_\times)_{A} = {\alpha_\times}_A = (\alpha_\times; 1_{\eta_{\bbX}})_A\]
    and therefore are equal. Similarly, for the second equality.
\end{itemize}

Secondly, let us define the counit $\boldsymbol{\epsilon}\c \B[-];\inc\Rarr 1_{\SMLDC}$. Let $(\bX, \ot, \top, \op, \bot)$ be a SMLDC, then define $\epsilon_{\bbX}\c\B[\bX]\rarr\bX$ to be the forgetful functor 
\[ \langle A, \Delta_A, e_A, \nabla_A, u_A\rangle \mapsto A \qquad f\mapsto f\]

It is immediate that $\boldsymbol{\epsilon}_{\bbX} =(\epsilon_{\bbX}, \epsilon_{\bbX})\c \B[\bX]\rarr \bX$ is a Frobenius medial linear functor when equipped with identity structure maps and linear strengths. Let $\boldsymbol{\epsilon}\c \B[-];\inc\Rarr 1_{\SMLDC}$ be the family of maps $\boldsymbol{\epsilon}_{\bbX} =(\epsilon_{\bbX}, \epsilon_{\bbX})\c \B[\bX]\rarr \bX$. It satisfies 1-cell and 2-cell naturality fairly trivially and therefore it is a 2-natural transformation.

To prove the desired 2-adjunction, it remains solely to show the triangle inequalities.
\[ \inc(\boldsymbol{\eta}); \boldsymbol{\epsilon}_{\inc} = 1_{\inc} \qquad \boldsymbol{\eta}_{\B[-]}; \B[\boldsymbol{\epsilon}] = 1_{\B[-]}\]
Consider a CLDC \bX, then \[ (\inc(\boldsymbol{\eta_{\bbX}}); \boldsymbol{\epsilon}_{\inc(\bbX)})_\times = (1_{\inc(\bbX)})_\times \qquad (\inc(\boldsymbol{\eta_{\bbX}}); \boldsymbol{\epsilon}_{\inc(\bbX)})_+ = (1_{\inc(\bbX)})_+\] as the as the left-hand side functors are defined by
\begin{align*}
(\inc(\boldsymbol{\eta_{\bbX}}); \boldsymbol{\epsilon}_{\inc(\bbX)})_\times = (\inc(\boldsymbol{\eta_{\bbX}}); \boldsymbol{\epsilon}_{\inc(\bbX)})_+ \c\qquad &A\mapsto \langle A, \langle 1_A, 1_A\rangle, t_A, [1_A, 1_A], b_A\rangle \mapsto A\\
&f\mapsto f\mapsto f
\end{align*}
Consider a SMLDC \bX, then 
\[ (\boldsymbol{\eta}_{\B[\bbX]}; \B[\boldsymbol{\epsilon}])_\ot = (1_{\B[\bbX]})_\ot \qquad (\boldsymbol{\eta}_{\B[\bbX]}; \B[\boldsymbol{\epsilon}])_\op = (1_{\B[\bbX]})_\op\] as the left-hand side functors are defined by
\begin{align*}
(\boldsymbol{\eta}_{\B[\bbX]}; \B[\boldsymbol{\epsilon}])_\ot=& (\boldsymbol{\eta}_{\B[\bbX]}; \B[\boldsymbol{\epsilon}])_\op \c\langle A, \Delta_A, e_A, \nabla_A, u_A\rangle \mapsto\,\langle \langle A\rangle, \langle 1_{\langle A\rangle}, 1_{\langle A\rangle} \rangle, t_{\langle A\rangle}, [1_{\langle A\rangle}, 1_{\langle A\rangle}], b_{\langle A\rangle}\rangle \\
&\mapsto\,  \langle \epsilon_{\bbX}(\langle A\rangle), \epsilon_{\bbX}(\langle 1_{\langle A\rangle}, 1_{\langle A\rangle} \rangle); 1_{A\ot A}, \epsilon_{\bbX}(t_{\langle A\rangle}); 1_\top, 1_{A\op A};\epsilon_{\bbX}([1_{\langle A\rangle}, 1_{\langle A\rangle}]), 1_\bot; \epsilon_{\bbX}(b_{\langle A\rangle}) \rangle \\
&= \langle A, \Delta_A, e_A, \nabla_A, u_A\rangle
\end{align*}
The final equality follows by Theorem \ref{thm:BX_CLDC}, since the unique maps induced by the universal properties of (co)products in $\B[\bX]$ are given by the structure maps of the medial bimonoids themselves, in other words
\[ \langle 1_{\langle A\rangle}, 1_{\langle A\rangle} \rangle, = \Delta_A \qquad  t_{\langle A\rangle} = e_A \qquad  [1_{\langle A\rangle}, 1_{\langle A\rangle}] = \nabla_A \qquad  b_{\langle A\rangle} = u_A\]
\end{proof}

\begin{remark}
With our main theorem proved, let us just reaffirm that a key reason it works is due to duoidal category theory. The theorem holds because CLDCs are canonically duoidal via their cartesian and co-cartesian structure, because MLDCs are coherently duoidal, because medial bimonoids are duoidal bimonoids for this structure and because medial linear functors/transformations are paired bilax duoidal functors/transformations which preserve duoidal bimonoids. In essence, buried within the linearly distributive Fox theorem is a duoidal Fox theorem. In fact, it would not be particularly difficult to state and prove this duoidal variant in full, the work in \cite[Chap 6]{Aguiar_Mahajan_2010} and Subsection \ref{sec:symmetry} having done most of the heavy lifting.
\end{remark}

Given that the component maps of the unit and counit of the adjunction are Frobenius, it is immediate that the 2-adjunction restricts to the Frobenius context:
\begin{corollary}
$\inc\vdash \B[-]\c\mathbf{FCLDC}\rarr \mathbf{SFMLDC}$ is a 2-adjunction.
\end{corollary}

Perhaps the most commonly used result of the traditional of the Fox theorem is that a SMC is cartesian if and only if it is isomorphic to its category of comonoids, as stated in Corollary \ref{cor:SMC_cartesian_iff}. We can of course state an analogous Corollary for SMLDCs.

\begin{corollary}
Given a SMLDC $\bX$, there exists a CLDC $\bY$ such that $\bX\cong \inc(\bY)$ if and only if there exists an isomorphism $\bX\cong \inc(\B[\bX])$.
\end{corollary}
\begin{proof}
$(\Rarr)$ Note first that since $\inc$ is fully faithful, the unit of the adjunction is an isomorphism and as such $\inc(\eta_{\bbY})$ is invertible for any $\bY\in\CLDC$. In particular, this means the adjunction is idempotent, in the sense of \cite[Sec 2.3]{Grandis_2020}. Now, consider $\bX\in\SMLDC$ such that there is a $\bY\in\CLDC$ with $\bX\cong \inc(\bY)$, then following is an isomorphism \[ \bX \cong \inc(\bY) \xrightarrow{\inc(\eta_{\bbY})} \inc(\B[\inc(\bY)]) \cong \inc(\B[\bX]) \]

$(\Leftarrow)$ Consider $\bX\in\SMLDC$ such that $\bX\cong \inc(\B[\bX])$ and let $\bY = \B[\bX]$, then $\bX\cong \inc(\bY)$.
\end{proof}

The above corollary can be unwrapped to give an explicit description in terms of medial bimonoid maps, as in the case of traditional cartesian categories.

\begin{corollary}
A SMLDC is cartesian if and only if there are natural transformations 
\[ \Delta_{A}\c A\rarr A\ot A \quad\quad e_{A}\c A\rarr \top \quad\quad \nabla_{A}\c A\op A\rarr A \quad\quad u_{A}:\bot\rarr A \]
such that, $\langle A, \Delta_{A}, e_{A}, \nabla_{A}, u_{A}\rangle$ determines a bicommutative medial bimonoid, and 
\begin{align*}
&\Delta_{A\ot B} = (\Delta_{A}\ot\Delta_{B}); \tau^\ot_{A, A, B, B} & e_{A\ot B} = (e_{A}\ot e_{B}); {u^{R}_{\ot}}_{\top}^{-1} \\
&\nabla_{A\ot B} = \mu_{A,B,A,B}; (\nabla_A \ot \nabla_B) & u_{A\ot B} = \nabla_\top; (u_A\ot u_B)\\
&\Delta_{A\op B} = (\Delta_A \op\Delta_B) ; \mu_{A,A,B,B} & e_{A\op B} = (e_A\op e_B);\Delta_\bot \\
&\nabla_{A\op B} = \tau^\op_{A, B, A, B};(\nabla_{A}\op\nabla_{B}) & u_{A\op B} = {u^R_\op}^{-1}_{\bot};(u_{A}\op u_{B}) 
\end{align*}

\[\Delta_{\top} = {u^{R}_{\ot}}_{\top} \qquad e_{\top} = 1_{\top} \qquad  \nabla_{\bot} = {u^{R}_{\op}}_{\bot} \qquad u_{\bot} = 1_{\bot} \]
\[ e_\bot = u_\top = m \]
\end{corollary}
This is essentially an amalgamations of Proposition \ref{prop:char_CLDC} and Corollary \ref{cor:char_CLDC_naive}, developed at the start to get a better understanding of the underlying structure of CLDCs. 

\begin{landscape}
\appendix
\section*{Appendix}
\addcontentsline{toc}{section}{Appendices}
\renewcommand{\thesubsection}{\Alph{subsection}}
\subsection{Large Commuting Diagrams}\label{app:commuting_diagrams}
\phantom{hello}
\vspace{5\baselineskip}
\begin{figure}[!h]
\resizebox{\linewidth}{!}{
\xymatrixrowsep{3.75pc}\xymatrixcolsep{4.75pc}\xymatrix@L=0.5pc{
\bot\ot\bot \ar[r]^-{1_\bot \ot {u_\op}_\bot^{-1}}\ar[rd]_-{{u^R_\op}^{-1}_{\bot\ot\bot}}^-{\quad(\mathrm{\bf LDC. 1})} & \bot\ot(\bot\op\bot)\ar[r]^-{1_\bot \ot({u^L_\ot}_\bot \op {u^R_\ot}_\bot)}\ar[d]^-{\delta^L_{\bot,\bot,\bot}}\ar@{}[rd]|{(\nat)} & \bot\ot((\top\ot\bot)\op(\bot\ot\top))\ar@{}[rdd]|{\eqref{cc:sym_medial_lindist}} \ar[r]^-{1_{\bot}\ot\,\mu_{\top,\bot,\bot,\top}} \ar[d]^-{\delta^L_{\bot,\top\ot\bot,\bot\ot\top}} & \bot\ot((\top\op\bot)\ot(\bot\op\top))\ar[rr]^-{1_\bot \ot({u^R_\op}_\top \ot {u^L_\op}_\top)} \ar[d]_-{{\alpha_\ot}_{\bot, \top\op\bot, \bot\op\top}^{-1}}\ar@{}[rd]|{(\nat)} && \bot\ot(\top \ot \top)\ar[dd]^-{1_\bot \ot {u_\ot}_\top^{-1}}\ar[ld]_-{{\alpha_\ot}^{-1}_{\bot,\top,\top}}\\
& (\bot\ot\bot)\op \bot\ar[r]_-{(1_\bot\ot{u_\ot^L}_\bot)\op {u_\ot^R}_{\bot}}\ar[rdd]_-{1_{\bot\ot\bot}\op {u_\ot^R}_{\bot}}& (\bot\ot(\top\ot\bot))\op(\bot\ot\top)\ar[d]^-{{\alpha_\ot}^{-1}_{\bot,\top,\bot}\op 1_{\bot\ot\top}} \ar@{}[ld]|{(\mathrm{\bf MC.2})} & (\bot\ot(\top\op\bot))\ot (\bot\op \top)\ar[d]_{\delta^L_{\bot,\top,\bot}\ot1_{\bot\op\top}}\ar[r]^-{(1_\bot\ot {u^R_\op}_\top)\ot {u^L_\op}_\top} & (\bot\ot\top)\ot\top\ar[rd]_-{{u^R_\ot}^{-1}_\bot\ot 1_\top}\ar@{}[r]|{\quad\quad(\mathrm{\bf MC.2})} &\\
\bot\ar[uu]^{\Delta_{\bot}}\ar[dd]_-{{u^R_\ot}_\bot}\ar[r]_-{{u_\op}^{-1}_\bot}\ar@{}[ru]|{(\nat)} \ar@{}[rd]|{(\nat)} &\bot\op\bot\ar[u]^-{\Delta_\bot\op 1_\bot}\ar[d]_-{1_\bot\op {u^R_\ot}_\bot} &((\bot\ot\top)\ot \bot)\op(\bot\ot\top)\ar[d]^-{({u^R_\ot}^{-1}_\bot \ot 1_\bot)\op 1_{\bot\ot\top}} \ar[r]^-{\mu_{\bot\ot\top,\bot,\bot,\top}} \ar@{}[rd]|{(\nat)}&((\bot\ot\top)\op\bot)\ot(\bot\op\top)\ar[d]_-{({u^R_\ot}^{-1}_\bot \op 1_\bot)\ot 1_{\bot\op\top}}\ar[ru]_-{{u_\op^R}_{\bot\ot\top}\ot {u_\op^L}_\top}^-{(\mathrm{\bf LDC.1})\quad\quad}\ar@{}[rr]|{(\nat)} && \bot\ot\top\\
& \bot\op(\bot\ot\top)\ar[r]_-{\Delta_\bot\op 1_{\bot\ot\top}}\ar[d]^-{{\sigma_\op}_{\bot,\bot\ot\top}}\ar@{}[ru]|{(\exch)}& (\bot\ot\bot)\op(\bot\ot\top)\ar[d]^-{{\sigma_\op}_{\bot\ot\bot,\bot\ot\top}} \ar[r]^-{\mu_{\bot,\bot,\bot,\top}} \ar@{}[rd]|{\eqref{cc:sym_medial_braiding}} &(\bot\op\bot)\ot(\bot\op\top)\ar[d]_{{\sigma_\op}_{\bot,\bot}\ot {\sigma_\op}_{\bot,\top}}\ar[rru]_-{{u_\op}_\bot\ot{u_\op^L}_\top}\ar@{}[rrd]|{(\mathrm{\bf BMC.2})} \\
\bot\ot\top \ar[r]_-{{u^R_\op}^{-1}_{\bot\ot\top}}\ar[ru]^-{{u^L_\op}^{-1}_{\bot\ot\top}}_-{\quad(\mathrm{\bf BMC.2})}& (\bot\ot\top)\op\bot \ar[r]_{1_{\bot\ot\top}\op\Delta_{\bot}} \ar@{}[ru]|{(\nat)}& (\bot\ot\top)\op(\bot\ot\bot)\ar[r]_{\mu_{\bot,\top,\bot,\bot}} & (\bot\op\bot)\ot(\top\op\bot)\ar[rr]_-{{u_\op}_\bot \ot {u^R_\op}_\top} && \bot\ot\top \ar@{=}[uu]
}}
\caption{Proof that $(\bot, \Delta_\bot, m)$ is a cocommutative $\ot$-comonoid in Proposition \ref{prop:alternative_SMLDC}}
\label{fig:bot_cocomm_ot_comonoid}
\end{figure}
\end{landscape}

\begin{landscape}
\begin{figure}
\resizebox{\linewidth}{!}{\xymatrixrowsep{4.75pc}\xymatrixcolsep{3.75pc}\xymatrix@L=0.7pc{
(A\ot B)\ot ((C\ot D)\op (E\ot F))\ar[rrrrrr]^-{\delta^L_{A\ot B, C\ot D, E\ot F}}\ar[d]_-{1_{A\ot B}\ot \mu_{C,D,E,F}}\ar[rd]^-{{\alpha_\ot}_{A,B,(C\ot D)\op (E\ot F)}} & & &\ar@{}[d]|{(\mathrm{\bf LDC.2})} & & & ((A\ot B)\ot (C\ot D))\op (E\ot F) \ar[d]^-{{\alpha_\ot}_{A,B,C\ot D}\op 1_{E\ot F}}\\
(A\ot B)\ot ((C\op E)\ot (D\op F))\ar@{}[r]|{(\nat)}\ar[d]_-{{\alpha_\ot}_{A, B, (C\op E)\ot (D\op F)}} & A\ot (B\ot ((C\ot D)\op (E\ot F)))\ar[rrrr]^-{1_A\ot \delta^L_{B, C\ot D, E\ot F}}\ar[d]^-{1_A\ot {\sigma_\ot}_{B, (C\ot D)\op (E\ot F)}}\ar[ld]^-{1_A\ot (1_B\ot \mu_{C,D,E,F})} & & \ar@{}[d]|{(\mathrm{\bf SLDC})}& & A\ot ((B\ot (C\ot D))\op (E\ot F))\ar@{}[rd]|{(\nat)}\ar[r]^-{\delta^L_{A, B\ot (C\ot D), E\ot F}}\ar[d]_-{1_A\ot ({\sigma_\ot}_{B,C\ot D}\op 1_{E\ot F})} & (A\ot (B\ot (C\ot D)))\op (E\ot F)\ar[d]^-{(1_A\ot \sigma_{B, C\ot D})\op 1_{E\ot F}}\\
A\ot (B\ot ((C\op E)\ot (D\op F)))\ar@{}[r]|{(\nat)}\ar[d]_-{1_A\ot {\sigma_\ot}_{B, (C\op E)\ot (D\op F)}} & A\ot (((C\ot D)\op (E\ot F))\ot B)\ar@{}[dd]|{(\mathrm{\bf BDUO})}\ar[ld]^-{1_A\ot (\mu_{C, D, E, F}\ot 1_B)}\ar[r]^-{1_A\ot ({\sigma_\ot}_{C\ot D, E\ot F}\ot 1_B)} & A\ot (((E\ot F)\op (C\ot D))\ot B)
\ar@{}[rrdd]|{(\mathrm{\bf MLDC.1})}\ar[d]^-{1_A\ot (\mu_{E,F,C,D}\ot 1_B)}\ar[rr]^-{1_A\ot \delta^R_{E\ot F, C\ot D, B}}& & A\ot ((E\ot F)\op ((C\ot D)\ot B))\ar@{}[rdd]|{(\nat)}\ar[d]_-{1_A\ot (1_{E\ot F}\op {\alpha_\ot}_{C,D,B})} & A\ot (((C\ot D)\ot B) \op (E\ot F))\ar[l]^-{1_A\ot {\sigma_\op}_{(C\ot D)\ot B, E\ot F}}\ar[r]^-{\delta^L_{A, (C\ot D)\ot B, E\ot F}}\ar[d]_-{1_A\ot ({\alpha_\ot}_{C,D,B}\op 1_{E\ot F})}\ar@{}[rd]|{(\nat)} & (A\ot ((C\ot D)\ot B)) \op (E\ot F) \ar[d]^-{(1_A\ot {\alpha_\ot}_{C,D,B})\op 1_{E\ot F}}\\
A\ot (((C\op E)\ot (D\op F))\ot B)\ar[dd]_-{1_A\ot {\alpha_\ot}_{C\op E, D\op F, B}}\ar[rd]^-{1_A\ot ((1_{C\op E}\ot {\sigma_\op}_{D,F})\ot 1_B)} & & A\ot (((E\op C)\ot (F\op D))\ot B)\ar[d]^-{1_A\ot {\alpha_\ot}_{E\op C, F\op D, B}} & & A\ot ((E\ot F)\op (C\ot (D\ot B)))\ar[d]^-{1_A\ot (1_{E\ot F}\op (1_C\ot {\sigma_\ot}_{D,B}))}\ar[ld]_-{1_A\ot \mu_{E,F,C,D\ot B}} & A\ot ((C\ot(D\ot B))\op (E\ot F))\ar[r]^-{\delta^L_{A, C\ot(D\ot B), E\ot F}}\ar[ddl]^-{1_A\ot ((1_C\ot {\sigma_\ot}_{D,B})\op 1_{E\ot F})} \ar@{}[rdd]|{(\nat)}& (A\ot (C\ot (D\ot B)))\op (E\ot F)\ar[dd]^-{(1_A\ot (1_C\ot {\sigma_\ot}_{D,B}))\op 1_{E\ot F}}\\
\ar@{}[r]|{(\nat)}& A\ot (((C\op E)\ot (F\op D))\ot B)\ar@{}[r]|{(\nat)} \ar[d]_-{1_A\ot {\alpha_\ot}_{E\op C, D\op F, B}}\ar[ru]^-{1_A\ot( ({\sigma_\op}_{C,E} \ot 1_{F\op D})\ot 1_B)}& A\ot ((E\op C)\ot ((F\op D)\ot B)) \ar@{}[d]|{(\exch)}\ar[r]^-{1_A\ot (1_{E\op C}\ot \delta^R_{F,D,B})} & A\ot ((E\op C)\ot (F\op (D\ot B)))\ar[d]_(.65){1_A\ot (1_{E\op C}\ot (1_F \op {\sigma_\ot}_{D,B}))}\ar@{}[r]|{(\nat)} & A\ot ((E\ot F)\op (C\ot (B\ot D))) \ar[d]_-{1_A\ot {\sigma_\op}_{E\ot F, C\ot(B\ot D)}}\ar[ld]_-{1_A\ot \mu_{E, F, C, B\ot D}} & \\
A\ot ((C\op E)\ot ((D\op F)\ot B))\ar@{}[rrdd]|{(\mathrm{\bf SLDC})}\ar[dd]_-{1_A\ot (1_{C\op E}\ot {\sigma_\ot}_{D\op F, B})} \ar[r]^-{1_A\ot (1_{C\op E}\ot ({\sigma_\op}_{D,F}\ot 1_B))}& A\ot((C\op E)\ot ((F\op D)\ot B))\ar[r]^-{1_A\ot (1_{C\op E}\ot \delta^R_{F, D, B})}\ar[ur]^-{1_A\ot ({\sigma_\op}_{C,E}\ot 1_{(F\op D)\ot B})} & A\ot ((C\op E)\ot (F\op (D\ot B)))\ar@{}[r]|{(\exch)}\ar[d]_-{1_A\ot (1_{C\op E \ot (1_F\ot {\sigma_\ot}_{D,B}))}}\ar[ru]^-{1_A\ot ({\sigma_\op}_{C,E}\ot 1_{F\op(D\ot B)})} & A\ot ((E\op C)\ot(F\op(B\ot D)))\ar@{}[d]|{(\mathrm{\bf BDUO})}\ar[ld]^-{1_A\ot ({\sigma_\op}_{E,C}\ot 1_{F\op(B\ot D)})} & A\ot ((C\ot (B\ot D))\op (E\ot F))\ar@{}[rrddd]|{(\mathrm{\bf MLDC.1})}\ar[rr]^-{\delta^L_{A, C\ot (B\ot D), E\ot F}}\ar[lldd]^-{1_A\ot \mu_{C,B\ot D, E, F}} && (A\ot (C\ot (B\ot D)))\op (E\ot F)\ar[d]^-{{\alpha_\ot^{-1}}_{A,C,B\ot D}\op 1_{E\ot F}}\\
& & A\ot ((C\op E)\ot (F\op (B\ot D)))\ar[d]_-{1_A\ot (1_{C\op E}\ot {\sigma_\op}_{F, B\ot D})} & &&& ((A\ot C)\ot (B\ot D))\op (E\ot F)\ar[dd]^-{\mu_{A,C,B,D}\op 1_{E\ot F}} \\
A\ot ((C\op E)\ot (B\ot (D\op F)))\ar@{}[rrd]|{(\nat)}\ar[d]_-{{\alpha_\ot^{-1}}_{A,C\op E, B\ot (D\op F)}}\ar[rr]^-{1_A\ot (1_{C\op E}\ot \delta^L_{B,D,F})} && A\ot ((C\op E)\ot ((B\ot D)\op F)))\ar[d]_-{{\alpha_\ot^{-1}}_{A, C\op E, (B\ot D)\op F}} \\
(A\ot (C\op E)) \ot (B\ot (D\op F))\ar[rr]_-{1_{A\ot(C\op E)}\ot \delta^L_{B,D,F}} & & (A\ot (C\op E))\ot ((B\ot D)\op F)\ar[rrrr]_-{\delta^L_{A,C,E}\ot 1_{(B\ot D)\op F}} &&&&((A\ot C)\op E) \ot ((B\ot D)\op F)
}}
\caption{Proof of the third equality in Proposition \ref{prop:interaction_canonical_flip_medial_maps_linear_dist}}
\label{fig:interaction_canonical_flip_medial_maps_linear_dist}
\end{figure}

\begin{figure}
\begin{equation*}\begin{gathered}
\resizebox{\linewidth}{!}{\xymatrixrowsep{5pc}\xymatrixcolsep{6pc}\xymatrix@L=0.5pc{
A\ot (B\op C)\ar@{}[rrdddd]|{\qquad(\mathrm{\bf MLDC.4})}\ar@{=}[ddddddd]\ar[rr]^-{\delta^L_{A,B,C}}\ar[rrd]^-{1_A\ot (1_B\op {u^L_\ot}_C)}\ar[rddd]^-{1_A\ot {u^L_\ot}_{B\op C}} && (A\ot B)\op C\ar[rr]^-{1_{A\ot B}\op {u^L_\ot}_C}\ar@{}[d]|{(\nat)} && (A\ot B)\op (\top \ot C)\ar[ddddddd]^-{\mu_{A,B,\top,C}}\ar@/^1.5pc/[lddd]^-{({u^R_\ot}_A\ot 1_B)\op 1_{\top\ot C}} \\
&& A\ot (B\op(\top \ot C))\ar[r]^-{\delta^L_{A,B,\top\ot C}}\ar[d]_-{1_A\ot ({u^L_\ot}_B\op 1_{\top\ot C})}\ar@{}[rd]|{(\nat)} & (A\ot B)\op (\top\ot C)\ar[d]_-{(1_A\ot {u^L_\ot}_B)\op 1_{\top\ot C}}\ar@{=}[ur] \ar@{}[rd]|{(\mathrm{\bf MC.2})\qquad}& \\
&& A\ot ((\top\ot B)\op (\top\op C))\ar[d]_-{1_A\ot \mu_{\top,B,\top,C}} \ar[r]^-{\delta^L_{A,\top\ot B, \top\ot C}} \ar@{}[rdd]|{(\mathrm{\bf MLDC.1})}& (A\ot (\top\ot B))\op (\top\ot C)\ar[d]_-{{\alpha_\ot}^{-1}_{A,\top,B}\op 1_{\top\ot C}} &\\
\ar@{}[rdd]|{(\mathrm{\bf MC.2})\qquad}& A\ot (\top\ot(B\op C))\ar[d]_-{{\alpha_\ot}^{-1}_{A, \top, B\op C}}\ar@{}[rd]|{(\nat)}& A\ot ((\top\op\top)\ot (B\op C))\ar[d]_-{{\alpha_\ot}^{-1}_{A,\top\op\top, B\op C}} \ar[l]_-{1_A\ot (\nabla_\top\ot 1_{B\op C}}& ((A\ot \top)\ot B)\op (\top\ot C)\ar[d]_-{\mu_{A\ot\top, B, \top, C}}\ar@{}[rd]|{(\nat)} & \\
& (A\ot\top)\ot (B\op C) & (A\ot (\top\op\top))\ot (B\op C)\ar@{}[d]|{(\nat)}\ar[r]^-{\delta^L_{A,\top,\top}\ot 1_{B\op C}}\ar[dl]^*+<8pt>{^{(1_A\ot (1_\top\op m^{-1}))\ot 1_{B\op C}}}_-{(\mathrm{\bf MLDC.2})\qquad}\ar[l]_-{(1_A\ot \nabla_\top)\ot 1_{B\op C}} & ((A\ot \top)\op\top)\ot (B\op C)\ar[dl]_*+<8pt>{_{(1_{A\ot\top}\op m^{-1})\ot 1_{B\op C}}}\ar@{}[ddd]|{(\exch)} & \\
& (A\ot (\top\op\bot)) \ot (B\op C)\ar[r]_-{\delta^L_{A,\top,\bot}\ot 1_{B\op C}}\ar[u]^-{(1_A\ot {u^R_\op}_\top)\ot 1_{B\op C}}\ar[d]_-{(1_A\ot {u^R_\op}_\top)\ot 1_{B\op C}} & ((A\ot \top)\op \bot)\ot (B\op C)\ar[ld]^-{{u^R_\op}_{A\ot \top}\op 1_{B\op C}}_-{(\mathrm{\bf LDC.1})\qquad} &\\
& (A\ot \top)\ot (B\op C)\ar@{}[rd]|{(\nat)} & & & \\
A\ot (B\op C)\ar[ru]_-{{u^R_\ot}_A\ot 1_{B\op C}} && (A\op \bot)\ot (B\op C)\ar[ll]^-{{u^R_\op}_A\ot 1_{B\op C}}\ar[uu]_-{({u^R_\ot}_A\op 1_\bot)\ot 1_{B\op C}} && (A\op \top)\ot (B\op C)\ar[ll]^-{(1_A\op m^{-1})\ot 1_{B\op C}}\ar[luuu]_-{({u^R_\ot}_A\op 1_\top)\ot 1_{B\op C}}
}}
\end{gathered}\end{equation*}
\caption{Proof that $\partial^R$ is the right inverse of $\delta^L$ in Proposition \ref{prop:isomix_MLDC}}
\label{fig:partial_R_inverse_delta_L}
\end{figure}

\begin{figure}
\resizebox{\linewidth}{!}{
\xymatrixrowsep{5pc}\xymatrixcolsep{3pc}\xymatrix@L=0.7pc{
(A\ot B)\op(C\ot D)\ar[r]^-{\mu_{A,B,C,D}}\ar[d]_-{(1_A\ot {u^L_\ot}_B^{-1})\op (1_C\ot {u^L_\ot}^{-1}_D)}\ar@{}[rdd]|{(\nat)} & (A\op C)\ot (B\op D)\ar[rr]^-{1_{A\op C}\ot {u^L_\op}^{-1}_{B\op C}}\ar[d]_-{1_{A\op C}\ot ({u^L_\op}^{-1}_B\op {u^L_\op}_D^{-1})}\ar@{}[rrd]|{\eqref{cc:canonical_flip_unitor}} && (A\op C) \ot (\bot\op (B\op D))\ar[dd]^-{1_{A\op C}\ot (m\op 1_{B\op D})}\\
(A\ot (\bot\op B)) \op (C\ot(\bot\op D))\ar[d]_-{(1_A\ot(m\op1_B))\op (1_C\ot (m\op 1_D))} & (A\op C)\ot((\bot\op B)\op (\bot\op D))\ar[r]^-{1_{A\op C}\ot \tau^\op_{\bot,B,\bot,D}}\ar[d]_-{1_{A\op C}\ot ((m\op 1_B)\op (m\op 1_D))}\ar@{}[rd]|{(\nat)}& (A\op C)\op ((\bot\op\bot)\op(B\op D))\ar[ru]^-{1_{A\op C}\ot ({u_\op}_\bot\op 1_{B\op D})}\ar[d]^-{1_{A\op C}\ot ((m\op m)\op 1_{B\op D})}\ar@{}[r]|{\qquad\qquad\qquad(\mathrm{\bf MLDC.2})}  & \\
(A\ot (\top\op B))\op (C\ot (\top\op B))\ar[r]^-{\mu_{A,\top\op B, C, \top\op D}}\ar[d]_-{\delta^L_{A,\top,B}\op \delta^L_{C,\top,D}}\ar@{}[rrd]|{\eqref{cc:interaction_canonical_flip_medial_maps_linear_dist}} & (A\op C)\ot ((\top\op B)\op (\top\op D))\ar[r]^-{1_{A\op C}\ot \tau^\op_{\top,B,\top,D}} & (A\op C)\ot ((\top\op\top)\op (B\op D))\ar[r]^-{1_{A\op C}\ot (\nabla_\top\op 1_{B\op D}))}\ar[d]^-{\delta^L_{A\op C, \top\op\top, B\op D}}\ar@{}[rd]|{(\nat)} & (A\op C)\ot (\top\op (B\op D))\ar[d]^-{\delta^L_{A\op C, \top, B\op D}}\\
((A\ot \top)\op B)\op ((C\ot \top)\op D) \ar@{}[rd]|{(\nat)}\ar[d]_-{({u^R_\ot}_A^{-1}\op 1_B)\op ({u^R_\ot}^{-1}_C\op 1_D)}\ar[r]^-{\tau^\op_{A\ot \top, B, C\ot \top, D}}& ((A\ot \top)\op (C\ot\top))\op (B\op D)\ar[rrd]_-{({u^R_\ot}_A^{-1}\op {u^R_\ot}_C^{-1})\op 1_{B\op D}}\ar[r]^-{\mu_{A,\top,C,\top}\op 1_{B\op D}} & ((A\op C)\ot (\top\op\top))\op (B\op D)\ar[r]^-{(1_{A\op C}\ot \nabla_\top)\ot 1_{B\op D}}\ar@{}[rd]|{(\mathrm{\bf MLDC.4})} & ((A\op C)\ot \top)\op (B\op D)\ar[d]^-{{u^R_\ot}^{-1}_{A\op C}\op 1_{B\op D}} \\
(A\op B)\op (C\op D)\ar[rrr]_-{\tau^\op_{A,B,C,D}} &&& (A\op C)\op (B\op D)
}}
\caption{Proof of the first diagram in Proposition \ref{prop:interchange_canonicalflip_mix}}
\label{fig:interaction_canonical_flip_medial_maps_mix}
\end{figure}

\begin{figure}
\resizebox{\linewidth}{!}{\xymatrixrowsep{3.75pc}\xymatrixcolsep{4.75pc}\xymatrix@L=0.7pc{
F_\ot(A\op B) \op F_\ot(C\op D)\ar@{}[rddd]|{(\mathrm{\bf SMLF.6})}\ar[dd]_-{{m_\op}_{A\op B, C\op D}}\ar[r]^-{{\nu^R_\ot}_{A,B}\op 1_{F_\ot(C\op D)}} & (F_\op(A)\op F_\ot(B)) \op F_\ot(C\op D)\ar[d]^-{{\alpha_\op}^{-1}_{F_\op(A), F_\ot(B), F_\ot(C\op D)}} \ar[rr]^-{1_{F_\op(A)\op F_\ot(B)}\op {\nu^R_\op}_{C,D}}  & \ar@{}[d]|{(\nat)}& (F_\op(A)\op F_\ot(B))\op (F_\op(C)\op F_\ot(D))\ar[d]^-{{\alpha_\op}^{-1}_{F_\op(A), F_\ot(B), F_\op(C), F_\ot(D)}}\\
&F_\op(A) \op (F_\ot(B)\op F_\ot(C\op D))\ar[dd]^-{1_{F_\op(A)}\op {m_\op}_{B, C\op D}}\ar[rd]^-{1_{F_\op(A)}\op {\sigma_\op}_{F_\ot(B), F_\ot(C\op D)}}\ar[rr]^{1_{F_\op(A)}\op (1_{F_\ot(B)}\op {\nu^R_\ot}_{C,D})} && F_\op(A) \op (F_\ot(B)\op (F_\op(C)\op F_\ot(D))\ar[d]^-{1_{F_\op(A)}\op {\sigma_\op}_{F_\ot(B), F_\op(C)\op F_\ot(D)}}\\
F_\ot((A\op B) \op (C\op D)) \ar[d]_-{F_\ot({\alpha_\op}^{-1}_{A,B,C\op D})} &\ar@{}[r]|{(\mathrm{\bf SMF})} & F_\op(A) \op (F_\ot(C\op D)\op F_\ot(B)))\ar@{}[u]|{(\nat)}\ar@{}[ddd]|{(\mathrm{\bf SMLF.6})}\ar[ldd]^-{1_{F_\op(A)}\op {m_\op}_{C\op D, B}}\ar[r]^-{1_{F_\op(A)}\op ({\nu^R_\ot}_{C,D}\op 1_{F_\op(B)})}& F_\op(A)\op ((F_\op(C)\op F_\ot(D))\op F_\ot(B))\ar[d]^-{1_{F_\op(A)}\op {\alpha_\op}^{-1}_{F_\op(C), F_\ot(D), F_\ot(B)}}\\
F_\ot(A\op (B\op (C\op D))\ar@{}[rdd]|{(\nat)}\ar[r]^-{{\nu^R_\ot}_{A, B\op(C\op D)}}\ar[d]_-{F_\ot(1_A \op {\sigma_\op}_{B, C\op D})} & F_\op(A) \op F_\ot(B\op (C\op D))\ar[d]_-{1_{F_\op(A)}\op F_\ot({\sigma_\op}_{B, C\op D})} & & F_\op(A)\op (F_\op(C)\op (F_\ot(D)\ot F_\ot(B)))\ar[d]^-{1_{F_\op(A)}\op (1_{F_\op(C)}\op {\sigma_\op}_{F_\ot(D), F_\ot(B)})}\ar[ldd]_*+<0.5pc>{_{1_{F_\op(A)}\op (1_{F_\op(C)}\op {m_\op}_{D,B})}}\\
F_\ot(A\op ((C\op D)\op B))\ar[d]_-{F_\ot(1_A\op {\alpha_\op}^{-1}_{C, D, B})} & F_\op (A) \op F_\ot((C\op D)\op B)\ar[d]_-{1_{F_\op(A)} \op F_\ot({\alpha_\op}^{-1}_{C,D,B})} && F_\op(A)\op (F_\op(C) \op(F_\ot(B)\op F_\ot(D)))\ar[d]^-{{\alpha_\op}_{F_\op(A), F_\op(C), F_\ot(B), F_\ot(D)}}\ar[ldd]^*+<0.5pc>{^{1_{F_\op(A)}\op (1_{F_\op(C)}\op {m_\op}_{B,D})}}\ar@{}[ld]|{(\mathrm{\bf SMF})}\\
F_\op(A\op (C\op (D\op B)))\ar[d]_-{F_\ot(1_A\op (1_C\op {\sigma_\op}_{D,B}))} & F_\op (A) \op F_\ot(C\op (D\op B))\ar[d]_-{1_{F_\op(A)}\op F_\ot(1_C\op {\sigma_\op}_{D,B})}^-{\qquad\qquad(\nat)}\ar[r]^-{1_{F_\op(A)}\op {\nu^R_\ot}_{C, D\op B}} & F_\op(A)\op (F_\op(C)\op F_\ot(D\op B))\ar[d]_-{1_{F_\op(A)}\op (1_{F_\op(C)}\op F_\ot({\sigma_\op}_{D,B}))} & (F_\op(A)\op F_\op(C))\op (F_\ot(B)\op F_\ot(D))\ar[d]^-{1_{F_\op(A)\op F_\op(C)}\op {m_\op}_{B,D}}\\
F_\op (A\op (C\op (B\op D)))\ar[d]_-{F_\ot({\alpha_\op}_{A,C,B\op D})} \ar[r]_-{{\nu^R_\ot}_{A,C\op(B\op D)}}& F_\op (A) \op F_\ot(C\op (D\op B))\ar[r]_-{1_{F_\op(A)}\op {\nu^R_\ot}_{C, B\op D}} & F_\op(A)\op (F_\op(C)\op F_\ot(B\op D))\ar[r]_-{{\alpha_\op}_{F_\op(A), F_\op(C), F_\ot(B\op D)}}\ar@{}[ru]|{(\nat)}& (F_\op(A)\op F_\op(C)) \op F_\ot(B\op D)\ar[d]^-{{n_\op}^{-1}_{A,C}\op 1_{F_\ot(B\op D)}}\\
F_\op((A\op C)\op (B\op D))\ar[rrr]_-{{\nu^R_\ot}_{A\op C, B\op D}}\ar@{}[rrru]|{(\mathrm{\bf LF.2})}& & & F_\op(A\op C)\op F_\ot(B\op D)
}}
\caption{Proof of the third equality in Proposition \ref{prop:interaction_canonical_flip_medial_linear_functor}}
\label{fig:proof_interaction_canonical_flip_medial_linear_functor}
\end{figure}
\end{landscape}

\begin{references*}

\bibitem{Aguiar_Mahajan_2010}
M.~Aguiar and S.~Mahajan.
\newblock \emph{Monoidal functors, species and {H}opf algebras}, volume~29 of \emph{CRM Monograph Series}.
\newblock American Mathematical Society, Providence, RI, 2010.

\bibitem{Balteanu_Fiedorowicz_Schwaenzl_Vogt_2003}
C.~Balteanu, Z.~Fiedorowicz, R.~Schw\"anzl, and R.~Vogt.
\newblock Iterated monoidal categories.
\newblock \emph{Adv. Math.}, 176\penalty0 (2):\penalty0 277--349, 2003.

\bibitem{Barr_1979}
M.~Barr.
\newblock \emph{{$\ast$}-autonomous categories}, volume 752 of \emph{Lecture Notes in Mathematics}.
\newblock Springer, Berlin, 1979.
\newblock With an appendix by Po Hsiang Chu.

\bibitem{Batanin_Markl_2012}
M.~Batanin and M.~Markl.
\newblock Centers and homotopy centers in enriched monoidal categories.
\newblock \emph{Adv. Math.}, 230\penalty0 (4-6):\penalty0 1811--1858, 2012.

\bibitem{Blute_Panangaden_Slanov_2012}
R.F. Blute, P.~Panangaden, and S.~Slavnov.
\newblock Deep inference and probabilistic coherence spaces.
\newblock \emph{Appl. Categ. Structures}, 20\penalty0 (3):\penalty0 209--228, 2012.

\bibitem{Bonchi_DiGiorgio_Haydon_Sobocinski_2024}
F.~Bonchi, A.~Di~Giorgio, N.~Haydon, and P.~Sobocinski.
\newblock Diagrammatic algebra of first order logic.
\newblock In \emph{Proceedings of the 39th Annual ACM/IEEE Symposium on Logic in Computer Science}, LICS '24, New York, NY, USA, 2024. Association for Computing Machinery.

\bibitem{Bonchi_DiGiorgio_Trotta_2024}
F.~Bonchi, A.~Di~Giorgio, and D.~Trotta.
\newblock When {L}awvere meets {P}eirce: an equational presentation of {B}oolean hyperdoctrines.
\newblock In \emph{49th {I}nternational {S}ymposium on {M}athematical {F}oundations of {C}omputer {S}cience}, volume 306 of \emph{LIPIcs. Leibniz Int. Proc. Inform.}, pages Art. No. 30, 19. Schloss Dagstuhl. Leibniz-Zent. Inform., Wadern, 2024.

\bibitem{Booker_Street_2013}
T.~Booker and R.~Street.
\newblock Tannaka duality and convolution for duoidal categories.
\newblock \emph{Theory Appl. Categ.}, 28\penalty0 (6):\penalty0 166--205, 2013.

\bibitem{Brunnler_2003_1}
K.~Br\"unnler.
\newblock Locality for classical logic.
\newblock \emph{Notre Dame J. Formal Logic}, 47\penalty0 (1):\penalty0 557--580, 2006.

\bibitem{Brunnler_Tiu_2001}
K.~Br\"unnler and A.~F. Tiu.
\newblock A local system for classical logic.
\newblock In \emph{Logic for programming, artificial intelligence, and reasoning ({H}avana, 2001)}, volume 2250 of \emph{Lecture Notes in Comput. Sci.}, pages 347--361. Springer, Berlin, 2001.

\bibitem{Bruscoli_Strassburger_2017}
P.~Bruscoli and L.~Stra{\ss}burger.
\newblock On the length of medial-switch-mix derivations.
\newblock In Juliette Kennedy and Ruy~J.G.B. de~Queiroz, editors, \emph{Logic, Language, Information, and Computation}, pages 68--79, Berlin, Heidelberg, 2017. Springer Berlin Heidelberg.

\bibitem{Carboni_Walters_1987}
A.~Carboni and R.~F.~C. Walters.
\newblock Cartesian bicategories. {I}.
\newblock \emph{J. Pure Appl. Algebra}, 49\penalty0 (1-2):\penalty0 11--32, 1987.

\bibitem{Cockett_Seely_1997}
J.~R.~B. Cockett and R.~A.~G. Seely.
\newblock Proof theory for full intuitionistic linear logic, bilinear logic, and {MIX} categories.
\newblock \emph{Theory Appl. Categ.}, 3\penalty0 (5):\penalty0 85--131, 1997.

\bibitem{Cockett_Seely_1997_LDC}
J.~R.~B. Cockett and R.~A.~G. Seely.
\newblock Weakly distributive categories.
\newblock \emph{J. Pure Appl. Algebra}, 114\penalty0 (2):\penalty0 133--173, 1997.

\bibitem{Cockett_Seely_1999}
J.~R.~B. Cockett and R.~A.~G. Seely.
\newblock Linearly distributive functors.
\newblock \emph{J. Pure Appl. Algebra}, 143\penalty0 (1-3):\penalty0 155--203, 1999.
\newblock Special volume on the occasion of the 60th birthday of Professor Michael Barr (Montreal, QC, 1997).

\bibitem{Cockett_Koslowski_Seely_2000}
J.~R.~B. Cockett, J.~Koslowski, and R.~A.~G. Seely.
\newblock Introduction to linear bicategories.
\newblock \emph{Math. Structures Comput. Sci.}, 10\penalty0 (2):\penalty0 165--203, 2000.
\newblock The Lambek Festschrift: mathematical structures in computer science (Montreal, QC, 1997).

\bibitem{Cockett_Comfort_Srinivasan_2021}
J.~R.~B. Cockett, C.~Comfort, and P.~V. Srinivasan.
\newblock Dagger linear logic for categorical quantum mechanics.
\newblock \emph{Log. Methods Comput. Sci.}, 17\penalty0 (4):\penalty0 Paper No. 8, 73, 2021.

\bibitem{Das_Strassburger_2016}
A.~Das and L.~Stra{\ss}burger.
\newblock On linear rewriting systems for boolean logic and some applications to proof theory.
\newblock \emph{Logical Methods in Computer Science}, Dec 2016.

\bibitem{Forcey_Siehler_Sowers_2007}
S.~Forcey, J.~Siehler, and E.~S. Sowers.
\newblock Operads in iterated monoidal categories.
\newblock \emph{J. Homotopy Relat. Struct.}, 2\penalty0 (1):\penalty0 1--43, 2007.

\bibitem{Fox_1976}
T.~Fox.
\newblock Coalgebras and {C}artesian categories.
\newblock \emph{Comm. Algebra}, 4\penalty0 (7):\penalty0 665--667, 1976.

\bibitem{Fuhrmann_Pym_2007}
C.~F\"{u}hrmann and D.~Pym.
\newblock On categorical models of classical logic and the geometry of interaction.
\newblock \emph{Math. Structures Comput. Sci.}, 17\penalty0 (5):\penalty0 957--1027, 2007.

\bibitem{Girard_1987}
J.-Y. Girard.
\newblock Linear logic.
\newblock \emph{Theoret. Comput. Sci.}, 50\penalty0 (1):\penalty0 101, 1987.

\bibitem{Garner_LopezFranco_2016}
R.~Garner and I.~L\'opez Franco.
\newblock Commutativity.
\newblock \emph{J. Pure Appl. Algebra}, 220\penalty0 (5):\penalty0 1707--1751, 2016.

\bibitem{Grandis_2020}
M.~Grandis.
\newblock \emph{Higher dimensional categories}.
\newblock World Scientific Publishing Co. Pte. Ltd., Hackensack, NJ, 2020.
\newblock From double to multiple categories.

\bibitem{Guglielmi_2007}
A.~Guglielmi.
\newblock A system of interaction and structure.
\newblock \emph{ACM Trans. Comput. Log.}, 8\penalty0 (1):\penalty0 Art. 1, 64, 2007.

\bibitem{Heunen_Vicary_2019}
C.~Heunen and J.~Vicary.
\newblock \emph{Categories for quantum theory}, volume~28 of \emph{Oxford Graduate Texts in Mathematics}.
\newblock Oxford University Press, Oxford, 2019.
\newblock An introduction.

\bibitem{Joyal_Street_1993}
A.~Joyal and R.~Street.
\newblock Braided tensor categories.
\newblock \emph{Adv. Math.}, 102\penalty0 (1):\penalty0 20--78, 1993.

\bibitem{Kudzman-Blais_Lemay_2025}
R.~Kudzman-Blais and J.-S.~P. Lemay.
\newblock Cartesian linearly distributive categories: Revisited.
\newblock \emph{Applied Categorical Structures}, 34\penalty0 (5):\penalty0 47, 2026.

\bibitem{Lamarche_1995}
F.~Lamarche.
\newblock Generalizing coherence spaces and hypercoherences.
\newblock In \emph{Mathematical foundations of programming semantics ({N}ew {O}rleans, {LA}, 1995)}, volume~1 of \emph{Electron. Notes Theor. Comput. Sci.}, page~15. Elsevier Sci. B. V., Amsterdam, 1995.

\bibitem{Lamarche_2007}
F.~Lamarche.
\newblock Exploring the gap between linear and classical logic.
\newblock \emph{Theory Appl. Categ.}, 18\penalty0 (17):\penalty0 473--535, 2007.

\bibitem{Lambek_1972}
J.~Lambek.
\newblock Deductive systems and categories. {III}. {C}artesian closed categories, intuitionist propositional calculus, and combinatory logic.
\newblock In \emph{Toposes, algebraic geometry and logic ({C}onf., {D}alhousie {U}niv., {H}alifax, {N}.{S}., 1971)}, Lecture Notes in Math., Vol. 274, pages 57--82. Springer, Berlin-New York, 1972.

\bibitem{Lambek_Scott_1988}
J.~Lambek and P.~J. Scott.
\newblock \emph{Introduction to higher order categorical logic}, volume~7 of \emph{Cambridge Studies in Advanced Mathematics}.
\newblock Cambridge University Press, Cambridge, 1988.
\newblock Reprint of the 1986 original.

\bibitem{MacLane_1998}
S.~Mac~Lane.
\newblock \emph{Categories for the working mathematician}, volume~5 of \emph{Graduate Texts in Mathematics}.
\newblock Springer-Verlag, New York, second edition, 1998.

\bibitem{Naeimabadi_2024}
S.~Naeimabadi.
\newblock \emph{Cartesian Linear Bicategories}.
\newblock Ph.d. thesis, University of Ottawa, Ottawa, Canada, 2024.

\bibitem{Seely_1989}
R.~A.~G. Seely.
\newblock Linear logic, {$*$}-autonomous categories and cofree coalgebras.
\newblock In \emph{Categories in computer science and logic ({B}oulder, {CO}, 1987)}, volume~92 of \emph{Contemp. Math.}, pages 371--382. Amer. Math. Soc., Providence, RI, 1989.

\bibitem{Spivak_Srinivasan_2024}
D.I. Spivak and P.V. Srinivasan.
\newblock What kind of linearly distributive category do polynomial functors form?, 2024.
\newblock (preprint).

\bibitem{Strassburger_2002}
L.~Stra{\ss}burger.
\newblock A local system for linear logic.
\newblock In \emph{Logic for programming, artificial intelligence, and reasoning}, volume 2514 of \emph{Lecture Notes in Comput. Sci.}, pages 388--402. Springer, Berlin, 2002.

\bibitem{Strassburger_2007_1}
L.~Stra{\ss}burger.
\newblock On the axiomatisation of {B}oolean categories with and without medial.
\newblock \emph{Theory Appl. Categ.}, 18\penalty0 (18):\penalty0 536--601, 2007a.

\bibitem{Strassburger_2007_2}
L.~Stra{\ss}burger.
\newblock A characterization of medial as rewriting rule.
\newblock In Franz Baader, editor, \emph{Term Rewriting and Applications}, pages 344--358, Berlin, Heidelberg, 2007b. Springer Berlin Heidelberg.

\bibitem{Tiu_2006}
A.~Tiu.
\newblock A local system for intuitionistic logic.
\newblock In \emph{Logic for programming, artificial intelligence, and reasoning}, volume 4246 of \emph{Lecture Notes in Comput. Sci.}, pages 242--256. Springer, Berlin, 2006.

\end{references*}

\end{document}